\documentclass[10pt]{amsart}

\usepackage[all]{xypic}

\usepackage{latexsym}

\usepackage{amssymb}
\usepackage{amsfonts}
\usepackage{amscd}
\usepackage{amsmath,amsthm}
\usepackage{verbatim}

\usepackage{tikz, tikz-cd}
\usetikzlibrary{matrix,arrows,decorations.pathmorphing, fit}
\tikzset{myboxgroup/.style={draw, densely dotted}} % style for the boxed groups

\newtheorem{lemma}{Lemma}[section]
\newtheorem{proposition}[lemma]{Proposition}
\newtheorem{theorem}[lemma]{Theorem}
\newtheorem{corollary}[lemma]{Corollary}

\newtheorem{lem}[lemma]{Lemma}
\newtheorem{prop}[lemma]{Proposition}
\newtheorem{thm}[lemma]{Theorem}
\newtheorem{cor}[lemma]{Corollary}

{
}

\theoremstyle{definition}

\theoremstyle{remark}

\numberwithin{equation}{section}

\newenvironment{pf}{\noindent{\bf Proof.}}{\hfill $\square$\medskip}

\def\CC{{\mathbb C}}

\def\GG{{\mathbb G}}

\def\NN{{\mathbb N}}

\def\PP{{\mathbb P}}

\def\ZZ{{\mathbb Z}}

\def\0ol{{\bar 0}}
\def\1ol{{\bar 1}}
\def\2ol{{\bar 2}}
\def\ol2{{\bar 2}}
\def\3ol{{\bar 3}}
\def\4ol{{\bar 4}}
\def\5ol{{\bar 5}}
\def\6ol{{\bar 6}}
\def\7ol{{\bar 7}}
\def\8ol{{\bar 8}}
\def\9ol{{\bar 9}}

\def\bold0{{\bf 0}}
\def\bold1{{\bf 1}}
\def\bold2{{\bf 2}} 
\def\bold3{{\bf  3}}
\def\bold4{{\bf 4}}
\def\bold5{{\bf 5}}
\def\bold6{{\bf 6}}
\def\bold7{{\bf 7}}
\def\bold8{{\bf 8}}
\def\bold9{{\bf 9}}

\def\P2Skly{\PP^2_{Skly}}

\def\coker{\operatorname {coker}}

\def\End{\operatorname {End}}
\def\Ext{\operatorname {Ext}}

\def\gr{\operatorname {gr}}

\def\ker{\operatorname {ker}}

\def\Ann{\operatorname{Ann}}
\def\Aut{\operatorname{Aut}}

\def\Br{\operatorname{Br}}

\def\det{\operatorname{det}}

\def\dim{\operatorname{dim}}

\def\End{\operatorname{End}}

\def\Ext{\operatorname{Ext}}
\def\Fdim{{\sf Fdim}}
\def\fdim{{\sf fdim}}

\def\Fract{\operatorname{Fract}}

\def\GKdim{\operatorname{GKdim}}

\def\gr{{\sf gr}}

\def\Gr{{\sf Gr}}

\def\id{\operatorname{id}}

\def\Im{\operatorname{Im}}

\def\mod{{\sf mod}}
\def\Mod{{\sf Mod}}

\def\\pdim{\operatorname{\pdim}}

\def\pdim{{\rm pdim}}

\def\Proj{\operatorname{Proj}}
\def\Projnc{\operatorname{Proj}_{nc}}

\def\QGr{\operatorname{\sf QGr}}

\def\Spec{\operatorname{Spec}}

\def\sup{\operatorname{sup}}
\def\Supp{\operatorname{Supp}}

\def\ul1{\operatorname{\underline{1}}}

\def\G{\mathop{\underline{\underline{\it \Gamma}}}\nolimits}

\def\l{\leftarrow}

\def\a{\alpha}
\def\b{\beta}
\def\c{\gamma}

\def\g{\gamma}
\def\l{\lambda}
\def\o{\omega}

\def\t{\tau}

\def\G{\Gamma}
\def\L{\Lambda}
\def\O{\Omega}

\def\fa{{\mathfrak a}}

\def\fgl{{\mathfrak g}{\mathfrak l}}

\def\fm{{\mathfrak m}}

\def\fp{{\mathfrak p}}
\def\fq{{\mathfrak q}}

\def\fsl{{\mathfrak s}{\mathfrak l}}

\def\wtS{{\widetilde{ S}}}

\def\cal{\mathcal}

\def\cA{{\cal A}}

\def\cL{{\cal L}}

\def\cO{{\cal O}}

\def\cV{{\cal V}}

\def\cZ{{\cal Z}}

\def\Qcoh{{\sf Qcoh}}

\def\dirlim{\mathop{\vtop{\baselineskip -100pt\lineskip -1pt\lineskiplimit 0pt
\setbox0\hbox{lim}\copy0\hbox to \wd0{\rightarrowfill}}}\limits}
\def\invlim{\mathop{\vtop{\baselineskip -100pt\lineskip -1pt\lineskiplimit 0pt
\setbox0\hbox{lim}\copy0\hbox to \wd0{\leftarrowfill}}}\limits}

\def\I11{{1 \kern -0.8pt \! \mbox{l}}}
\def\mumu{{\mu\kern-4.2pt\mu}}
\def\bfmu{{\mu\kern-4.2pt\mu}}
\def\2slash{\backslash \! \backslash}

\def\boxtimes{\setbox0\hbox{$\Box$}\copy0\kern-\wd0\hbox{$\times$}}

\begin{document}

\title[ The 4-dimensional Sklyanin Algebras at points of finite order]{Simple modules over the 
4-dimensional Sklyanin Algebras at points of finite order}

\author{S. Paul Smith }
\address{Department of Mathematics, Box 354350, Univ.  Washington, Seattle, WA 98195}
\email{smith@math.washington.edu}



\begin{abstract}
 In 1982 E.K. Sklyanin defined a family of graded algebras $A(E,\tau)$, 
depending on an elliptic curve $E$ and a point $\tau \in E$ that is not 4-torsion. 
The present paper is concerned with the structure of $A$ when $\tau$ is a point of 
finite order, $n$ say. It is proved that every simple $A$-module has 
dimension $\le n$ and that ``almost all'' have dimension precisely
$n$. There are enough finite dimensional simple modules to separate elements of $A$; 
that is, if $0\ne a \in A$, then there exists a 
simple module $S$ such that $a.S \ne 0.$ Consequently $A$
satisfies a polynomial identity of degree $2n$ (and none of lower
degree). Combined with results of Levasseur and Stafford it follows that 
$A$ is a finite module over its center. Therefore one may
associate to $A$ a coherent sheaf, ${\cal A}$ say,  of finite ${\cal O}_S$ 
algebras where $S$ is the projective 3-fold determined by the center of 
$A$. We determine where ${\cal A}$ is Azumaya, and prove that
the division algebra $\Fract({\cal A})$ has rational center. Thus, for each $E$
and each $\t \in E$ of order $n \ne 0,2,4$ one obtains a division algebra of
degree $s$ over the rational function field of $\PP^3$, where $s=n$ if $n$ is
odd, and $s={{1} \over {2}} n$ if $n$ is even.

The main technical tool in the paper is the notion of a ``fat point''
introduced by M. Artin. A key preliminary result is the classification of
the fat points: these are parametrized by a rational 3-fold.
\end{abstract}

\maketitle

\section{Introduction}

\subsection{}

Throughout this paper $A(E,\tau)$ denotes the 4-dimensional Sklyanin associated to a complex elliptic curve $E$ and a 
translation automorphism $\tau:E \to E$ that is not 4-torsion. These algebras were introduced by E.K. Sklyanin  in 1982  
(\cite{Skly11}, \cite{Skly12}) in order to study certain exactly solvable models of quantum and statistical physics. 

They form a flat family of noetherian algebras and are, in some sense, the most generic 
non-commutative deformations of the polynomial ring in 4 variables.  
This paper builds on earlier work by Artin, Levasseur, Smith, Stafford, Staniszkis, Tate, and Van den Bergh.

%of results in the papers  \cite{Skly11}, \cite{Skly12}, \cite{SS15}, \cite{LS8}, and \cite{Stan94}. 
%We assume the reader has some familiarity with those papers.

 The present paper examines their representation theory when $\t$ has finite order: our main result is a 
 classification of their finite dimensional simple modules. When $\tau$ has infinite order their finite dimensional simple modules 
 were classified in \cite{SS16}.   There is a sharp dichotomy between the finite and infinite order cases. 
 In particular,  $A(E,\tau)$ is a finitely generated module over its center if and only if $\tau$ has finite order.

\begin{theorem}
(Theorem \ref{thm.7.8})
%(Theorem \ref{thm3.8})
\label{thm.main1}
If $\tau$ has order  $n<\infty$, then $A(E,\tau)$  
\begin{enumerate}
\item{} is a finitely generated module over its center,
\item{} is a maximal order in its division ring of fractions, and
\item{} satisfies the same polynomial identities as the ring of $n \times n$ matrices.
\end{enumerate}
\end{theorem}

When $\tau$ has finite order we define
\begin{align*}
n & \; =\; \text{the order of $\tau$ and}
\\
s & \; =\; \text{the order of $2\tau$.}
\end{align*}

Suppose $\tau$ has finite order and let $Z$ denote the center of $A(E,\tau)$. 
The center is a finitely generated algebra and $\Spec(Z)$ is an irreducible variety of dimension 4. 
Theorem \ref{thm6.7} provides an explicit element 
$c \in Z$ such that $A[c^{-1}]$ is an Azumaya algebra of rank $n^2$ over the open set $\{c \ne 0\} \subseteq \Spec(Z)$.
To prove this one must first show that ``almost all'' simple $A(E,\tau)$-modules have dimension $n$.

%The simple modules of dimension $<n$ are particularly interesting. There 
 
\subsection{} 
\label{sect.cA}
 Suppose $\tau$ has finite order. Because $A(E,\tau)$ is a graded algebra so is its center, $Z$ say.  
The projective variety $S=\Proj(Z)$ is irreducible of dimension 3. 
Associated to $A(E,\tau)$ is a sheaf $\cA$ of coherent $\cO_S$-algebras: if $z$ is a homogeneous element in $Z$, then the sections of $\cA$ over the open set $\{z \ne 0\}$ is   the degree-zero component of
 the localization $A(E,\tau)[z^{-1}]$. In many ways, $\cA$ is of more interest than $A(E,\tau)$ and many of the results in this
 paper are, in effect, results about $\cA$.  For example, the fact that $A(E,\tau)$ is a maximal order implies that 
 $\cA$ is a sheaf of maximal orders. Over its Azumaya locus, $\cA$ is a sheaf of Azumaya algebras of rank $s^2$.\footnote{We haven't been  entirely honest here.  When $n$ is odd, the sheaf of centers of $\cA$ is $\cO_S$ but when $n$ is even the sheaf of centers is
 strictly larger than $\cO_S$  so one should in that case consider $\cA$ as a sheaf of algebras over the 
 projective variety ${\bf Spec}(\cZ)$ where $\cZ$ is that sheaf of centers. The variety ${\bf Spec}(\cZ)$ is $\Proj(Z(A^{(2)})$ where $Z(A^{(2)}$ denotes the center of the even-degree subalgebra of $A$. This issue is examined in \cite{ST94}.}

Much more is known about graded $A$-modules than modules are not graded, and most of the results in this paper require, first ,
a good understanding of graded $A$-modules.

\subsection{} 
It was shown in \cite{SS15} that the homogeneous components of $A(E,\tau)$ have the same dimension as the homogenous 
components of the polynomial ring on 4 variables. It was also shown there that $A(E,\tau)$ is a noetherian domain 
(i.e., it has no zero-divisors) having  the same homological properties as the polynomial ring on 4 variables.  There are similar deformations of the polynomial rings on $n$ variables for all $n \ge 3$ (\cite{OF}, \cite{TV}). Those other deformations, which were discovered/invented by Artin-Schelter, Cherednik, and Odesskii-Feigin (\cite{AS},  \cite{Ch}, \cite{OF}), also come in families parametrized by the data $(E,\tau)$.

Apart from their interest as algebraic objects, Sklyanin algebras also provide rich examples of non-commutative algebraic geometry.
In particular, graded deformations of the polynomial ring on $n$ variables play the role of homogeneous coordinate rings for non-commutative analogues of the projective space $\PP^{n-1}$.  
%The 4-dimensional Sklyanin algebras also appear naturally in the work of Connes and Dubois-Violette on 
%non-commutative 3-spheres {\bf refs} .

\subsection{}
When $\t$ has infinite order,  $A(E,\tau)$ resembles the enveloping algebra of the Lie algebra $\fgl_2$.
It also resembles  a homogenized version of the quantized enveloping algebra $U_q(\fsl_2)$ of $\fsl_2$. 
Indeed, the homogenized  $U_q(\fsl_2)$'s are ``degenerate'' Sklyanin algebras with the 
elliptic curve degenerating to a pair of conics in $\PP^3$; and $\tau$ degenerates to the parameter $q$. 
These degenerate Sklyanin algebras and their relation with $U_q(\fsl_2)$ are the subject of  \cite{CSW};
the degeneration proceedure is discussed in detail in  \cite[\S6]{CSW}. 

There are interesting similarities with the classical representation theory of simple Lie algebras, though it is the differences that are most
fascinating. For example, when $\tau$ has infinite order the finite dimensional representations of $A(E,\tau)$ are realized through
difference operators acting on spaces of theta functions which are, essentially, the same things as sections of line bundles on 
the elliptic curve $E$, and this provides a nice parallel, and contrast,  with the Borel-Weil-Bott theorem 
that realizes the irreducible representations of 
the Lie algebra as differential operators on sections of line bundles on the flag variety which is, in the case of $\fsl(2)$, the projective line
$\PP^1$.

\subsection{}
In \cite{Skly12}, Sklyanin described a set finite dimensional representations of $A(E,\tau)$.
In \cite{SS16}, it was shown that these are simple when $\tau$ has infinite order. Using the fact that over a 
connected graded $\CC$-algebra non-trivial simple modules come in 1-parameter families parametrized by $\CC^\times$, 
Sklyanin's simple modules give rise to all finite dimensional simple $A(E,\tau)$-modules when $\tau$ has infinite order. 
The results in \cite{SS16} make use of various infinite dimensional graded $A(E,\tau)$-modules in an essential way, particularly
point modules and line modules (see \S\ref{sect.linear.modules} for the definitions). 
These graded modules, which were studied in great detail in  \cite{LS8}, also play an essential role in this paper. 
A dominant theme in all these results is that the representation theory of $A(E,\tau)$ is 
intimately related to, indeed, controlled by,  the projective geometry of $E$ embedded as a quartic curve in  $\PP^3$  
and the interaction of that geometry with the translation automorphism $p \mapsto p+\tau$.   

By \cite{SS16}, the point modules (more precisely, their isomorphism classes) are parametrized by the points on $E$ and the 
4 points that are the singular loci of the four singular quadrics in $\PP^3$  that contain $E$.  By \cite{LS8}, the line modules are in 
natural bijection with the secant lines to $E \subseteq \PP^3$.  We write $M(p)$ for the point module corresponding to 
a point $p$ and $M(p,q)$ or $M(\ell)$ for the line module corresponding to the secant line $\ell$ 
whose scheme-theoretic intersection with $E$ is the divisor $(p)+(q)$.

\subsection{The strategy}
\label{sect.strategy}
For brevity, we write $A=A(E,\tau)$.

Most of the results in this paper rely on an understanding of the 1-critical $A$-modules, also called (fat) point modules,
these being the graded $A$-modules of Gelfand-Kirillov dimension 1 all of whose proper quotient modules have finite dimension. 
Every finite dimensional simple $A$-module is the quotient of such a module. 
Accordingly, \S\ref{sect.3} is devoted to some results for a large class of algebras 
that relate finite dimensional simple modules to 1-critical modules. Some of these results 
These preliminary results are applied to the $A(E,\tau)$ in later sections.

Similar problems to those addressed in this paper are solved in \cite{ATV2} and \cite{A1} for the 3-dimensional Sklyanin algebras; actually 
the results there apply to a larger class of algebras, namely the 3-dimensional Artin-Schelter regular algebras.
The algebras $A(E,\tau)$ are 4-dimensional Artin-Schelter regular algebras (\cite{AS}, \cite{SS15}). 
We make use of several results in  \cite{ATV2} and \cite{A1}.

We first show that $A$ satisfies a polynomial identity. This done by showing, in Proposition \ref{lem3.5}, that if $M(\ell)$ is a
line module, then $A/ \! \Ann(M(\ell))$ satisfies the polynomial identities of $2n \times 2n$ matrices, and also showing that the diagonal map
$$
A \; \longrightarrow \; \bigoplus_{\text{all secant lines } \ell} \; \frac{A}{ \! \Ann(M(\ell))}
$$
is injective; it follows immediately from this that $A$ also satisfies the polynomial identities of $2n \times 2n$ matrices; this is a 
weak version of Theorem \ref{thm.main1}(3). Once we know that $A$ satisfies a polynomial identity, a 
result of Stafford \cite{JTS17} tells us that $A$ is a maximal order. It is therefore equal to its own trace ring;
trace rings are finitely generated modules over their centers.

Once we know that $A(E,\tau)$ is finite over its center, we turn to the classification of its  finite dimensional simple modules. 
After  the results in  \S\ref{sect.3}  this is more or less equivalent to determining the equivalence classes of
1-critical graded modules where two such modules are said to be equivalent if they have isomorphic non-zero submodules 
(morphisms between graded modules are required to preserve degree). The equivalence classes are called {\it fat points}; 
actually this is not quite correct, but suffices for the introduction (fat points are those of
multiplicity $>1$, and points are those of multiplicity $=1$). Thus the problem becomes that of determining all the fat points for 
$A$. This problem is addressed in \S\S5-8.

In \S \ref{sect.fat.pts} we consider fat points and fat point modules and show amongst other things that
every fat point lies on a secant line. By definition, this means that the fat point has a representative that is a quotient of a line module. 
The 1-critical modules of multiplicity 1 are the point modules. 
There are exactly 4 fat points of each multiplicity $2,3,\ldots, s-1$ and all other fat points have multiplicity $s$
 (it is not until \S7 that this proof is completed). 
One consequence of the fact that most fat points are of multiplicity $s$ is that $A$ satisfies the polynomial identities of 
$n \times n$ matrices  (and does not satisfy the identities of any smaller matrix algebra). 

The classification of the fat points of multiplicity $s$ allows one to determine the center of the division algebra
$\Fract(\cA)$ where $\cA$ is the sheaf of algebras defined in \S\ref{sect.cA}. Fat points are essentially the same things as the simple
$\cA$-modules, i.e., the simple objects in the category $\Qcoh(\cA)$ of quasi-coherent $\cA$-modules. 
In fact, fat points are the simple objects in the quotient category $\QGr(A)$ of the category of graded left $A$-modules ($\QGr(A)$ 
is defined in \S\ref{sect.cats}).

\subsection{Fat points of intermediate multiplicity}
The fat points of multiplicity other than $1$ or $s$, are said to have 
{\it intermediate multiplicity}. These are classified in \S6. There are exactly 4 fat points of intermediate multiplicity $k+1$ 
for $0\le k<s-1$ and these are labelled $F(\o+k\t)$ where $\o \in E[2]$, the
$2$-torsion subgroup of $E$.  This labelling is such that the
$F(\o)$ are actually the point modules $M(e_i)$ corresponding to the vertices of the singular quadrics  containing $E$. 
These 4 point modules have more in common with the fat points of intermediate multiplicity than with
the point modules $M(p)$ for $p\in E$. 
The labelling is such that $F(\o+k\t)$ lies on all the secant lines $\ell_{pq}$ such that $p+q=\o+k\t$; this is expressed more 
precisely through the exact sequences 
$$
0 \to M(p-(k+1)\t,q-(k+1)\t)(-k-1) \to M(p,q) \to F(\o+k\t) \to 0.
$$

\subsection{Multiplicity-$s$ fat points}
Section 7 begins the classification of the fat points having multiplicity $s$. 
They are classified in terms of the lines on which they lie. The
first step is to show that every fat point does lie on a line (this is done
in \S\ref{sect.fat.pts}). Much of Section 7 is rather technical: we must analyse exactly 
which lines a fat point can lie on, what fat points a given line can contain, 
and when two lines can contain a common fat point. As a culmination of these
technicalities, we are able to prove that a fat point has multiplicity at
most $s$, and that every secant line contains 
infinitely many multiplicity-$s$ fat points; in fact the multiplicity-$s$ fat points  that lie on the secant line $\ell_{pq}$  arise as the
cokernels of (most of) the maps in  ${\rm Hom}(M(p-s\t,q-s\t)(-s),M(p,q))$. It is
shown that this Hom-space is 2-dimensional, whence the multiplicity-$s$ fat points on a given secant line are 
naturally parametrized by $\PP^1$ 
minus some finite set. As a consequence of all this it follows that there is
a central element $c$, specified in Theorem \ref{thm6.7}, such that $A[c^{-1}]_0$ is
Azumaya of rank $s^2$ over its center.

It can happen that two distinct secant lines
contain infinitely many fat points in common. In this case we say
that the lines are {\it equivalent}. A secant line is
equivalent to a finite number of other lines, and Theorem \ref{thm7.6} gives a
precise description of the equivalence classes. It is nearly true that 
equivalence
classes of lines are in bijection with the points on $S^2(E')$ the second
symmetric power of $E'= E/\langle 2\t \rangle.$ More precisely, the
equivalence classes of secant lines are nearly in bijection with the secant
lines to $E'$ (for an appropriate embedding of $E'$ in another $\PP^3$).
One of the main algebraic consequences of such considerations is that $A$ is
a polynomial identity ring of p.i.degree $n$; that is, $A$ has a faithful
family of $n$-dimensional simple modules.

In general, a multiplicity-$s$ fat point  lies
on exactly two equivalence classes of lines; compare this to the geometric
statement that a general point in $\PP^3\, - \, E$ lies on exactly two 
secant lines of $E$. The exceptions involve the points that lie on the
singular quadrics containing $E$; if $p$ is a smooth point of such a
singular quadric then $p$ lies on a unique secant line, and if $p$ is
a singular point then $p$ lies on infinitely many secant lines. There
are analogous exceptions involving fat points. One such exception has
already been mentioned: $F(\o+k\t)$ lies on
infinitely many secant lines (recall that $F(\o)$ is the point module
associated to one of the singular points). Hence if one tries to have
a geometric picture of the point modules and fat points, then the fat
points of intermediate multiplicity play the role of singular points
of various quadrics. The other exception involves the multiplicity-$s$ fat points 
that lie on a line that contains one of the $F(\o+k\t)$. These lie
on only one equivalence class of lines.

The analysis of multiplicity-$s$ fat points splits into two cases depending
on whether or not the fat point is annihilated by a certain special central
element $c$ (see Theorem \ref{thm6.7}). The easiest fat points to deal with are those not annihilated by $g$. This is closely related to the fact that $A[c^{-1}]_0$ is an Azumaya algebra of rank $s^2$ over its center. 
By Theorem \ref{thm7.10},  such fat points are parametrised by a dense open subset of
the variety $E' \times E' \times E$ modulo the action of a certain
finite group. Theorem \ref{thm7.11} shows that this variety is rational. It follows
that the division algebra $D=\Fract(\cA)$ has center the rational function field
of $\PP^3$ (the fact that most fat points have multiplicity $s$ implies that
the degree of $D$ is $s$). Thus each Sklyanin algebra $A(E,\t)$ gives rise
to such a division algebra. It is probably an interesting project to
understand these division algebras in detail.

\subsection{}
An early draft of this paper was circulated in the early 1990's. Shortly after that an explicit description of the center of $A(E,\tau)$ was given in \cite{ST94} in the case $\tau$ has finite order. Some of the results in the earlier draft appeared in \cite[\S 11]{SPS-survey} and 
\cite[\S3]{LB}, and some were used in \cite{ST94}.  

\smallskip

{\bf Acknowledgements.} This work was supported by NSF grant DMS-9100316.
I am very grateful to Mike Artin,  Lieven LeBruyn, Thierry Levasseur, 
Toby Stafford, Joanna Staniszkis, John Tate, Michaela Vancliff, Michel Van den Bergh, 
and James Zhang for useful discussions and feedback on the earlier draft.

%%%%%%%%%%%%%%%%%%%%%%%%%%%%%%%%%%%%%%%%%%%%%%%%%%%%%%%%%%%%%%%%%
%%%%%%%%%%%%%%%%%%%%%%%%%%%%%%%%%%%%%%%%%%%%%%%%%%%%%%%%%%%%%%%%%
%%%%%%%%%%%%%%%%%%%%%%%%%%%%%%%%%%%%%%%%%%%%%%%%%%%%%%%%%%%%%%%%%
  
    \section{Preliminaries}
    \label{sect.prelims}

\subsection{Definition}
Once and for all we fix a lattice $\L \subseteq \CC$. Throughout this paper 
$E=(E,+,0)$ denotes the abelian group $\CC/\L$. We always think of $E$ as an
algebraic variety, i.e., as an elliptic curve.
We write $E[2]$ for its 2-torsion subgroup and $E[4]$ for its 4-torsion subgroup.

Following the notation and definitions in \cite{Web19}, let $\{\theta_{ab} \;|\; ab=00,01,10,11\}$ denote Jacobi's four theta functions 
with period lattice $\L$. The function $\theta_{11}(z)$ is an odd function and the other three theta functions are even functions.

The 4-dimensional Sklyanin algebra is the free algebra $\CC\langle x_0,x_1,x_2,x_3\rangle$ modulo the 6 relations
$$
[x_0,x_\alpha] = iJ_{\beta \gamma} \{x_\beta,x_\gamma\}
\qquad \text{and} \qquad
[x_\alpha,x_\beta]  \; =\;   i\{x_0,x_\gamma\},
$$
where $(\a,\b,\gamma)$ runs over the cyclic permutations of $(1,2,3)$,
$[x,y]=xy-yx$, $\{x,y\}=xy+yx$, and 
$$
J_{12} = 
{{
\theta_{11}(\tau)^2
\theta_{01}(\tau)^2
}\over{
\theta_{00}(\tau)^2
\theta_{10}(\tau)^2
}},
\qquad
J_{23} = 
{{
\theta_{11}(\tau)^2
\theta_{10}(\tau)^2
}\over{
\theta_{00}(\tau)^2
\theta_{01}(\tau)^2
}},
\qquad
J_{31} = -
{{
\theta_{11}(\tau)^2
\theta_{00}(\tau)^2
}\over{
\theta_{01}(\tau)^2
\theta_{10}(\tau)^2
}}.
$$

\subsection{The quartic elliptic curve $E \subseteq \PP^3$} 
The divisor $4(0)$ on $E$ is very ample. 
Almost always,  we will consider $E$ as a quartic curve in $\PP^3=\PP^3_{x_0,\ldots,x_3}$ embedded via the complete linear 
system $|4(0)|$.
If we do this, four points on $E$, counted with multiplicity, lie on a common hyperplane in $\PP^3$ if and only if their sum is 0.

This embedding $i:E \to \PP^3_{x_0,x_1,x_2,x_3}$ 
can be described in terms of Jacobi's theta functions: thinking of $z$ as a point in $\CC$ modulo $\L$,
$$
i(z) \;=\; (g_{11}(z),g_{00}(z),g_{01}(z),g_{10}(z)) \in \PP^3_{x_0,\ldots,x_3}
$$
where 
$$
g_{ab}(z) \; = \;  \gamma_{ab}\theta_{ab}(2z)(\t)\theta_{ab}(\t)  \quad \hbox{ and } \quad
 \g_{ab}= \begin{cases}
\sqrt{-1}, & \text{if $ab=00,11$}
\\
1, &\text{if $ab=01,10$}
\end{cases}
$$
and $x_{11}=x_0$, $x_{00}=x_1$, $x_{01}=x_2$, and $x_{10}=x_3$.

If $p,q \in E$ we write $\overline{pq}$ or $\ell_{pq}$ for the secant line to $E \subseteq \PP^3$ whose scheme-theoretic intersection 
with $E$ is the divisor $(p)+(q)$.  For each $z \in E$ we define  
$$
L(z) \; :=\;  \{ \ell_{pq} \; | \; p+q=z \}
$$ 
and call it a {\sf family of lines}.   
The union of the lines in $L(z)$ is a quadric that we denote by $Q(z)$.
It obviously contains $E$. There is a pencil of quadrics containing $E$ and each of those quadrics is equal to
$Q(z)$ for some $z \in E$.  

If $\ell$ is a line in $L(z)$  and $\ell'$ a line $L(-z)$, then the sum of the points in $\ell \cap E$, counted with multiplicity,
and the points in $\ell' \cap E$, also counted with multiplicity, is 0 so $\ell$ and $\ell'$ belong to a common hyperplane and hence 
$\ell \cap \ell'  \ne \varnothing$.  It follows that  $Q(z)=Q(-z)$. These are the only equalities among the quadrics so the pencil of 
quadrics containing $E$ is parametrized by the projective line $E/\pm$.

The smooth $Q(z)$'s are isomorphic to $\PP^1 \times \PP^1$ and the lines in $L(z)$ and $L(-z)$ provide the two rulings on $Q(z)$.
If $Q(z)$ is not smooth it has a single ``ruling'' and a unique singular point. It is clear that $Q(z)$ has a single ruling exactly when
$L(z)$ and $L(-z)$. It follows from what we have already said that this happens if and only if $z \in E[2]$. 
 Since $E[2] \cong (\ZZ/2)^2$, there are exactly four  singular quadrics. 
 % Of course, there is a more elementary way to see this. If $q$ and $q'$ are two generic degree-$2$ elements in the polynomial ring on 4 
 % variables, then $\det(\l q+\l' q')$ is a quartic polynomial on the projective line $\PP^1_{\l,\l'}$ so has four zeroes; in fact, four distinct 
 % zeroes, and the rank of $\l q+\l' q'$ is 3 at those four points.
We denote by $e_0,\ldots,e_3$ the singular loci of these four quadrics. We label these in such a way that 
$e_i$ is the common zero-locus of $\{x_0,x_1,x_2,x_3\}-\{x_i\}$.

\subsection{Categories}
\label{sect.cats}
Some of the preliminary results in \S\ref{sect.3} apply to a wider class of algebras than just the $A(E,\tau)$'s.
With that in mind we now let $\Bbbk$ denote an arbitrary field and assume that $R$ is a left and right noetherian connected graded 
$\Bbbk$-algebra. $A(E,\tau)$ has these properties.
  
 We write $\Mod(R)$ for the category of left $R$-modules, and $\Gr(R)$ for the category of $\ZZ$-graded left $R$-modules with degree-preserving homomorphisms. We write $\mod(R)$ and $\gr(R)$ for the full subcategories of  $\Mod(R)$ and $\Gr(R)$ consisting of the noetherian modules. \footnote{All rings in this paper are both left and right noetherian.}
 
The category of graded left $R$-modules is denoted by $\Gr(R)$. Morphisms in $\Gr(R)$ are the $R$-module homomorphisms that preserve degree. We write $\gr(R)$ for the full subcategory of $\Gr(R)$ consisting of finitely generated $R$-modules.

Unless we say otherwise all modules in this paper are left modules.

If $M \in \Gr(R)$ and $n \in \ZZ$, we write $M(n)$ for the graded $R$-module that is $M$ as an $R$-module and degree-$i$ component 
$M(n)_i=M_{n+i}$. If $f:M \to M'$ is a morphism in $\Gr(R)$ we define $f(n):M(n) \to M(n)$ to be the homomorphism $f$.
We call the auto-equivalence $M \rightsquigarrow M(n)$ of $\Gr(R)$ a {\sf degree-shift} or {\sf Serre twist}.

We write  $\Fdim(R)$, resp.,  $\fdim(R)$, for the full subcategory of $\Gr(R)$,  resp.,  $\gr(R)$, 
consisting of those modules that are the sum of their  finite-dimensional submodules.
Since $\Fdim(R)$ is a Serre subcategory, i.e., closed under  sub-quotients and extensions,  there is a quotient category
$$
\QGr(R) \; :=\; \frac{\Gr(R)}{\Fdim(R)}.
$$
The objects in $\QGr(R)$ are the objects in $\Gr(R)$ but the morphisms differ. 
We write $\pi^*:\Gr(R) \to \QGr(R)$ for the quotient functor. 
Because $R$ is left noetherian $\Fdim(R)$ is closed under direct limits so $\pi^*$
has a right adjoint that we denote by $\pi_*$.

If $R$ is a quotient of a polynomial ring with its standard grading, then $\QGr(R)$ is equivalent to the category of quasi-coherent sheaves on $\Proj(R)$. For that reason we think of $\QGr(R)$ as playing the role of the category of  ``quasi-coherent sheaves'' on a
  ``non-commutative scheme'' that we denote by $\Projnc(R)$.  
  Many of the results in this paper are best thought of as results about the ``geometry'' of   $\Projnc(R)$.

% The 4-dimensional Sklyanin algebras have global homological dimension 4 and have he same Hilbert series as the polynomial ring on 4 variables. namely $(1-t)^{-4}$. Thus, every finitely generated graded left $A(E,\tau)$-module $M$ has a finite resolution by 
%finitely generated free left $A$-modules, whence $H(M;t)=f(t)(1-t)^{-4}$ for some $f(t) \in \ZZ[t^{\pm 1}]$. 

%The 4-dimensional Sklyanin algebras $A(E,\tau)$ have other good homological properties: they are Cohen-Macaulay and Gorenstein
%in the sense of Levasseur. See \cite{L7} and \cite{LS8} for details. We remind the reader that non-commutative algebras
%having finite global homological dimension need not have properties analogous to the usual commutative Cohen-Macaulay and Gorenstein properties. Those definitions need to be modified in suitable ways.

 \subsection{Equivalent and shift-equivalent graded modules}  
 \label{sect.equiv} 
Finitely generated graded $R$-modules $M$ and $N$ are {\sf equivalent} if they satisfy the following equivalent conditions:
\begin{enumerate}
  \item
 $M$ and $N$  are isomorphic as objects in $\QGr(R)$, 
  \item 
there are graded submodules $M^\prime$ and $N^\prime$ such that 
$M^\prime \cong N^\prime$ in $\Gr(A)$ and   $M/M^\prime$ and $N/N^\prime$ are finite dimensional;
  \item 
  $M_{\ge n} \cong N_{\ge n}$ for some $n \ge 0$.
\end{enumerate} 
If $M$ and $N$ are equivalent we write  $M \sim N$. This {\it is} an equivalence relation.
 
More generally, $M$ and $N$ are {\sf shift-equivalent} if $M(n) \sim N$ for some $n \in \ZZ$.

\subsection{Point modules and line modules}  
\label{sect.linear.modules}
Let $M \in \Gr(A)$. We call   $M$ a 
\begin{itemize}
  \item 
  {\sf point module} if it is cyclic and  $H(M;t)=(1-t)^{-1}$, and a
  \item 
 {\sf line module} if it is cyclic and $H(M;t)= (1-t)^{-2}$.
\end{itemize} 
By \cite{LS8}, $M$ is a 
\begin{itemize}
  \item 
point module if and only if $M \cong A/Ap^\perp$ for some point $p \in E \cup \{e_0,\ldots,e_3\}$; 
  \item 
 line module if and only if $M \cong A/A\ell^\perp$ for some secant line $\ell$.
 \end{itemize}
 If $p \in E \cup \{e_0,\ldots,e_3\}$ we write $M(p)$ for $A/Ap^\perp$. 
 If $\ell$ is a secant line to $E$ we write $M(\ell)$ for $A/A\ell^\perp$.
If the scheme-theoretic intersection of $\ell$ with $E$ is the divisor $(p)+(q)$ we also write $M(p,q)$ for $A/A\ell^\perp$. 

Point modules $M(p)$ and $M(p')$ are isomorphic if and only if $p=p'$. Line modules $M(\ell)$ and $M(\ell')$ are isomorphic if and only if 
$\ell=\ell'$.

The quotient of $A(E,\tau)$ by the ideal generated by $\{x_0,x_1,x_2,x_3\}-\{x_i\}$ is isomorphic to the polynomial ring $\CC[x_i]$.   
Thus, the point modules $M(e_i)$ corresponding to the singular loci of the four singular quadrics are these four quotient rings of 
$A(E,\tau)$. 

\subsubsection{Notation}
We write $\fa_{pq}$ and/or $\fa_{\ell}$ for the annihilator of the line module $M(\ell)=M(p,q)$.
 
\subsection{Degree-2 central elements in $A(E,\tau)$}

Sklyanin discovered two linearly independent central elements $\O_1$ and $\O_2$ in $A(E,\tau)_2$. 
It was shown in \cite{SS15} that the subalgebra $\CC[\O_1,\O_2]$ is a polynomial ring in two variables 
and that the quotient $A/\Omega_1,\O_2)$ is a twisted homogenous coordinate ring of $E$ 
in the sense of Artin and Van den Bergh \cite{AV}. More precisely, if $\cL=\cO_E(4(0))$, i.e., if $\cL=i^*\cO_{\PP^3}(1)$ where
$i:E \to \PP^3$ is the embedding above, then 
$$
\frac{A(E,\tau)}{(\Omega_1,\O_2)} \; \cong \; B(E,\cL,\tau)
$$
in the notation of \cite{AV}. By \cite{AV}, there is an equivalence of categories
\begin{equation}
\label{eq.QGr.B}
\QGr(B(E,\cL,\tau)) \;  \equiv \; \Qcoh(E).
\end{equation}
This equivalence plays a key role in understanding $A(E,\tau)$.

The homomorphism $A(E,\tau) \to B(E,\cL,\tau)$ leads to functors $i^*:\QGr(A(E,\tau))\to \Qcoh(E)$ and 
 $i_*:\Qcoh(E) \to \QGr(A(E,\tau))$ that play the role of inverse and direct image functors associated to a morphism 
 $i:E \to  \Projnc(A)$. There is a notion of closed subscheme for non-commutative algebraic geometry and $E$  is a closed 
 subscheme of $ \Projnc(A)$ in this sense \cite{SPS2004}.

The irreducible objects in $\Qcoh(E)$ are the skyscraper sheaves $\cO_p$ at the points $p \in E$.
The corresponding objects in $\QGr(B(E,\cL,\tau))$ are the point modules $M(p)$, $p \in E$. 
Indeed, up to isomorphism, the point modules $M(p)$, $p \in E$, are the only finitely generated 
1-critical $B(E,\cL,\tau)$-modules (the notion of a 1-critical module is introduced in \S\ref{sect.gkdim.1-crit}).

For the purposes of understanding $A(E,\tau)$, the projective line $\Proj(\CC[\O_1,\O_2])$ 
is best thought of as the quotient of $E$ by the group generated by the involution $z \mapsto -z-2\tau$. 
With this in mind, the non-zero elements in $\CC\O_1+\CC\O_2$, up to scalar multiples, are labelled $\O(z)$, $z \in E$,
with $\Omega(z)=\Omega(-z-2\tau)$. This can be done in such a way that  the line module $M(p,q)$ is annihilated by $\Omega(p+q)$
\cite[Cor. 6.6]{LS8}. By \cite[Prop. 6.2]{LS8}, $A/( \Omega(z) )$ is a domain of GK-dimension 3.

For each $z \in E$ we call the set of line modules $\{M(p,q) \; | \; p+q=z\}$ a {\it family of line modules.}  

The parallel between the elements $\Omega(z)$ and the quadrics $Q(z)$ should be apparent. 
The results in \cite{SVdB} show that the closed subschemes (in the non-commutative sense) 
$\Projnc\!\big(A(E,\tau)/(\Omega(z))\big)$ of $\Projnc(A(E,\tau))$ really do behave like quadrics in $\PP^3$. 
There is one difference that will appear often: the lines on $Q(z)$ are the lines $\ell_{pq}$ such that $p+q=\pm z$ whereas the line modules annihilated by $\Omega(z)$, which play the role of ``structure sheaves'' of ``noncommutative lines'' on 
$\Projnc\!\big(A(E,\tau)/(\Omega(z))\big)$, are $M(p,q)$ with $p+q=z$ and $p+q=-z-2\tau$.  

Another analogy that will play a role later on is this. The four singular quadrics are the $Q(\omega)$'s for $\omega \in E[2]$. The points
in $E[2]$ are the fixed points under the involution $z \mapsto -z$.
There are four central elements  that annihilate only one family of line modules, namely $\O(\o-\t)$ for $\o \in E[2].$  
The points $\o-\t$ for $\o \in E[2]$ are the fixed points under the involution $z \mapsto -z-2\tau$. 
The main result in \cite{SVdB} shows that exactly four  of the non-commutative  quadrics $\Projnc(A/(\Omega(z))$ are singular, 
namely those for which $z=\o-\t$ for  some $\o \in E[2]$; $\Projnc(A/(\Omega(z))$ is defined to be singular if $\Ext^i \ne 0$ on
$\QGr((A/(\Omega(z)))$ for $i \ge 3$; equivalently, it is singular if there are objects in $\QGr((A/(\Omega(z)))$ that do not have a 
finite injective resolution.

 \subsection{A basis for certain line modules}   
 \label{sect.good.basis}
 We frequently assume that $p-q \notin 2\ZZ\tau$. Under this assumption  $M(p,q)$ has a basis $\{e_{ij} \; | \;i,j \ge 0\}$
 such that $e_{ij} \in M(p,q)_{i+j}$ and 
$$
Ae_{ij} \cong M(p+(j-i)\tau, q+(i-j)\tau).
$$ 
See  \cite[5.6]{LS8}.  We will need an explicit description of the action of $A$ on such a basis. For this we use \cite[Prop. 3.1]{Stan94}. 

To do this we must $p$ and $q$ by points in $\CC-(2\ZZ\tau + \Lambda)$ at which we can  evaluate  
 theta functions. Accordingly, fix $p,q \in \CC$ such that  $p-q \notin 2\ZZ\tau + \Lambda.$ By  \cite[Prop. 3.1]{Stan94},
$M(p,q)$ has a basis $\{e_{ij} \; | \;i,j \ge 0\}$, which we fix for the remainder of the paper, on which $X\in A_1$ acts by 
$$
X.e_{ij} \quad = \quad
 {{X(q+(i-j)\tau)} \over {\theta_{11}(p-q-2(i-j)\tau)}}e_{i,j+1} \quad
- \quad  {{X(p+(j-i)\tau)} \over {\theta_{11}(p-q-2(i-j)\tau)}}e_{i+1,j}.
$$
The hypothesis on $p-q$ ensures that the denominator is non-zero for all $i,j$. 
In this expression $X(p \pm q + k \tau)$ is interpreted as the evaluation of $X$ at a point of $\CC$ in the following way 
(cf., \cite[\S3]{SS16}): write $x_{11}=x_0$, $x_{00}=x_1$, $x_{01}=x_2$, and $x_{10}=x_3$;
 if $X=\sum_{ab}\mu_{ab} x_{ab}$, then 
$X(z)= \sum_{ab}\mu_{ab} g_{ab}(z)$ for $z \in \CC$, where $g_{ab}$ is the function 
$$
g_{ab} = \gamma_{ab}\theta_{ab}(\t)\theta_{ab}(2\t)  \quad
\hbox{ and } \quad
 \g_{ab}= \begin{cases}
\sqrt{-1}, & \text{if $ab=00,11$}
\\
1, &\text{if $ab=01,10$.}
\end{cases}
$$

%%%%%%%%%%%%%%%%%%%%%%%%%%%%%%%%%%%%%%%%%%%%%%%%%%%%%%%%%%%%%%%%%
%%%%%%%%%%%%%%%%%%%%%%%%%%%%%%%%%%%%%%%%%%%%%%%%%%%%%%%%%%%%%%%%%
%%%%%%%%%%%%%%%%%%%%%%%%%%%%%%%%%%%%%%%%%%%%%%%%%%%%%%%%%%%%%%%%%

  \section{Modules over  $\NN$-graded $\Bbbk$-algebras}
  \label{sect.3}
  
In this section we gather some general results that have wider applicability. 
For the sake of brevity we do not state them in maximum generality. Certainly, many of them hold in much greater
generality than they are stated.
%Many of them are routine, and perhaps well known, 
 % so the only value to this section is that it organizes them and gathers them in a single place. 
  
 One theme running through this paper is the relation between modules that may not be graded, in particular finite-dimensional simple
 modules, and graded modules. Almost all results in this section are related to this theme, some directly, some indirectly. 
 A more geometric interpretation of this theme is discussed in the final subsection of this section, \S\ref{sect.final.rmks}.
  
Non-commutative rings fall into two distinct classes: those that are not finitely generated modules over their centers and those that are.
Actually, the division is a little more subtle than this: the latter class of algebras should be replaced by those that satisfy a 
polynomial identity. Every ring that is a finite module over its center satisfies a polynomial identity but the converse is false.
For example, subrings of a matrix ring over a commutative ring need not be finite modules over their centers but they satisfy a 
polynomial identity because the matrix ring does. Most rings satisfying a polynomial identity embed in rings that are finite modules over their centers, and most rings satisfying a polynomial identity have a localization that is a finite module over its center. 
When $\tau$ has finite order $A(E,\tau)$ is a finite module over its center. With applications to $A(E,\tau)$ in mind, many of the 
results in \S\ref{sect.3} concern rings that satisfy a polynomial identity. The very nicest such rings are the Azumaya algebras so
some results in this section concern Azumaya algebras, but almost always with the theme in the previous paragraph as
motivation.
  
  \subsection{}
  \label{sect.3.1}
Let $\Bbbk$ be an arbitrary field and $R=R_0+R_1+\cdots$ a left and right noetherian connected graded $\Bbbk$-algebra.
The connected hypothesis says that $R_0=\Bbbk$. We write $\fm$ for the ideal $R_1+R_2+\cdots$. Thus $R/\fm=\Bbbk$. An
$R$-module is said to be {\sf trivial} if it is annihilated by $\fm$. 

We also assume that $R$ is a graded quotient of an Artin-Schelter regular algebra, the latter being a connected graded $\Bbbk$-algebra having finite global homological dimension, finite Gelfand-Kirillov dimension, and satisfying the Gorenstein property. Basic 
properties of Artin-Schelter regular algebras can be found in \cite[\S2]{ATV2}. 

These assumptions on $R$ remain in force throughout \S\ref{sect.3}.

The 4-dimensional Sklyanin algebra is Artin-Schelter regular \cite{SS15}. 

%Since $R$ is noetherian it is a finitely generated $\Bbbk$-algebra: to see this observe that $\fm/\fm^2$ is a finitely generated
%$R/\fm$-module. 

We write $R_{\ge i}$ for  $R_i+R_{i+1}+\cdots$. This is an ideal of $R$. (The word ``ideal'' on its own always means ``two-sided ideal''.)
Since $R_{\ge i}/R_{\ge i+1}$ is an $R/\fm$-module, the noetherian hypothesis implies that each $R_i$ has finite dimension. 
Thus, if $M$ is a finitely generated graded $R$-module, then $\dim_\Bbbk(M_i)<\infty$ for all $i$, and  $M_i=0$ for $i\ll 0$. 
Associated to such a module is its {\sf Hilbert series} 
$$
H(M;t) \;=\; \sum_{i \in \ZZ} \dim_\Bbbk(M_i)t^i,
$$
the right-hand side being an element in the ring $\ZZ[t^{\pm 1}]$ of Laurent polynomials.

\subsection{Gelfand-Kirillov dimension and 1-critical modules}
\label{sect.gkdim.1-crit}
Let $R$ be as in the previous section. The {\sf Gelfand-Kirillov dimension} of $R$ is the number
 $$
 \GKdim (R)  \; :=\; \lim \sup_V \{\log_n (\dim_\Bbbk V^n)\}
 $$
 where the limit is taken as $n \to \infty$ and 
 the supremum is taken over all finite dimensional subspaces $V \subseteq R$ that contain $\Bbbk$ and generate $R$ as  a 
 $\Bbbk$-algebra and
 $$
 V^n={\rm span}\{v_1\ldots v_n \;| \; v_1,\ldots,v_n \in V\}.
 $$
It is clear that  $\GKdim (R) = 0$ if and only if $\dim_\Bbbk(R)<\infty$ and that  $\GKdim (R) \ge 1$ otherwise. 

Since $R$ is noetherian, if $\fp \subset \fq$ are distinct prime ideals, then $d(R/\fq)<d(R/\fp)$. 

 \begin{thm}
 [Small-Warfield,    \cite{SW13}]
 \label{thm.SW}
 Let $R$ be a finitely generated $\Bbbk$-algebra. If $R$ is prime and $\GKdim(R)=1$, then 
 $R$ is a finitely generated module over its center.
 \end{thm}

 The {\sf Gelfand-Kirillov dimension} of a finitely generated left $R$-module $M$ is  the number
 $$
d(M) \;=\;  \GKdim (M)  \; :=\; \lim \sup_{V,D} \{\log_n (\dim_\Bbbk V^nU)\}
 $$
 where the supremum is taken over all finite dimensional subspaces $V \subseteq R$ that generate $R$ as  a $\Bbbk$-algebra and 
 all finite dimensional subspaces $U \subseteq M$ that generate $M$ as an $R$-module.

\subsubsection{}
The hypotheses on $R$ imply that the Hilbert series of every finitely generated graded $R$-module is 
of the form 
$$
H(M;t) \;=\;  \frac{q_M(t)}{\prod_{1\le i \le s}(1-t^{d_i})}
$$ 
for some non-negative integers $d_i $ and some $q_M(t) \in \ZZ [t,t^{-1}]$ such that $q_M(1) \ne 0$ if $M \ne 0$, 
and only a finite number of $d_i$ occur as $M$ varies. With this assumption, the GK-dimension of every finitely generated 
graded $R$-module is the order of the pole of $H(M;t)$ at $t=1$.

\subsubsection{}
The {\sf multiplicity} of a graded $R$-module $M$ is the number
$$
e(M) \;:= \; q_M(1) \;=\; (1-t)^{d(M)} H(M;t)\vert_{t=1}.
$$
If $0 \to L \to M \to N \to 0$ is an exact sequence of finitely generated graded $R$-modules such that
$d(L)=d(M)=d(N)$, then $e(M)=e(L)+e(N)$. 

It is not difficult to see that $e(M)= q_M(1)/\prod_{1\le i\le
s}d_i$  when $H(M;t)$ has the above form. 
In particular, $e(M) \ge  \prod_{1 \le i \le s}d_i^{-1}$ and,  since only a finite number of
$d_i$ occur, there is a positive number $e_0$ such that
$e(M)\ge e_0$ for all  non-zero $M.$ 
It follows that if $M=M^0 \supset M^1 \supset M^2 \ldots$ is a strictly decreasing chain of finitely generated
graded modules, then 
$d(M^j/M^{j+1}) < d(M)$ for some $j$. This good property of multiplicity is used in Proposition \ref{prop.1.5}. 

Shift-equivalent $R$- modules have the same GK-dimension, and the same multiplicity when they are not finite 
dimensional. 

\subsubsection{}
An $R$-module  is {\sf $d$-critical} if it has of GK-dimension $d$ and every proper quotient of it has GK-dimension $<d$.

\begin{prop}
\cite[Prop. 2.30(vi)]{ATV2}
The annihilator of every finitely generated critical graded $R$-module is a prime ideal.
\end{prop}

It is easy to see that a finitely generated graded $R$-module $M$ is {\sf 1-critical} if 
and only if $\dim_\Bbbk (M)=\infty$ and $\dim_\Bbbk (M/M')<\infty$
for all non-zero graded $R$-submodules $M' \subseteq  M$.
The following facts are well-known and clear:
\begin{enumerate}
  \item 
  The only finite dimensional graded submodule of a 1-critical module is the zero submodule.
  \item 
Every non-zero homomorphism  between 1-critical modules is injective.
%  \item 
\end{enumerate}
Proposition \ref{prop.old.2.2} describes the close relationship between 1-critical graded $R$-modules and finite dimensional simple $R$-modules.

\subsection{Twisted modules and equivalent simple modules}
 The forgetful functor $F:\Gr(R) \to \Mod(R)$ that simply forgets the grading   has a right adjoint $G:\Mod(R) \to \Gr(R)$ 
 given by $GN=N \otimes \Bbbk[t^{\pm 1}]$ where $N \otimes \Bbbk[t^{\pm 1}]$ is graded by declaring that
  $\deg(N \otimes t^j)=j$, and  the action of a homogeneous element $a\in R_i$ on 
 a degree-$j$ element  $x \otimes t^j$ is $a.(x \otimes t^j)=(ax) \otimes t^{i+j}$. 
 The functor $G$ is exact. We will often want to use its (exact) subfunctor $N \mapsto \widetilde{N}$ where
 $$
 \widetilde{N} \;  :=\; N \otimes k[t] \; \subseteq \; N \otimes \Bbbk[t^{\pm 1}].
 $$
 
Since $R$ is a $\ZZ$-graded $\Bbbk$-algebra, $\Bbbk^\times$ acts
as automorphisms of $R$ by having $\lambda \in \Bbbk^\times$ act on $a \in R_i$ by $a \to \lambda^i a.$ 

This gives rise to an action of $\Bbbk^\times$
on the category of $R$-modules: if $M$ is a left $R$-module, we define the
{\sf twisted} module $M^\lambda$ to be $M$ as a $\Bbbk$-vector space 
with a new $R$-action given by $a*m=\lambda^i a.m$ for $a \in R_i$. 
(If $M$ is a graded module, then $M \cong M^\l$ for all $\l \in \Bbbk^\times.$)

Finite dimensional simple left $R$-modules $S$ and $S^\prime$ are 
{\sf equivalent} if $S^\prime \cong S^\lambda$ for some $\lambda \in \Bbbk^\times$.

We write $\Ann(M)$ for the annihilator of a module $M$.

%\subsection{The module $\widetilde{S}$}

%Let $S$ be a non-trivial finite dimensional   simple $R$-module. Define a graded left $R$-module 
%$$
%\wtS \; :=\; S \bigotimes \Bbbk[t]
%$$
%by declaring that $\deg (S \otimes \Bbbk t^i)=i$ and that $a \in R_n$ act by 
%$$
%a.(s\otimes t^i) \;=\;   (a.s)\otimes t^{i+n}.
%$$

\begin{prop}
\label{prop.old.2.2} 
Let $S$ be a non-trivial finite dimensional simple left $R$-module. 
\begin{enumerate}
  \item  
  $\wtS$ is a finitely generated %1-critical 
  graded $R$-module. 
  \item 
 The map $\pi:\wtS \to S$ defined by $\pi(s \otimes t^j) :=s$ is an $R$-module homomorphism with the following universal property: if $M$ is a graded left $R$-module such that 
 $M_{<0}=0$ and $\psi : M \to S$ is an $R$-module map, then there exists a unique degree 0 map 
 $\tilde \psi :M \to 
\wtS$ such that $\psi=\pi \circ \tilde \psi.$ 
  \item 
If $N$ is a 1-critical graded module  mapping onto $S$ and $N_{<0}=0$, then   $N$ embeds in $\wtS.$
\item{} 
Every  $S^\lambda$,  $\lambda\in \Bbbk^\times$, is a quotient 
of $\wtS.$
  \item 
 $\Ann( \wtS) $ is the unique prime graded ideal $\fp$ such that 
$\GKdim (R/\fp)=1$ and $\fp \subseteq \Ann (S).$  
  \item 
 $\Ann( \wtS)$ is the largest graded ideal contained in $\Ann (S)$.
\end{enumerate}
\end{prop}
\begin{pf}
(1)
% It is trivial to check that $\wtS$ is a graded  $R$-module. 
Since $\fm S \ne 0$, $\fm S=S$ and, since $\dim_\Bbbk S<\infty$, 
$$
R_1S+ R_2S+\cdots +R_n S \;=\;  S
$$
for some $n \ge 1$. We fix such an $n$ and will show that $\wtS$ is generated as an $R$-module by
$\wtS_0+\cdots+\wtS_{n-1}$.

\underline{Claim:} 
$R_{i+1}S+ \cdots +R_{i+n}S=S$ for all $i \ge 0$. 
\underline{Proof:} 
We chose $n$ so that the claim is true for 
$i=0$ and will argue by induction on $i$. To that end, suppose the claim is true for $i$. Now
$R_{i+2}S+ \cdots +R_{i+n}S +R_{i+1+n}S $ contains
\begin{align*}
&  R_{i+2}S+ \cdots +R_{i+n}S + R_{i+1}(R_1S+ R_2S+\cdots +R_n S)
\\
& \qquad \;=\;  R_{i+2}S+ \cdots +R_{i+n}S +R_{i+1}S 
\\
&  \qquad  \;=\;  R_{i+1}S  + R_{i+2}S+ \cdots +R_{i+n}S 
\\
&  \qquad  \;=\;  S.
\end{align*}
The claim therefore holds for all $i \ge 0$.

The truth of the claim implies that 
$$
R_{i+1}\wtS_{n-1}+ \cdots +R_{i+n}\wtS_0 \;=\;  (R_{i+1}S+ \cdots +R_{i+n}S) \otimes t^{i+n} \;=\;  \wtS_{i+n}
$$
for all $i \ge 0$. Hence $\wtS$ is generated by
$\wtS_0+\cdots+\wtS_{n-1}$.

(2) Define $\tilde \psi(m) = \psi(m) \otimes t^n$ for $m \in M_n.$

(3) This is a consequence of (2).

(4) If $\lambda\in \Bbbk^\times$, then $\wtS(t-\lambda)$ is a 
submodule of $\wtS$ and  $\wtS \tilde / S(t-\lambda) \cong  S^\lambda$. 

(5) 
Let $\fp$ be the sum of the graded 2-sided ideals that annihilate $S$. 
It is clear that $\fp$ annihilates $\wtS$. But $\Ann( \wtS)$ is a 
graded ideal contained in $\Ann ( S)$ so $\fp =\Ann( \wtS).$ 
Since $\Ann S$ is prime, it is easy to prove that $\fp$ is also prime. 
Replace $R$ by $R/\fp$ so we can assume that $0$ is the only homogeneous element in $\Ann( S)$. 

Suppose that  $\dim_\Bbbk( R_n) \to \infty$ as $n \to \infty.$ Since $\Ann(S)$ 
is of finite codimension in $R$, it follows that $\Ann(S) \cap R_n 
\ne 0$ for $n \gg 0.$ This is a contradiction, so we must have a bound 
on $\dim (R_n)$. But this forces $d(R) \le 1.$ If $d(R)=0$, then $R$ is
finite dimensional, whence the only simple $R$-module is the trivial
module. Since $S$ is not trivial, it follows that $\GKdim(R)=1$. 

If $\fq$ were another such ideal, then by definition of $\fp$, $\fq \subseteq 
\fp$, so $\GKdim(R/\fp)<\GKdim(R/\fq)$ which is absurd.
\end{pf}

{\bf Remarks. 1.} 
If $S$ is a finite dimensional simple $R$-module, then $e\big(\wtS\big)=\dim_\Bbbk( S)$. 
It therefore follows from Proposition \ref{prop.old.2.2}(3) that if  $S$ is a quotient of a  
1-critical graded module $F$, then $e(F) \le \dim S$. 

{\bf 2.} 
Suppose $R$ is generated by $R_1$ as a $\Bbbk$-algebra.
If $\dim_\Bbbk( S)=1$ and $S\not\cong R/\fm $, then $R/\Ann(\wtS)$ is a polynomial ring in one variable and $\wtS$ is isomorphic as a graded left 
$R$-module to  $R/\Ann(\wtS)$  (cf., \cite[Prop. 5.9]{LS8}).
\underbar{Proof}: 
Let $e$ be a non-zero element of $S$ and $\fa$ the left annihilator of $e \otimes 1 \in \wtS$. Clearly, $\fa$ is a graded left ideal. 
Every homogeneous element in $\fa$, and hence every element in $\fa$, annihilates {\it every} homogeneous component of $\wtS$ 
so $\fa \subseteq \Ann(\wtS)$; the reverse inclusion is obvious so $\fa=\Ann(\wtS)$. As graded left $R$-modules $R/\fa$ and 
$\wtS$ are isomorphic. Since $S$ is not annihilated by $\fm$ and $R$ is generated by $R_1$, there is an element,  $x$ say, in $R_1-\fa$. 
Since $R/\fa$ is a prime ring, the subalgebra of $R/\fa$ generated by $x$ is a polynomial ring in one variable. 
Thus $\Bbbk[x]$ has the same Hilbert series as $\wtS$ and therefore the same Hilbert series as $R/\fa$ so we conclude that the natural 
map $\Bbbk[x] \to R/\fa$ is an isomorphism of graded $\Bbbk$-algebras. 

{\bf 3.}
In general, $\wtS$ need not be a critical module. 
%\underbar{Proof:}
%Consider the example at the end of Section 1.
%Let $S$ be the simple left $B/( X-1)$ module, and consider $S$ as
%an $R$-module. Then $S$ remains simple as an $R$-module, and 
%$\wtS$ is as described in that example. In particular, $\wtS$ contains 
%a direct sum of the inequivalent 1-critical modules $C$ and $D$. 
 
{\bf 4.}
Sometimes $R/\fm$ is the only simple quotient of a 1-critical module.  For example,
consider the algebra $B=B(E,\sigma, {\cal L})$ (see \cite{ATV2} for the
definition) where $E$ is an elliptic curve, 
$\sigma \in \Aut E$, and $\cL$ is an invertible $\cO_E$-module of degree $\ge 3.$
If $\sigma$ has infinite order, then by \cite{ATV2}  the only finite dimensional simple $B$-module is the 
trivial module, and the only 1-critical $B$-modules are the point modules. By \cite[Lem. 5.8(d)]{LS8} 
such a point module does not have a non-trivial simple quotient.

{\bf 5.}
We will make frequent use of the observation that if $S$ is a non-trivial 
simple quotient of a 1-critical graded module $N,$ then $S$ is also a 
quotient of every non-zero submodule of $N$.

\subsection{}

We remind the reader that we are assuming $R$ is noetherian.

\begin{prop}
\label{prop1.1} 
Suppose $R$ is as above, and is also a prime ring.
If $M$ and $N$ are finitely generated, critical, graded left $R$-modules such that 
$$
\GKdim(M) \;=\;  \GKdim (N) \;=\; \GKdim(R),
$$ 
then $M(-i)$ embeds in $N$ for some $i \in \ZZ$.
Moreover, $i$ can be chosen arbitrarily large.
\end{prop}
\begin{pf}
Since $\GKdim (M) =  \GKdim(R)$,  $M$ is torsion-free in the sense that regular elements of $R$ do not annihilate any non-zero 
element in $M$.
Since $M$ and $N$ are critical they are uniform modules, i.e., non-zero submodules of $M$ have non-zero
intersection and similarly for $N$.  It now follows from \cite[Prop. 3.4.3(iv)]{MR9} that both  $M$ and $N$
are isomorphic to left ideals of $R$. Hence by \cite[Lem. 3.3.4(ii)]{MR9}
$M$ and $N$ contain isomorphic copies of one another. However, these
embeddings of $M$ and $N$ in one another may not respect the gradings.
Nevertheless ${\rm Hom}_R(M,N)$ is the sum of its homogeneous components, 
so there is some non-zero degree-preserving homomorphism  $f:M\to N.$ Since $M$ and $N$ are
critical of the same GK-dimension, $f$ is injective. If $\deg( f)=i$, then 
there is a degree zero injective map $M(-i)\to N.$

To see that $i$ may be chosen arbitrarily large, fix $i$ such that $M_i \ne 0$,
and choose $\ell\in \ZZ$ arbitrarily large. The hypotheses on $N$ ensure that
$\dim_\Bbbk(N/N_{\ge \ell})<\infty$ so $N_{\ge \ell}$ satisfies the same hypotheses
as $N$. By the first part of the proof, there is an $i$ such that $M(-i)$
embeds in $N_{\ge \ell}$. Since $M(-i)_{i+k} \ne 0$, it follows that
$(N_{\ge \ell})_{i+k} \ne 0$. Thus $i+k \ge \ell$. Since $i$ was fixed, and 
$\ell$ was arbitrary, this shows that we may choose $k$ arbitrarily large.
\end{pf}

\begin{prop}
\label{prop.1.2}
Suppose $R$ is as above, and is also a prime ring.
Let $F$ and $F^\prime$ be finitely generated critical graded $R$-modules such that 
$$
\GKdim \!  \left( \frac{R}{\Ann\,(F)}\right)   \;=\;  \GKdim \! \left( \frac{R}{\Ann\,(F')} \right)  \;=\; 1.
$$
Then $\Ann\,(F)= \Ann\,(F^\prime)$ if and only if $F$ and  $F^\prime$ 
are shift-equivalent.  
\end{prop}
\begin{pf}
$(\Rightarrow)$ 
Let $\fp$ be the common annihilator of $F$ and $F'$. 
Since $F$ is critical, $\fp$ is a prime ideal. We may replace $R$ by $R/\fp$ and assume that $\fp=0$. Thus 
$R$ is a  prime noetherian ring   and $\GKdim( F)=
\GKdim( F^\prime)=\GKdim(R)=1$. Therefore Proposition \ref{prop1.1} applies: there is 
a non-zero degree-zero map $F(n) \to F^\prime$  for some $n \in \ZZ.$ 
Since $F^\prime$ is 1-critical, the cokernel is finite dimensional, so 
$F(n)$ and $F^\prime$ are equivalent.

$(\Leftarrow)$ 
First, the annihilator of $F$ depends only on 
the equivalence class of $F$. To see this suppose 
$M$ is a non-zero submodule of $F$. If $\Ann(M)$ were strictly larger
than $\Ann( F)$, then $R/\! \Ann(M)$ would have finite dimension which would imply 
that $\dim(M) <\infty$ contradicting the fact that $F$ is 1-critical. 
Thus  $\Ann(M)=\Ann( F)$.

If $F(n) \sim F^\prime$  there is a non-zero graded module, $N$ say,
such that both $F$ and $F^\prime$ contain a copy of suitable shifts of $N$. 
By the previous paragraph it follows that $\Ann( F)=\Ann (N)=\Ann (F^\prime)$. 
\end{pf}

\begin{thm}
\label{thm.old.1.3}
Suppose $R$ is as above, and is also a prime ring.
If $\GKdim(R)=1$ there is a finitely generated 1-critical graded $R$-module $M$ such that
\begin{enumerate}
  \item 
  $M$ is generated in degree zero;
  \item 
up to isomorphism, the simple objects in $\QGr(R)$ are $\{M(n) \, \vert \, 0\le n < p \}$ for some integer $p$; % i.e., there is a unique 1-critical $R $-module up to shift equivalence;
  \item 
 the Hilbert series of $M$ is periodic.
\end{enumerate}
Moreover, the degree of every homogeneous central element in $R$ is a multiple of $p$.
\end{thm}
\begin{pf}
Let $F$ be a 1-critical $R $-module and $M$ a 1-critical $R$-module that is generated in degree 0.
Such an $M$ exists: for example, if $I$ is a graded left ideal that is maximal subject to the condition that $\dim_\Bbbk(R/I)=\infty$
we could take $M$ to be $R/I$. The noetherian hypothesis implies that every proper quotient ring of $R$ has finite dimension, 
so the annihilators of $M$ and $F$ are zero. By Proposition \ref{prop.1.2}, $F \sim M(n)$ for some $n$; hence 
$R $ has a unique 1-critical module up to shift-equivalence.  

By Proposition \ref{prop1.1}, there is an injective degree 0 map $\varphi:M \to M_{\ge 1}(k)$ 
for some $k \ge 1,$ whence $M \sim M(k)$.  It follows that $M(j)\sim
M(i+k)$ for all $i \in \ZZ.$ Let $p\in \NN$ be minimal such that $M
\sim M(p).$ Then $\dim_\Bbbk (M_n)=\dim_\Bbbk (M_{n+p})$ for all $n \gg 0$. Now replace
$M$ by $M_{\ge n}(-n)$ for some suitably large $n$, to obtain a module with periodic Hilbert series.

If $z$ is a homogeneous central element of degree $d$, then multiplication by $z$ gives an injective homomorphism 
$C \to C(d)$ so $C \sim C(d)$, whence $d \in \ZZ p.$
\end{pf}

\subsubsection{}
{\bf Remarks.} 
{\bf 1.} 
We call the smallest integer $p$ in Theorem \ref{thm.old.1.3}(b) the {\sf period} of $R $; the 
irreducible objects of $\QGr(R )$ are in bijection with the elements of $\ZZ/p$ and they are all shift-equivalent to each other.

{\bf 2.}
The example in \S\ref{MV.example} shows it is possible for non-isomorphic point modules
to become isomorphic in $\QGr(R)$. However, it is an easy consequence of Proposition \ref{prop1.1} that this phenomenon 
does not occur over a noetherian ring satisfying a polynomial identity. The following more general statement holds.

\begin{lemma}
\label{lem.pt.mods}
Let $M$ and $N$ be equivalent point modules over $R$.  
If $R/\! \Ann(M)$ or $R/\! \Ann(N)$ has GK-dimension 1, then $M \cong N$. 
\end{lemma}
\begin{pf}
Since $M$ and $N$ are critical and $R$ is 2-sided noetherian, 
$\Ann( N_{\ge k})=\Ann( N)$ and $\Ann(M_{\ge k}) = \Ann(M)$ for $k$ sufficiently 
large. Thus $\Ann (N)=\Ann(M)$ since they have a common
submodule. The ring $R/\! \Ann(M)$  is prime and $\GKdim(R/\! \Ann(M))=1=\GKdim(M)=\GKdim( N)$ by
hypothesis so, by Proposition \ref{prop1.1}, suitable shifts of $N$ and $M(1)$ embed in $M$ via
degree zero maps. Thus for some $i,j>0$ both $N(-i) $ and $M(-j) $ embed in
$M$. Therefore, $M(-ij)$ and $N(-ij) $ embed in $M$. But the image of both
these embeddings is $M_{\ge ij}$ whence $M(-ij)  \cong N(-ij) $ as required.
\end{pf}

\subsubsection{Example} (M. Vancliff)
\label{MV.example}
Without the 2-sided noetherian hypothesis most of the foregoing fails
(essentially because $R $ need not be left ideal invariant for
GK-dimension).
If $R =\CC\langle x,y\rangle/( x^2, xy)$, then $M=R /R x$ and $N=R /R y$ are 
non-isomorphic point modules that are equivalent. It is easy to see that
$M_{\ge 1} \cong M(-1) \cong N_{\ge 1}.$ The ring $R $ is left
but not right noetherian: the 2-sided ideal $R x$ is not finitely generated
as a right ideal. The proof in Lemma \ref{lem.pt.mods} does not apply to this example
because $R x=\Ann (N_{\ge 1}) \ne \Ann( N)$. In particular, the annihilator of
the critical module $R/Ry$ is not prime.

\begin{prop}
Suppose $R$ is as above, and also prime.  If $\GKdim(R)=1$, then the center of $R $ is isomorphic to a 
subalgebra of the polynomial ring $\Bbbk[t]$ generated by various $t^i$s. 
In particular, its ring of fractions $\Fract(R )$ contains a copy of the rational function field $\Bbbk(t).$
\end{prop}
\begin{pf}
Because $R $ is prime, its center is a domain. Because
$\GKdim(R )=1$, $R $ is a finitely generated module over its center (Theorem \ref{thm.SW}).
Let $m$ and $n$ be positive integers and suppose that $u \in R _m$ and $v \in R _n$ 
are non-zero central elements. Let
$X$ and $Y$ be commuting indeterminates, both of degree $mn$, and
let $\psi:\Bbbk[X,Y] \to R $ be the homomorphism from the polynomial ring given by $\psi(X)=u^n$ and $\psi(Y)=v^m.$
The kernel of $\psi$ is a non-zero homogeneous prime ideal hence of the form $(\a Y-\beta X)$ 
for some $0 \ne \alpha, \beta \in \Bbbk$. Therefore 
$\alpha u^n = \beta v^m$. A similar argument shows that in any fixed
degree the center of $R $ is either $0$ or 1-dimensional. The result
follows.
\end{pf}

\subsubsection{Example}
 Let $\CC[t]$ be a polynomial ring  with $\deg t=p.$ Define a 
grading on the $p \times p$ matrix algebra $M_p(\CC)$ by declaring that $\deg e_{ij}=i-j$ for each
matrix unit $e_{ij}$
and give $M_p(\CC) \bigotimes \CC[t]$ the tensor 
product grading. Let $R $ be the $\CC$-subalgebra generated by the 
degree 1 part of $M_p(\CC) \bigotimes \CC[t].$ Then $R $ is generated over
$\CC$ by  $x_i:=e_{i+1,i}$ for $1 \le i \le p-1$ and $x_p:=te_{1p}.$
As a subalgebra of $M_p(\CC[t])$, $R $ consists of those matrices whose
strictly upper triangular entries are multiples of $t$, whose diagonal 
entries have the same constant term, and whose lower triangular 
entries are arbitrary.

It is not hard to show that  $R $ is prime noetherian with $\GKdim(R )=1$, that
the center of $R $ is $\CC[t]$, and that $gl.dim(R )=\infty$ if $p>1$. 
Furthermore, 
$R $ is a monomial algebra defined by the relations $\{x_ix_j=0 \, | \,
i \ne j+1 (\mod\, p)\}$. A basis for $R $ is given by words of the form
$\ldots x_3x_2x_1x_px_{p-1}\ldots x_2x_1 \ldots$. The augmentation 
ideal of $R $ is the direct sum $R x_1 \oplus \ldots \oplus R x_p.$

There are $p$ different point modules, namely the modules
$$
M^i \;=\;  R /R x_1+\ldots + R  \hat x_i+\ldots +R x_p
$$
where $\hat x_i$ means that $x_i$ is omitted. These give the distinct 
irreducible objects in $\QGr(R )$. In terms of the generators 
$$
t \;=\;  x_px_{p-1}\ldots x_2x_1 \, +\, x_{p-1}x_{p-2} \ldots x_1x_p \, + \, 
\ldots \, + \, x_1x_p\ldots x_3x_2.
$$ 
Furthermore, $R /( t)$ is isomorphic to the ring
of lower triangular matrices over $\CC$ with constant diagonal entries.
Fix $\omega = e^{2\pi i/p}$. Let $G \subseteq  \Aut(R )$ be the subgroup 
generated by the automorphisms $\rho:x\to \omega x$
for $x \in R _1$ and $\tau:x_i \to x_{i+1(mod\,p)}$. Then the invariant
subring $R ^G$ is  the center of $R $.

\subsection{``Composition series'' of modules having GK-dimension one}
Since 1-critical modules are irreducible objects in $\QGr(R)$ they
behave rather like simple modules. For example, if $M$ is a 1-critical graded $R$-module, then
$\End_{\Gr(R)}(M)$ is a division ring: a non-zero degree-preserving endomorphism $\varphi$ is
injective because $M$ is critical, and therefore  surjective because
each $M_n$ has finite dimension. 

Another example of the similarity is the existence of ``composition series'', not in the usual sense, of course. 
Let $M$ be a finitely generated graded $R$-module of
GK-dimension 1. A {\sf composition series} for $M$ is a finite chain
of graded submodules $M=M^0\supseteq M^1 \supseteq \cdots \supseteq M^k =0$ such
that each $M^j/M^{j+1}$ is either finite dimensional or 1-critical.
Such a chain exists because multiplicity is additive on short exact
sequences in which all the modules have the same GK-dimension. The
{\sf composition factors} associated to this series are those
$M^j/M^{j+1}$ that are 1-critical; we consider the composition
factors as objects of $\QGr(R)$.

\begin{prop}
\label{prop.1.5}
If $M$ is a finitely generated graded
$R$-module of GK-dimension 1, then $M$ has a finite composition
series, and the composition factors are uniquely determined up to isomorphism
in $\QGr(R)$ by $M$.
\end{prop}
\begin{pf}
Suppose $M=M^0\supseteq M^1 \supseteq \cdots \supseteq M^k =0$ and
$M=N^0\supseteq N^1 \supseteq \cdots \supseteq N^\ell=0$ are composition series.
By the Butterfly Lemma,
$$
H_{ij}  \; :=\; 
{{M^j\cap N^i + M^{j+1}} \over {M^j \cap N^{i+1} +M^{j+1}}}  \; \cong \;
{{N^i \cap M^j + N^{i+1}} \over {N^i \cap M^{j+1} + N^{i+1}}}.
$$
(One must check that all the maps in the proof of the Butterfly Lemma have degree 0 so that this is actually an isomorphism in $\Gr(R)$.)
In particular, $H_{ij}$ is a subquotient of both $M^j/M^{j+1}$ and
$N^i/N^{i+1}$. Fix $j$ such that $M^j/M^{j+1}$ is 1-critical.
Then $H_{ij}$ is infinite dimensional for some 
$i$, whence it is a submodule, not just a subquotient, of
$M^j/M^{j+1}$ since the latter is 1-critical. Similarly, $H_{ij}$ is
also a submodule of $N^i/N^{i+1}$, whence $M^j/M^{j+1}$
is equivalent to $N^i/N^{i+1}.$ For each such $j$ there is a
unique $i$ such that $\dim_\Bbbk( H_{ij})=\infty$, so we obtain a bijection between
the composition factors arising from the two composition series, and
corresponding composition factors are equivalent (i.e. isomorphic in 
$\QGr(R)$). 
\end{pf}

\begin{cor}
\label{cor.1.5}
If $L$ be a 2-critical graded $R$-module  having
infinitely many inequivalent 1-critical quotients, $\{ C^\a \, \vert
\, \a \in \L\}$ say. Then
$$
\Ann (L) \;=\;  \bigcap_{\a \in \L} \Ann (C^\a).
$$
\end{cor}
\begin{pf}   
Write $C^\a=L/L^\a$.
If  the intersection of all the $L^\a$'s is non-zero, then 
$$
\GKdim \left( \frac{L}{\bigcap_\a L^\a} \right) \;=\;1. 
$$
This quotient of $L$ has infinitely many inequivalent quotients so has
infinitely many composition factors thereby contradicting Proposition \ref{prop.1.5}.
We conclude that $\bigcap_{\a \in \L} L^\a=0$.
Hence there is a diagonal embedding $L\to \prod_{\a \in \L}
C^\a$. It follows that $\Ann (L)= \bigcap_{\a \in \L} \Ann (C^\a)$.
\end{pf}

\begin{prop}
\label{prop.old.2.3} 
Assume $R$ is as above.
 Two 1-critical graded $R$-modules have a common non-trivial simple quotient if and only if they are shift-equivalent. 
In this case, they have exactly the same non-trivial simple quotients.
\end{prop}
\begin{pf}
Let $S$ be a non-trivial simple quotient of 1-critical graded
modules $F$ and $F^\prime$. Let $\fp=\Ann (\wtS)$.
Since   $F$ and $F^\prime$ embed in $\wtS$, it follows that 
$\fp=\Ann( F)=\Ann (F^\prime)$.
 Hence $F(n)\sim F^\prime$ by Proposition \ref{prop.1.2}.

If $F$ and $F'$ are shift equivalent 1-critical modules, then there
is a 1-critical graded module $N$ such that $N(i)$ embeds in $F$
and $N(j)$ embeds in $F'$ for some $i,j \in \ZZ$. Therefore $N(i)$,
and consequently $N$, has the same non-trivial simple
quotient modules as $F$. The same argument applies to
$N(j)$ and $F'$, so $F$ and $F'$ have the same non-trivial
simple quotients.
\end{pf}

\begin{prop}
\label{2.4}
Assume $R$ is as above.
A finitely generated 1-critical graded $R$-module $C$ has a non-trivial finite dimensional simple quotient if and only
if $d(R/\! \Ann( C))=1.$
\end{prop}
\begin{pf} 
$(\Rightarrow)$
 Let $S$ be a non-trivial finite dimensional simple quotient of $C$. 
 By  Proposition \ref{prop.old.2.2}(3), $C$ embeds in $ \wtS$, whence $\Ann(S) \supseteq \Ann (C) \supseteq
\Ann( \widetilde S)$.  But $\Ann( \widetilde S)$ is the largest graded ideal contained in $\Ann(S)$ by Proposition \ref{prop.old.2.2}(6), 
and $\Ann (C)$ is a graded ideal so $\Ann (C) = \Ann( \widetilde S)$.
By Proposition \ref{prop.old.2.2}(5), $R/\Ann(\widetilde S)=1$. 

$(\Leftarrow)$ 
Consider $C$ as a module over $R^\prime=R/\! \Ann( C)$.
 By Theorem \ref{thm.SW}, $R^\prime$ has a homogeneous central element,  $z$ say, of positive degree.
The ring $R'$ is prime so its center is a domain. The subalgebra $\Bbbk[z]$ of $R'$ generated by $z$ is therefore a polynomial ring in one variable. Since $C$ is a torsion-free graded $\Bbbk[z]$-module and $C_{\le n}=0$ for $n \gg 0$,  $C$ must be a free $\Bbbk[z]$-module. 
It follows that  $C/(z-\nu)C$ is non-zero for all $ \nu \in \Bbbk$. But $C/(z-\nu)C$ has finite dimension since $C$ is 1-critical. 
Thus, for all non-zero $\nu$, $C/(z-\nu)C$ has a finite dimensional simple quotient.
\end{pf}

\begin{thm}
%{\bf Theorem 2.5.} 
\label{2.5}
Assume $R$ is as above, and is also prime and has GK-dimension 1.
If $\Bbbk$ is algebraically closed, then all non-trivial finite dimensional simple $R$-modules are equivalent.
\end{thm}
\begin{pf} 
There is a central element $0 \ne z \in R$ such that the localization $R_z=R[z^{-1}]$ is an Azumaya algebra over its center  
\cite[Prop. 13.7.4]{MR9}. 
Since Formanek's central polynomial is homogeneous \cite[\S\S13.5-13.7]{MR9} we can assume that $z$ is homogeneous.
We will do that.

If $S$ is a simple $R$-module, then $\dim_\Bbbk(S)<\infty$ so $\End_R( S)=\Bbbk$.
Therefore, either $z$ annihilates $S$ or acts on $S$ as multiplication by a non-zero scalar. 
The first possibility occurs if $S$ is trivial; conversely if $zS=0$ then $S$ is an $R/( z )$-module, but this 
ring is finite dimensional and graded so its only simple module is the 
trivial one. If $S$ is not trivial then $z$ acts as a non-zero scalar, so
the $R$-action extends to $R_z$, making $S$ a simple $R_z$-module. More
precisely $R_z \otimes_R S \cong S$ as $R$-modules. 
Hence the functor $R_z\otimes_R -$ induces a bijection between
non-trivial simple $R$-modules and simple $R_z$-modules. To see 
that every simple $R_z$ module arises in this way, let $V$ be a simple
$R_z$-module, choose $S\subseteq  V$ a simple $R$-submodule, and check that
$R_z \otimes S \cong V$.

If $0\ne \lambda \in \Bbbk$, then we can also twist $R_z$-modules by $\lambda$
and $(R_z\otimes_R S)^\lambda \cong R_z\otimes_R S^\lambda$. Hence the proof of the
theorem will be complete once we show that $R_z$ has only one equivalence class
of simple modules. 
	
The center of $R_z$ is a $\ZZ$-graded 1-dimensional commutative domain 
containing units that are not in $\Bbbk$, and is a
finitely generated $\Bbbk$-algebra since $R_z$ is.  But the only such
commutative domain is $\Bbbk[x,x^{-1}]$ with grading induced by $\deg( x) =d$
say. However, $\Bbbk$ is algebraically closed so, by Tsen's Theorem,
the Brauer group of $\Bbbk(x)$ is trivial. By \cite[Thm 7.2]{AG4}
 there is an injection of Brauer groups  $\Br(\Bbbk[x,x^{-1}]) \to \Br(\Bbbk(x))$ whence every
Azumaya algebra over $\Bbbk[x,x^{-1}]$ is  a matrix algebra. Thus
$R_z \cong M_n(\Bbbk) \otimes \Bbbk[x,x^{-1}].$ It is easy to see that
$R_z/(x-\lambda^d) \cong \bigl(R_z/(x-1)\bigr)^\lambda$ for all $\l \in \Bbbk^\times$ from which it
follows that $R_z$ has only one equivalence class of simple modules.
\end{pf}

\begin{cor}
\label{cor.old.remark}
Assume $R$ is as above and  that $\Bbbk$ is algebraically closed. 
If a 1-critical graded $R$-module $C$ has a non-trivial finite dimensional simple quotient, then 
it has exactly one equivalence class of such quotients. 
\end{cor}
\begin{pf}
The point to be checked is that Lemma \ref{2.4} implies that Theorem \ref{2.5} can be applied to $R/\! \Ann (C)$.
\end{pf}

\begin{cor}
\label{cor.old.2.6} 
Assume $R$ is as above and  that $\Bbbk$ is algebraically closed. 
Let
\begin{equation}
\label{set.of.graded.primes}
{\cal P}=\{\fp \in \Spec (R) \; \vert \; \fp \hbox { is graded and } \GKdim (R/\fp)=1\}.
\end{equation}
There are bijections between the following sets:
\begin{enumerate}
\item{} equivalence classes of non-trivial finite dimensional simple
$R$-modules;
\item{} shift-equivalence classes of 1-critical graded modules that
have a non-trivial simple quotient;
\item{} shift-equivalence classes of 1-critical graded modules 
$C$ such that 
\newline 
$\GKdim(R/\! \Ann( C))=1;$
\item{} elements of ${\cal P}$.
\end{enumerate}
\end{cor}
\begin{pf}
If $S$  is a non-trivial finite dimensional simple there is a
1-critical module $C$ that maps onto $S$. By Proposition \ref{prop.old.2.3}, the shift
equivalence class of $C$ is uniquely determined by $S$. On the other hand,
if $C$ is a 1-critical module having a non-trivial finite dimensional simple quotient
then, by Corollary \ref{cor.old.remark}, $C$ determines a unique equivalence class of
simple modules. This proves the bijection between the sets in (1) and (2).

By Theorem \ref{2.5}, the sets in (2) and (3) are  the same.

If $\fp\in {\cal P}$,  then $\fp$ is the annihilator of a graded 1-critical 
$R$-module: in fact, $\fp$ is the annihilator of every 1-critical composition factor of $R/\fp$. The bijection between the sets in 
(3) and (4) now follows from Proposition \ref{prop.1.2}.
\end{pf}

The set in Corollary \ref{cor.old.2.6}(b) consists of {\it all} 1-critical graded $R$-modules if $R$ also 
satisfies a polynomial identity.

\begin{prop}
\label{prop.2.7}
Assume $R$ is as above and  that $\Bbbk$ is algebraically closed. 
If  
\begin{equation}
\label{enough.simples}
 \bigcap_{
 \begin{subarray}{c} 
  \hbox{all finite diml} \\ \hbox{simple $S$}
  \end{subarray}}
  \Ann (S) \; = \; 0,
\end{equation}
 then
\begin{enumerate}
\item{} 
$\bigcap_{\fp \in {\cal P}} \fp =\bigcap \Ann( C)=0$ where the second 
intersection is taken over all 1-critical graded $R$-modules $C$;
\item{}
 if $\GKdim(R) \ge 2$, then the set ${\cal P}$ in (\ref{set.of.graded.primes}) is infinite;
\item{}
 if $\GKdim(R) \ge 2$, then $R$ has infinitely many 
shift-inequivalent graded 1-critical modules.
\end{enumerate}
\end{prop}
\begin{pf}
 (1) 
Let $\fp\in {\cal P}$. Since  $\fp$ is the annihilator of a graded 1-critical 
$R$-module, $\bigcap_{\fp\in {\cal P}} \fp \supseteq \bigcap \Ann (C)$ where the right-hand intersection is taken over all 1-critical graded
$R$-modules $C$. 
Suppose that $x \in \bigcap \fp.$ Then for every finite dimensional
simple module $S$, we have $x. \wtS =0$ by Proposition \ref{prop.old.2.3}. Thus
$x$ annihilates every finite dimensional simple $R$-module.  It follows from (\ref{enough.simples}) that $x=0.$ 

(2) If ${\cal P}$ were finite, then $0=\fp_1 \cap \ldots \cap \fp_n$ with
$d(R/\fp_j)=1$ for all $j$, so $\GKdim(R) =1.$ Since this is not the case ${\cal P}$ must be infinite.

(3) 
If $\fp\in {\cal P}$ there is a 1-critical graded module $C$ such that $\fp=\Ann( C)$. 
 It now follows from (2) and Proposition \ref{prop.1.2} that there 
are infinitely many shift-inequivalent 1-critical graded modules.
\end{pf}

\begin{lem} 
\label{2.8}
Assume $R$ is as above and is also prime. Let $Z$ be the center of $R$ and 
suppose $R$ is a finitely generated $Z$-module. Let $M$ be a finitely generated torsion-free $R$-module
and $\fq \in \Spec( Z).$
Then there exists $\fp \in \Spec( R)$ such that $\fp\cap Z =\fq$ and $M/\fq M$
has a critical quotient, $C$ say, such that $\Ann( C)=\fp$.
\end{lem}
\begin{pf}
By the Artin-Tate Lemma, $Z$ is a finitely generated $\Bbbk$-algebra.

We consider supports of modules. Since $M$ is a finitely generated
torsion-free $Z$-module, $\Supp(M)=\Spec( Z).$ Thus 
$$
\Supp\Bigg(\frac{M}{\fq M}\Bigg)   \;=\;  \Supp\Bigg(\frac{Z}{\fq} \otimes_Z M\Bigg)  \;=\;  \Supp\Bigg(\frac{Z}{\fq} \Bigg)
\cap\Supp(M) \;=\;  \cV(\fq)
$$
where $\cV(\fq)=\{\fq'\in \Spec(Z) \; | \; \fq'\supseteq\fq\}$.
Thus $d(M/\fq M)=d(Z/\fq)$. Hence $M/\fq M$ has an $R$-module
quotient, $C$ say, that is $d(Z/\fq)$-critical. Let $\fp=\Ann C$, which is
a prime ideal because $C$ is critical. Since $R$ satisfies a
polynomial identity, $d(Z/\fq)=\GKdim( C)=\GKdim(R/\fp)$.  Since $\fq \subseteq  \fp \cap Z$
it follows that $R/\fp$ is a finitely generated $Z/\fq$-module. Thus
$\GKdim(R/\fp) \le \GKdim(Z/\fq)$ so it follows that $\GKdim(R/\fp)=\GKdim(Z/\fq)$. Therefore
$\fp\cap Z=\fq$.
\end{pf}

\begin{cor}
\label{2.9}
Assume $R$ is as above and that it is a finitely generated module over its center $Z$.
Let $M$ be a finitely generated 2-critical graded $R$-module.
Then $M$ has infinitely many shift-inequivalent 1-critical graded 
quotients.
\end{cor}
\begin{pf}
Replacing $R$ by $R/\! \Ann(M)$, we may assume that $M$ is torsion-free
and $R$ is prime. By Proposition \ref{prop.2.7}(2), $R$ has infinitely many graded
height-1 prime ideals, as does $Z$. Furthermore, if $\fq$ is a
graded height-1 prime ideal of $Z$ then there is only a finite
number of graded primes $\fp \in\Spec(R)$ such that $Z\cap \fp=\fq$.

Now follow the proof of Lemma \ref{2.8} making use of the fact that if $N$
is a graded module of GK-dimension $d$, then there is a graded $d$-critical
quotient of $N$. Thus infinitely many different
graded height 1 primes arise as the annihilators of 1-critical
quotients of $M$. The result now follows from Corollary \ref{cor.old.2.6}.
\end{pf}

\begin{prop}
\label{2.10} 
Assume $R$ is as above, and $\Bbbk$ is algebraically closed. Suppose $R$ is also a prime ring  
of GK-dimension 1.
 If $R$ has a unique 1-critical module up to shift-equivalence, $F$ say, then
\begin{enumerate}
\item{} $\dim_\Bbbk( F_n)$ is constant for $n\gg 0$, say $\dim F_n=e$, and
\item{} every non-trivial simple $R$-module has dimension $e$.
\end{enumerate}
\end{prop}
\begin{pf}
(1) 
Since $F$ is shift-equivalent to $F(1)$, $F(1)_{\ge n} \cong F_{\ge n}$ for $n\gg 0$ whence $\dim(F_n)=\dim( F_{n+1})$ for $n \gg 0.$ 

(2) 
By Proposition \ref{prop.old.2.3}, shift-equivalent modules have the same simple quotients so we can
replace $F$ by a suitable $F_{\ge r}(r)$ and assume that 
$\dim( F_n)=e$ for all $n \ge 0.$ Since $F\sim F(1)$ there is an isomorphism $F_{\ge r} \cong F_{\ge r+1}(1)$ for $r \gg 0$. Replacing
$F$ by $F_{\ge r}$ there is a non-zero degree-one injective $R$-module homomorphism $\varphi:F \to F$.
Clearly
$$
\dim_\Bbbk \left( \frac{F}{\Im(\varphi -\id_F)} \right) \;=\; e
$$
so $F$ module has a simple quotient, $S$ say, of dimension $\le e.$ However, $F$ embeds in $\tilde
S$ so $e \le \dim_\Bbbk(  S)$. Hence $\dim_\Bbbk(  S)=e$. By Theorem \ref{2.5}, every other non-trivial
simple quotient of $F$ also has dimension $e$. Since every non-trivial simple $R$-module is a quotient of $F$, the result follows.
\end{pf}

The argument in Proposition \ref{2.10}(2) may be extended to prove that if $F$ is 1-critical
and $F \cong F(n)$ for some $n \ge 1$, then every non-trivial simple
quotient $S$, of $F$, satisfies $e(F) \le \dim (S) \le n.e(F).$

\begin{lem} 
Let $F$ be a finitely generated 1-critical graded $R$-module.
If  there is a central homogeneous element $z$ of positive degree such that $zF
\ne 0$, then $F$ has a non-trivial finite dimensional simple
quotient.
\end{lem}
\begin{pf}
The ring $\End_R(F)$ of {\it all} $R$-module endomorphisms of $F$, not just
the degree-preserving ones, is a graded ring with $\deg(\psi)=d$ if  $\psi(F_i) \subseteq  F_{i+d}$ for all $i$.
The map $\Phi:\Bbbk[z] \to \End_R (F)$ that sends $w$ to the map ``multiplication by $w$''
  is a homomorphism of graded $\Bbbk$-algebras.

Since $F$ is critical, $\End_R (F)$ is a domain. Hence $\ker(\Phi)$ is
a homogeneous prime ideal of $\Bbbk[z]$. By hypothesis $\Phi(z) \ne 0$ so $\ker \Phi=\{0\}$. 
	
Since $F$ is zero in sufficiently negative degree, $(z-\mu)F \ne F$
for all $\mu \in \Bbbk$. Thus, if $\mu \ne 0$, then $F/(z-\mu)F$ is a non-zero finite
dimensional $R$-module  having a non-trivial simple quotient module.
\end{pf}

\subsection{Remarks}
\label{sect.final.rmks}
The results in \S\ref{sect.3} are analogues of elementary relationships between a projective variety $X$ and the
an affine variety $Y$ that is a cone over it. 
If $R$ is the coordinate ring of $Y$, the $\GG_m$-action on $Y$ corresponds to an $\NN$-grading on $R$.   
Points on $Y$ are of two types: the vertex $0 \in Y$ which corresponds to the trivial $R$-module $R/\fm$, and the
others which correspond to the non-trivial simple $R$-modules. If $S$ is a non-trivial simple $R$-module, 
then the twists $S^\l$, $\l \in \Bbbk^\times$, correspond to the points in the $\GG_m$-orbit of the point corresponding to $S$.

Now drop the hypothesis that $R$ is commutative. 

Finitely generated 1-critical graded $R$-modules play the role of the coordinate rings of the
closures of the $\GG_m$-orbits of the points in $Y-\{0\}$. The notion of a critical module is analogous to the notion of an irreducible
topological space. One fundamental difference from the commutative case is that if $\fp$ is a 
graded prime ideal of $R$ such that $\GKdim(R/\fp)=1$ there may be more that one 1-critical graded $R$-module whose annihilator is 
$\fp$. 

Every 1-critical graded $R$-module becomes a simple object in $\QGr(R)$. Such simple objects are the analogues of the skyscraper
sheaves and therefore analogues of the points on $X$. 

$R$-modules of GK-dimension play the role of curves on $Y$ and 
graded modules of GK-dimension 1 play the role of $\GG_m$-stable curves in $Y$. The image of such a module in $\QGr(R)$ has finite
length and the ``composition factors'' of the  graded module of GK-dimension 1 correspond in a precise way to the composition factors of
that finite length object.

%%%%%%%%%%%%%%%%%%%%%%%%%%%%%%%%%%%%%%%%%%%%%%%%%%%%%%%%%%%%%%%%%
%%%%%%%%%%%%%%%%%%%%%%%%%%%%%%%%%%%%%%%%%%%%%%%%%%%%%%%%%%%%%%%%%
%%%%%%%%%%%%%%%%%%%%%%%%%%%%%%%%%%%%%%%%%%%%%%%%%%%%%%%%%%%%%%%%%

  \section{The Sklyanin Algebra is finite over its center when $\vert \tau \vert < \infty$}
 % {Section 3. }
 \label{old.sect.3}

In this section $A$ denotes  $A(E,\tau)$ and $p$ and $q$ denote points on $E$. Recall that $\fa_{pq}$ denotes the annihilator of 
$M(p,q)$.

\subsection{} 
In this section we show that $A$ is a finitely generated module over its center when $\tau$ has finite order. 
We described the general strategy of the proof in \S\ref{sect.strategy}.
 We first show $A$ satisfies a polynomial identity by showing that it   embeds in a matrix algebra over a commutative ring; 
 we do this by showing that the intersection of the annihilators of all modules of dimension $\le 2n$ is $\{0\}$.
 % that É and  then Stafford's Theorem  \cite{JTS17}. 

The first step towards this is to show that $A/ \fa_{pq}$ satisfies a polynomial identity
when $p-q \notin 2\ZZ\t$. This restriction on $p-q$ allows us to use the basis $\{e_{ij} \;|\; i,j \ge 0\}$ for $M(p,q)$ 
that was described in \S\ref{sect.good.basis}. 
Of course, we eventually wish to show that $A/ \fa_{pq}$ satisfies a polynomial identity for {\it all} $p$ and $q$, 
not just for these ``general'' ones.  It is possible to do this. A good comparison to keep in mind is the following. 
To prove that a  ``reasonable''  prime algebra
of GK-dimension $d$ satisfies a polynomial identity of degree $\le
2r$ it suffices to show that there exists a $d$-dimensional variety
parametrizing some simple modules of dimension $r$; in particular,
one does not first attempt to prove that {\it all} simple modules
are of dimension $\le r$, but just that ``most'' are. One then
concludes via the Amitsur-Levitski Theorem \cite[Thm. 13.3.3]{MR9} that all simple modules have
dimension $\le r$. Several of our arguments are in this spirit.

\begin{proposition}
\label{prop3.1}
Suppose  $\tau$ has order $n<\infty.$
If $p-q \notin 2\ZZ\tau$, then
\begin{enumerate}
\item{} $d(A/ \fa_{pq})=2$,
\item{} $e(A/ \fa_{pq}) \le n^2$, and
\item{} $A/ \fa_{pq}$ satisfies a polynomial identity.
\end{enumerate}
\end{proposition}
\begin{pf}
Since $p-q \notin 2\ZZ\tau$, we can use the basis $\{e_{ij} \;|\; i,j \ge 0\}$ for $M(p,q)$ 
 described in \S\ref{sect.good.basis}.   Let 
$$
\fa \; = \; \bigcap_{0 \le i,j<n}\Ann(e_{ij}).
$$ 

(1)
We will show that $\fa_{pq} \;=\;\fa$. 
 
Since $Ae_{0n} \cong Ae_{n0} \cong M(p,q)(-n),$ there are degree-$n$ $A$-module isomomorphisms  
$\psi_1,\psi_2 : M(p,q)\to M(p,q)$ such that $\psi_1(e_{00})=e_{n0}$ and $\psi_2(e_{00})=e_{0n}.$ 
Hence, for all $i,j \ge 0$,  $\psi_1(e_{ij})$ and $\psi_2(e_{ij})$ are non-zero scalar multiples of $e_{i+n,j}$  and $e_{i,j+n}$, respectively. 
To see this we argue by induction on $i,j.$ We will only do this for $\psi_1$. Suppose 
$\psi_1(\CC e_{ij}) = \CC e_{i+n,j}.$ Then 
\begin{align*}
 \psi_1(\CC e_{i+1,j}   \, \oplus   \,  \CC  e_{i,j+1})  &  \;=\; \psi_1(A_1e_{ij})
 \\
 & \;=\;  A_1\psi_1(e_{ij}) 
 \\
 & \;=\;  A_1e_{i+n,j} 
 \\
 & \; = \;  \CC e_{i+n+1,j} \, \oplus \, \CC e_{i+n,j+1}.
\end{align*}
Furthermore, ignoring the gradings, 
$$
Ae_{i+1,j}  \; \cong \; M(p+(j-i-1)\t,q+(i+1-j)\t)  \;\cong \; Ae_{i+n+1,j}
$$
and
$$
Ae_{i,j+1}   \; \cong \;  M(p+(j+1-i)\t,q+(i-j-1)\t)  \; \cong \;  Ae_{i+n,j+1}.
$$
Since $2\t \ne 0$, the two line modules just displayed are not isomorphic. 
It follows that $\psi_1(\CC e_{i+1,j}) = \CC e_{i+n+1,j}$ and $\psi_1(\CC e_{i,j+1})=\CC e_{i,j+n+1}.$

Thus $M(p,q)$ is a free $\CC[\psi_1,\psi_2]$-module with basis $\{e_{ij} \; | \;0 \le i,j \le n-1  \}.$

If $e_{k\ell}$ is one of the basis elements for $M(p,q)$, then $e_{k\ell}=\theta(e_{ij})$  for some  $0\le i,j\le  n-1$ 
and some $\theta \in \CC[\psi_1,\psi_2]$ so $\fa  e_{k\ell}= \fa  \theta(e_{ij}) =  \theta(\fa e_{ij})=0$.
Hence $\fa  M(p,q)=0$.  This completes the proof that $\fa_{pq} \;=\;\fa$. 

Since $Ae_{ij}$ is a line module, $d(A/ \! \Ann(e_{ij}))=2$ and $e\bigl( A/ \! \Ann(e_{ij}) \bigr)=1$.
Hence $d(A/\fa )=2$ and $e(A/\fa )\le n^2.$ 
This completes the proof of (1) and (2). 

(3) 
Since $(\O_1,\O_2) M(p,q)\ne 0$, $\Omega M(p,q)\ne 0$ for some 
$\Omega \in \CC \O_1 \oplus \CC \O_2$. The center of $A/ \fa_{pq}$ is  therefore strictly larger 
than $\CC$ so, by \cite[Cor. 2]{SZ14},  $A/ \fa_{pq}$ satisfies a polynomial identity.
\end{pf}

\begin{lemma} 
\label{lem3.2}
Assume  $p-q \notin 2\ZZ\tau.$ If 
$$
{\rm Hom}_{\Gr(A)}(M(p^\prime,q^\prime)(-m),M(p,q)) \ne 0,
$$
 then 
\begin{enumerate}
\item{} $\{ p^\prime,q^\prime\}=
\{p+(d-2k)\tau,q-d\tau\}$ for some $d,k \in \ZZ$ with $0 \le k \le m$;
\item{} if $2k\tau \ne 0,$ then $p+q \in E[2] + (k-1)\tau$.
\end{enumerate}
\end{lemma}
\begin{pf}
(1)
 Because line modules are 2-critical, non-zero homomorphisms between them are
injective.
Let  $0\ne e \in M(p,q)_m$ be such that $Ae \cong
M(p',q')(-m)$. Using the basis for $M(p,q)$ described in \S\ref{sect.good.basis}, write
$$
e \;=\;  \alpha_ie_{ij} \, +\,  \alpha_{i+1} e_{i+1,j-1}  \, +\,  \cdots  \, +\,  \alpha_{i+k} e_{i+k,j-k}
$$ 
with $\alpha_i \alpha_{i+k} \ne 0.$
Let $u,v \in A_1$ be such that $\ell_{p^\prime q^\prime}$ is the line $\{u=v=0\}$. 
Then $ue=ve=0$ so the component of $ue_{ij}$ in $\CC e_{i,j+1}$
is zero, and the component of $ue_{i+k,j-k}$ in $\CC e_{i+k+1,j-k}$ is
also zero. By \S\ref{sect.good.basis}, it follows that 
$$
u(q+(i-j)\tau) \;=\;  u(p+(j-k-i-k)\tau)\;=\; 0.
$$
The same argument applies to $v$ so $\{p^\prime,q^\prime\}=
\{p+(d-2k)\tau,q-d\tau\}$ where $d=j-i$. 

(2)
Since $M(p',q')$ embeds in $M(p,q)$,
$\Omega(p+q)=\Omega(p^\prime+q^\prime).$
It follows that either $p+q=p^\prime+q^\prime$, in which case 
$2k\tau=0$, or $p^\prime+q^\prime =
 -(p+q+2\tau),$ in which case  $p+q \in E[2] + (k-1)\tau$.
\end{pf}

\begin{corollary} 
\label{cor3.3}
Assume  $\tau$ has order $n<\infty$. Suppose  $p-q \notin 2\ZZ\tau$. % and let $\fp=\Ann(M(p,q))$.
\begin{enumerate}
\item{} 
There are only  finitely many $(p',q') \in E \times E$ such that $\fa_{p'q'}=\fa_{pq}$.
\item{}
 If $p+q \notin E[2] + \ZZ\tau$, then $\fa_{p'q'}=\fa_{pq} $ if and only if $(p',q')$ or $(q',p')$ is in 
 $\{(p+d\tau,q-d\tau) \; | \; d\in \ZZ\} \subseteq E \times E$; there are $n$ elements in this set.
\end{enumerate}
\end{corollary}
\begin{pf}
(1)
Suppose  $\fa_{p'q'}=\fa_{pq}$.
Set $R=A/\fa_{pq}$. Both $M(p,q)$ and $M(p^\prime,q^\prime)$ are critical 
$R$-modules so, by Proposition \ref{prop1.1},  there is a non-zero homomorphism $M(p' , q')(-m)
\to M(p,q)$ for some $m \in \ZZ$.  Lemma \ref{lem3.2} therefore applies and says that $p'$ and $q'$ belong to the finite  set
$\{p+(d-2k)\tau,q-d\tau \; | \; d,k \in \ZZ\}$.

(2)
Suppose $p+q \notin E[2] +\ZZ\tau$. Then $2k\tau=0$ so the only line
modules having the same annihilator as $M(p,q)$ are the line modules
$M(p+d\tau,q-d\tau)$ for $d\in \ZZ$. Since a suitable shift of each of these line modules embeds in each of the other ones,
they all have the same annihilator.  
\end{pf}

\begin{proposition}
\label{prop3.4}
Assume  $\tau$ has order $n<\infty$.
If  $p-q \notin 2\ZZ\tau$, then $A/ \fa_{pq}$ has
infinitely many shift-inequivalent 1-critical graded modules of 
multiplicity $\le n$.
\end{proposition}
\begin{pf}
We adopt the notation used in Proposition \ref{prop3.1}.

Let $R$ be the subalgebra $\CC [\psi_1,\psi_2]$ of $\End_A(M(p,q))$ generated by $\psi_1$ and $\psi_2$. 

Consider $M(p,q)$ as a left $\CC [\psi_1,\psi_2]$-module. Since $\psi_1 \psi_2(e_{00})$ and $\psi_2\psi_1(e_{00})$ are non-zero scalar
multiples of $e_{nn}$,  $\psi_2\psi_1 = \nu \psi_1\psi_2$ for some $\nu \in \CC^\times$. Thus, 
$R$ is a quotient of $\CC\langle x,y\rangle/(yx-\nu xy)$. 
Since the set $\{\psi_1^k \psi_2^{\ell}(e_{00}) \; | \; k,\ell \ge 0\}$ is linearly independent,
$\{\psi_1^k \psi_2^{\ell } \; | \; k,\ell \ge 0\}$ is a linearly independent subset of $\End_A(M(p,q))$.
Hence $R \cong \CC\langle x,y\rangle/(yx-\nu xy)$. (One can show that $\nu=1$ but we won't need that fact.) 

For each $\lambda=(\lambda_1,\lambda_2)\in \PP^1$ define $\psi_\lambda
=\lambda_1\psi_1 + \lambda_2\psi_2$ and $N^\lambda=M(p,q)/
{\rm im}(\psi_\lambda).$ Then $d(N^\lambda)=1$ and $e(N^\lambda)=n.$

It follows from well-known properties of  $\CC\langle x,y\rangle/(yx-\nu xy)$ that $\bigcap_{\lambda \in \PP^1}
R\psi_\lambda =0$. It follows from this and the fact that $M(p,q)$ is a free $R$-module  that $\bigcap_\lambda {\rm im}(\psi_\lambda) =0$.
Hence if $a \in \bigcap_\lambda \Ann(N^\lambda)$ then $a.e_{ij} \in
{\rm im}(\psi_\lambda)$ for all $i,j$ and all $\lambda$, so $a.M(p,q)=0.$
That is
$$ 
 \fa_{pq}\;=\; \bigcap_{\lambda \in \PP^1} \Ann(N^\lambda).$$

Fix $\lambda$. Then $N^\lambda$ has a finite filtration $N^\lambda=
N^0 \supseteq N^1 \supseteq \ldots \supseteq N^k =0$ such that each factor
$C^j:=N^{j-1}/N^j$ is 1-critical, and $k=k(\lambda) \le n$ since 
$e(N^\lambda)=n.$ 
Set $\fp^{\lambda,j}=\Ann(C^j)$. Then the $\fp^{\lambda,j}$ are  
prime ideals containing $\Ann(M(p,q))$ that satisfy 
$\fp^{\lambda,k}\ldots \fp^{\lambda,1}.N^\lambda =0.$ Thus 
\begin{equation}
\label{finite.intersection}
\bigcap_{\lambda \in \PP^1}  ( \fp^{\lambda,k(\lambda)}\ldots \fp^{\lambda,1}) \; \subseteq \;  \fa_{pq}
\end{equation}

Suppose the proposition is false. Then, up to 
shift-equivalence, there are only a finite number of 1-critical modules $C^j$
that occur above. Hence by Proposition \ref{prop.1.2}, for each $j$ there are only finitely many different $\fp^{\lambda,j}$'s as $\l$ varies over $\PP^1$. 
Therefore only finitely many different products $\fp^{\lambda,k(\lambda)}\ldots \fp^{\lambda,1}$ arise; 
thus the intersection in (\ref{finite.intersection}) is finite. It follows that 
$\Ann(M(p,q))$ contains the intersection of a finite number of ideals $I$ with the property that  $d(A/I)=1$. 
But this is absurd since $d(A/  \fa_{pq})=2$, so we conclude that the proposition is true. 
\end{pf}

\begin{proposition}
\label{lem3.5}
Assume  $\tau$ has order $n<\infty$.
If $p-q \notin 2\ZZ\tau$, then $A/ \fa_{pq}$ satisfies the identities of $2n \times 2n$ matrices.
\end{proposition}
\begin{pf}
By Proposition \ref{prop3.4}, $A/ \fa_{pq}$ has infinitely many finitely generated
shift-inequivalent 1-critical modules $C$ such that $e(C) \le n.$
If all but a finite number of these $C$ were point modules, then it would 
follow that $ \fa_{pq}\subseteq \bigcap_{x}\Ann(M(x))$ 
where the intersection is taken over infinitely many $x\in E.$ 
But such an intersection is equal to $(\Omega_1, \Omega_2)$.\footnote{This uses the fact that, first,  $A/(\O_1,\O_2)$ is a 
twisted homogeneous coordinate ring of $E$ and, consequently, that $\QGr(A/(\O_1,\O_2))$ is equivalent to $\Qcoh(E)$
and, second, that the intersection in $\cO_E$ of the ideals $\fm_x$ vanishing at any infinite set of points $x \in E$ is zero.}
all point modules are annihilated by
$(\O_1,\O_2)$ and that $A/(\O_1,\O_2)$ 
It would follow that $\fa_{pq} \subseteq (\Omega_1, \Omega_2)$ and 
hence that there is equality because
both are prime ideals of $A$ with quotient rings of GK-dimension 2. But 
this contradicts \cite[6.1]{LS8}which says  there is a unique 1-dimensional
subspace of $\CC \Omega_1 \oplus \CC \Omega_2$ that annihilates $M(p,q).$
Hence at
most a finite number of these 1-critical modules $C$ can be shifted
point modules. Thus there exist infinitely many shift-inequivalent
$C$ with $1<e(C) \le n.$ 

Suppose $C$ is a finitely generated 1-critical graded $A/ \fa_{pq}$-module with $1<e(C) \le n$. 
It follows from Artin and Van den Bergh's equivalence of categories (\ref{eq.QGr.B}) that the only finitely generated 1-critical 
modules are the point modules so there is an element 
$\Omega \in \CC \O_1 \oplus \CC \O_2$ such that $\Omega.C \ne 0.$ Since
$A/ \! \Ann(C)$ is a polynomial identity ring, $C$ has a non-trivial finite 
dimensional simple quotient, $S$ say. 
However, $S$ is annihilated by $\Omega-\nu$ for some $\nu \in \CC,$ so $S$
is actually a quotient of $C/(\Omega - \nu)C.$ An argument like the one in the proof of Proposition \ref{2.4} shows that $C$ is a 
free $\CC[\Omega]$-module;  since $\dim_\CC(C_i)=e$ for $e \gg 0$, $C$ is a free $\CC[\Omega]$-module of rank $2e(C)$;
hence $\dim(S) \le 2e(C).$
By Theorem \ref{2.5}, the (non-trivial) simple $A/ \! \Ann(C)$-modules are all twists of $S$
so also have dimension $\le 2n.$ Since $\Ann(C)$ is the intersection of the
annihilators of these simple modules, it follows that $A/ \! \Ann(C)$ 
satisfies the identities of $2n \times 2n$ matrices. 

Let $I=\bigcap \Ann(C)$ where the intersection is taken over the
infinitely many shift-inequivalent 1-critical modules $C$ with $1<e(C) \le n$ that
are $A/ \fa_{pq}$-modules. If $I \ne  \fa_{pq}$ then $d(A/I) \le 1$.
By passing to some (graded) minimal prime over $I$ we obtain a prime noetherian
ring of GK-dimension 1 having infinitely many shift-inequivalent 1-critical modules. 
This contradicts Theorem 1.3, so we conclude that $ \fa_{pq} = \bigcap \Ann(C).$
Hence, by the last paragraph, it follows that $A/ \fa_{pq}$ also
satisfies the identities of $2n \times 2n$ matrices. 
\end{pf}

\begin{theorem}
\label{thm3.6} 
If $\tau$ has order $n<\infty$,   then $A(E,\tau)$ satisfies a polynomial identity: it satisfies the identities of $2n \times 2n$ matrices.
\end{theorem}
\begin{pf}
 \underbar{Claim}: $0=\bigcap ( \Omega(z) )$ 
where the intersection is taken over all $z \in E$ such that
$z \notin E[2] +\ZZ\tau.$ \underbar{Proof}: 
There are infinitely many points $z_1,\,z_2,\,\ldots \in E-(E[2] +\ZZ\tau)$ such that the lines 
  $\CC\Omega(z_1), \CC \Omega(z_2),\ldots $ are pairwise distinct. 
  Using the fact that each $A/( \Omega(z_j))$ is a domain, an induction on $k$ shows that 
$$
\bigcap_{1 \le j \le k} ( \Omega(z_j) )\; = \; (\Omega(z_1),\ldots , \Omega(z_k)).
$$ 
In
particular, 
$$
\bigcap_{1 \le j \le k}(\Omega(z_j)) \; \subseteq  \;( A_{2k} ). 
$$ 
The intersection of the ideals $(A_{2k})$ is obviously zero, so the claim follows.
 
Now fix $z \in E$ such that $z \notin E[2] +\ZZ\tau.$
 
\underbar{Claim}:  $( \Omega(z) )=\bigcap \Ann \, M(p,q)$
 where the intersection is taken over all $p,q \in E$
such that $p+q=z$ and $p-q \notin 2\ZZ\tau.$ 
\underbar{Proof}: By [8, 6.2], $A/( \Omega(z) )$ is a domain
 of GK-dimension 3. By Proposition \ref{prop3.1}, $d(A/ \! \Ann(M(p,q)))=2$ for each such $p$ and 
$q$, so $Ann\,M(p,q)$ gives a height 1 prime in $A/(\Omega(z) ).$ By Corollary \ref{cor3.3}, $\{ \Ann \, M(p,q) \; | \;p+q=z\}$ is 
infinite. Since  the intersection of infinitely many height
one primes in a prime ring is zero, the truth of the claim follows. 

It follows from the two claims that $0= \bigcap \Ann(M(p,q))$ where the 
intersection is taken over all $p,q \in E$
such that $p+q\notin E[2] + \ZZ\tau$ and $p-q \notin 
2\ZZ\tau.$ 
However, for each such pair $(p,q),\,$ Proposition \ref{lem3.5} proves that $A/ \! \Ann(M(p,q))$ 
satisfies the identities of $2n \times 2n$ matrices. The result follows.
\end{pf}

 \subsubsection{Remark}
We now improve Proposition \ref{prop3.1} by showing that $A/  \fa_{pq}$
satisfies a polynomial identity for all $p,q \in E.$ Hence, by Proposition \ref{prop.2.7}(3),  $A/ \fa_{pq}$ has infinitely many shift-inequivalent
1-critical graded modules: in the terminology of \S4 and, after Lemma \ref{lem4.1}, it
follows that $\ell_{pq}$ contains infinitely many different fat points.

In the next result, $K_0(\ZZ)$ denotes the Grothendieck group of finitely generated projective $A$-moduiles

\begin{theorem} 
[J.T. Stafford] \cite{JTS17} 
\label{thm3.7}
Let $R$ be an Auslander-regular ring satisfying the Cohen-Macaulay 
property. If $K_0(R)=\ZZ$ then $R$ is a maximal order in its ring of 
fractions.
\end{theorem}

\begin{theorem}
\label{thm3.8}
If  $\tau$ has order $n<\infty$, then $A(E,\tau)$  
\begin{enumerate}
\item{} is a maximal order;
\item{} satisfies the polynomial identities of $2n \times 2n$ matrices;
\item{} is a finite module over its center.
\end{enumerate}
\end{theorem}
\begin{pf}
By  \cite{L7}, $A$ is Auslander-regular.
Since $A$ is a connected graded noetherian ring of finite global dimension, $K_0(A)=\ZZ$. 
The other hypotheses of Theorem \ref{thm3.7} 
were verified in \cite{LS8}, so Theorem \ref{thm3.7} applies to give (1). We have already 
proved (2).  
By \cite[13.9.8]{MR9}, (1) and (2) imply that  $A$ is its own trace 
ring. But the trace ring is always finitely generated as a module over 
its center so (3) holds.
\end{pf}

\subsubsection{} 
By \cite[Lem. 5.8c]{LS8}, if the order of $\tau$ is $n<\infty$ and $p\in E$, then $M(p)$ has a simple quotient module having dimension $n$. 
Hence  $A$ does not satisfy the identities of $m \times m$ matrices for any $m< n$.
Eventually, in Theorem \ref{thm.7.8}, we will improve Theorem \ref{thm3.8} and show that $A$ satisfies the identities of $n \times n$ matrices.

 \subsection{}
Recall that $n$ denotes the order of $\tau$ and $s$ the order of $2\tau$.

Since $A$ is a domain its ring of fractions, $\Fract(A)$, is a division ring.
By Theorem \ref{thm3.8}, $\Fract(A)$ is finite over its center and that center has transcendence degree 4 over $\CC$.
 It would be interesting to understand this division algebra. 
 We are unable to do this at present. 

We are able to say quite a lot
about another division algebra associated to $A$, namely
$$
D \; :=\;   (\Fract A)_0 \;=\; 
\{ ab^{-1} \; | \;a,b \in A \hbox { are homogeneous of the same
degree} \}. 
$$
The division ring $D$ is finite over its center which is of 
transcendence degree 3 over $\CC.$ In Theorem \ref{thm6.7} it is proved that $D$ is
of degree $s$ (= the order of $2\t$) over its center. We now  construct
 a sheaf of algebras that will be used in studying $D$.

 \subsubsection{}
Before doing that we note that $D$ is a more natural object than $\Fract(A)$ from the point of view of non-commutative algebraic geometry. The non-commutative projective variety $\Proj_{nc}(A)$ is an integral scheme in the sense of \cite{SPS-integral} so has a  division ring that 
plays the role of its function  field. That division ring is $D$.

  \subsection{}

Let $Z(A)$ denote the center of $A$ and define
$$
S \;=\; \Proj (Z(A)).
$$
Since $A$ is a finitely generated $Z(A)$-module the Artin-Tate Lemma tells us that  $Z(A)$ is a finitely generated (commutative graded) 
$\CC$-algebra. Since $A$ is a domain of GK-dimension 4, $Z(A)$ is an integral domain of Krull dimension 4. Hence $S$  
is an irreducible projective variety of dimension 3. We will now construct a sheaf $\cA$ of coherent ${\cal O}_S$-algebras.

 If $0\ne u \in Z(A)_d$ we write $Z(A)_{(u)}$ for the degree 0 subalgebra of $Z(A)[u^{-1}]$.
 We define $\cA$ by declaring that the sections of ${\cal A}$ over $S_{(u)}=\Spec(Z(A)_{(u)})$ is the degree
zero subalgebra of $A[u^{-1}]$, which is denoted by $A_{(u)}$ or
$\G(S_{(u)},\cA).$  

Thus $\cA$ is an order in $D$.
The center of $\cA$ may be strictly larger than $\cO_S$ and we denote that center by $\cZ.$ We denote by $Z$ the scheme
{\bf Spec}$_{\cO_S}\cZ$ as defined in \cite[Chap. II, Ex. 5.17]{H5}.  There is a finite morphism $Z \to S$. 
The inclusion of the polynomial ring $\CC[\O_1,\O_2] \subseteq Z(A)$  induces a rational map $S \dashrightarrow \PP^1$.

Although we do not have an explicit description of the center of $A$ we 
will study $\cA$ on $S$ and on $Z$ by looking at the fibers of the maps $Z \to S \dashrightarrow  \PP^1.$  
First, $S$ is the union, over all $\o \in E[2]$ and $0 \le k <n$, of the 
closed subschemes 
$$
X_{\o+k\t} \; :=\;  \Proj\left( \frac{Z(A)}{(\O(\o+k\t))}\right)
$$
and the dense open set that is their complement, namely 
$$
S_{(c)} \;=\;  \Spec Z(A)[c^{-1}]_0
$$ 
where 
\begin{equation}
\label{c}
c=\prod_{\omega \in E[2]} \prod_{k=0}^{\vert 2\tau \vert -1} \Omega(\omega+k\tau).
\end{equation}
Theorem \ref{thm6.7} shows that $\L_0:=\G(S_{(c)},\cA)=
A[c^{-1}]_0$ is Azumaya of rank $s^2$ over its center. This is
analogous to \cite[Thm. 7.3]{ATV2}.

Let $B=A/( \O_1, \O_2 )$. By \cite{AV}, $\Projnc(B)$ is isomorphic to $E$ and, by \cite{SPS_maps}, is 
a closed subscheme of $\Projnc (A)$ \ 

Clearly $\Projnc( B)$ is contained in each of the closed subschemes
$X_{\o+k\t}$, so along each of these closed sets we must concentrate our
attention on the complement of $E$, and hence on the rings
$A/ ( \O(\o+k\t) )[\O^{-1}]_0$ where $\O$ is chosen so that
$\CC \O_1 \oplus \CC \O_2 = \CC \O \oplus\CC\O(\o+k\t)$. The simple modules over
$A/ ( \O(\o+k\t) )[\O^{-1}]_0$ 
are in bijection with the $\O$-torsion-free fat point modules annihilated by
$\O(\o+k\tau)$. These are precisely the fat
point modules annihilated by $\O(\o+k\t)$.

This motivates the next three sections in which we classify the fat points.

\section{Generalities on Fat Points}
\label{sect.fat.pts}
%  Section 4. 

From now on suppose that $\tau$ is of finite order, $n$ say, and let $s$ denote the order of $2\t.$ Thus, $s$ is the smallest positive integer such that $s\t \in 
E[2].$

\subsection{Normalized modules}
\cite[p. 361]{ATV1}
%{\bf tidy this up}
\label{sect.normalized}
A finitely generated graded $A$-module $M$ is {\sf normalized} if its Hilbert series is equal to $e(1-t)^{-1}$ for some 
integer $e>0$ 
%it has a periodic Hilbert series, and 
and is Cohen-Macaulay (equivalently, it has zero socle).
For example, a point module is normalized. 

We call $e$ the {\sf multiplicity} of $M$ and denote it by $e(M)$.

%The GK-dimension of a normalized module is 1. The Hilbert series of a module of 
%GK-dimension 1 is equal to $q(t)(1-t)^{-1}$ for some $q \in \ZZ [t,t^{-1}]$. 
%Thus, for the Sklyanin algebra
%the periodicity requirement is equivalent to the condition that the Hilbert series be constant. 
%Thus a normalized $A$-module is a Cohen-Macaulay module 
%with Hilbert series is equal to $et^k(1-t)^{-1}$ for some $e\in \NN$ and some integer $k$.

My thanks to T. Levasseur for showing me a proof of the next result and for allowing me to include it here.

\begin{lemma}
\label{lem4.1} 
Every normalized $A$-module is generated by its degree-zero component. 
\end{lemma}
\begin{pf}
Let $M$ be a normalized $A$-module 
%Since the degree-shift is an auto-equivalence of $\Gr(A)$,
%it suffices to prove the result  when $M_0\ne0$ and $M_{<0}=0$. We therefore make that assumption. 
%Thus $H_M(t)=e(1-t)^{-1}$ for some $e \in \NN.$ 
and $P_\bullet \to M \to 0$ a minimal projective resolution  for it. Write $P_i=\oplus_j
A(-j)^{c_{ij}}$. Since the resolution is minimal, and $M_k=0$ for
$k<0$, $c_{ij}=0$ for all $j<0$ and 
$c_{i0}=c_{i1}=\ldots = c_{i,i-1}=0$. Since $M$ is Cohen-Macaulay, the remark after \cite[Lem. 1.12]{LS8} shows 
that $P_\bullet^\vee$ is a projective resolution of $M^\vee$. Thus, since the projective dimension of $M$ is 3, 
$c_{3j}$ is the number of elements of degree $-j$ in a minimal generating set for $M^\vee.$ By \cite[Prop. 1.10]{LS8}, 
$H_{M^\vee}(t)=et^{-3}(1-t)^{-1}$ so $c_{3j}=0$ for $j \ge 4.$ Hence
$c_{33}$ is the only non-zero $c_{3j}$. Thus $M^\vee$ is generated
in degree $-3$.

By \cite[Prop. 1.10]{LS8}, there is a duality $N \rightsquigarrow N^\vee$ between left and right
Cohen-Macaulay modules of a fixed projective dimension. Hence $M^{\vee\vee} \cong M$. 
It follows that $M$ is generated by $M_0$.
\end{pf}

\subsection{Fat points and fat point modules}

Let $\cO_p$ be a simple object in $\QGr(A)$. 

There is a 1-critical module $M \in \gr(A)$ such that $\cO_p \cong \pi^*M$.
If $M' \in \gr(A)$ also has the property that $\cO_p \cong \pi^*M'$, then $M_{\ge n} \cong M_{\ge n}$ for $n \gg 0$ so $e(M)=e(M')$.
We therefore define $e(\cO_p):=e(M)$ and call it the {\sf multiplicity} of $\cO_p$.
We say $\cO_p$ is a {\sf point} if $e(\cO_p)=1$ and a {\sf fat point} if $e(\cO_p)>1$.

We will later see that 
$$
e(\cO_p) \;= \; \dim_\CC\! \big( H^0(\Projnc(A),\cO_p)\big),
$$ 
which is, by definition, $\dim\big( {\rm Hom}_{\QGr(A)}(\cO,\cO_p)\big)$.

A 1-critical normalized $A$-module $F$ of multiplicity $>1$ is called a {\sf fat point module.} 
As an object in $\QGr(A)$, $F$ is irreducible, or simple. We often call the isomorphism class of $\pi^*F$ a {\sf fat point} and say that  {\sf  it lies on the (secant) line} $\ell$ if there is a non-zero homomorphism $M(\ell ) \to F$ in $\Gr(A)$.  

The main result in this section is that every fat point lies on a secant line. If we write $\cO_\ell$ for $\pi^* M(\ell)$ we can state this in 
familiar geometric terms: if $\cO_p$ is a point or a fat point, there is an epimorphism $\cO_\ell \to \cO_p$ in $\QGr(A)$ for some secant line
$\ell$.

Theorem \ref {thm4.2} shows if $\cO_p$ is a fat point, then $\cO_p \cong \pi^*F$ for some fat point module $F$. 

We also express this by saying that $\ell$ {\it contains} the fat point $F$. 
More generally, we say that {\it a 1-critical module $C$ 
lies on the secant line $\ell$} if the fat point determined by $C$ lies on
$\ell.$

By \cite[Thms. 5.3, 5.8]{LS8},   every point module is a quotient of a line module.

 \subsection{}
The next result can be proved as in  \cite[Prop. 6.6]{ATV2} by using the  Auslander-Gorenstein property of $A$ (see \cite{L7} and \cite{LS8}).

\begin{theorem}
\label{thm4.2} 
{\rm (cf. \cite[Prop. 6.6]{ATV2})}
Let $C$ and $D$ be graded $A$-modules with
$d(C)=d(D)=1.$ We write $C^\vee:=\Ext^3_A(C,A).$ 
\begin{enumerate}
\item{} 
If  $f \in {\rm Hom}_{\Gr(A)}(C,D)$ becomes an isomorphism in $\QGr(A)$, then the dual map $f^\vee : D^\vee \to C^\vee$ 
provides
an equivalence between $D^\vee$ and $C^\vee.$
\item{}
 If $C$ is Cohen-Macaulay, then $C$ is contained in an equivalent
module $D$ such that for some $k \ge 0, \, D(k)$ is normalized;
 that is, $D$ is a negative shift of a normalized module.
\item{} 
If $C$ is equivalent to a normalized module $F$, then there is a
unique homomorphism $C_{\ge 0} \to F$ that extends this equivalence.
In particular, if $C=C_{\ge 0}$ is 1-critical, then $C$ embeds in an
equivalent normalized module. 
\item{} 
A 1-critical module in $\gr(A)$ is equivalent to a unique normalized module.
In particular, every fat point is represented by a unique fat point module.
\item{}
 If $C=C_{\ge 0}$ is 1-critical, and $F$ is either a point module or a fat point module equivalent to $C$, then 
$\dim_\CC( {\rm Hom}_A(C,F)_0)=1.$
\end{enumerate}
 \end{theorem}

\begin{theorem}
\label{thm4.3} 
Every fat point lies on a secant line.
\end{theorem}
\begin{pf}
Let $D$ be a fat point module of multiplicity $e$. We begin by replacing $D$ 
by another fat point module. By Theorem \ref {thm4.2},  $D(-1)$ embeds in a unique fat point 
module, $F$ say.  Let  $\fp=\Ann(F)$. Since $F$ 
is critical, $\fp$ is a prime ideal.
If $0\ne m \in F$, then $\CC[\Omega_1,\Omega_2].m$ has GK-dimension 1,
so $\fp \cap \CC[\Omega_1,\Omega_2] \ne 0.$ This intersection is a homogeneous prime 
ideal of $\CC[\Omega_1,\Omega_2]$ so contains a non-zero $ \Omega
\in \CC\Omega_1 \oplus \CC\Omega_2$. 
By construction,  $\Omega.F=0.$

By \cite[Prop. 6.12]{LS8}, there is a secant line $\ell$  such that $\Omega.M(\ell)=0$. 
Let $\ell^\perp \subseteq A_1$ be the subspace vanishing on $\ell$.
The homogeneous function 
$$
\ell^\perp \; \longrightarrow \;  {\rm Hom}_{\CC}(F_0,F_1) \; \cong  \; {\rm End}_{\CC}(\CC^2) \; \stackrel{\rm det}{\longrightarrow} \; \CC
$$
 has a non-trivial zero so $u.m=0$ for some $0 \ne u \in \ell^\perp$ and some $0 \ne m \in F_0$.  But  $\Omega m=0$ also, so there is a non-zero homomorphism $\varphi : A/Au+A \Omega \to F$. 
Consider the diagram
$$
\xymatrix{ 
0 \ar[r] & L\ar[r] & A/Au+A \Omega  \ar[d]^{\varphi}  \ar[r] &  M(\ell) \ar[r] & 0
\\
&& F
}
$$
in which the row is exact.

By \cite[6.2]{LS8},   $A/(\Omega)$ is a domain. Since $A$ is a domain
$H_{A/( \Omega )} (t) = (1-t^2)(1-t)^{-4},$ and 
$H_{A/Au+A \Omega}(t) = (1-t)(1-t^2)(1-t)^{-4}.$ It follows that $H_L(t)=
t(1-t)^{-2}.$ But $L$ is cyclic, generated by the image of $\ell^\perp/\CC u$, so 
$L(1)$ is a line module.  Write $M(\ell^\prime):=L(1).$ 

There are two possibilities. 

(1) 
If $L \subseteq \ker(\varphi)$, then there is an induced map 
${\overline \varphi}:M(\ell) \to F$ that is non-zero. In this case the fat 
point represented by $F$ lies on $\ell.$ 
Let $p\in E \cap \ell$, and  consider the diagram
$$
\xymatrix{
0  \ar[r] & M(\ell^{\prime \prime})(-1)  \ar[r] & \ar[d] M(\ell) \ar[r] & M(p) \ar[r] & 0
\\
0 \ar[r] & D(-1) \ar[r] & F
}
$$
If the image of  
$ M(\ell^{\prime\prime})(-1)$ in $F$ is non-zero, then the fat point 
represented by $D$ lies on $\ell^{\prime\prime}.$ However, if that image is
zero, then $F \cong M(p)$ and consequently $D \cong M(p-\t).$ In this
case $D$ also lies on a secant line.

(2)
If $L \not\subseteq \ker(\varphi),$ then $\varphi(L) \ne 0$ so  there is a non-zero map $M(\ell^\prime)(-1) \to F$ in $\QGr(A)$, 
the image of which is in $F_{\ge 1}=D(-1).$ Hence there is a non-zero map $M(\ell^\prime) \to D,$
so the fat point  corresponding to $D$ lies on $\ell^\prime.$
\end{pf}

\begin{proposition}
\label{prop4.4}
Let $F$ be a fat point module of multiplicity $e(F)=k>1$.
\begin{enumerate}
\item{} 
$F_{\ge 2}(2) \cong F$ and $F_{\ge 2s}(2s) \cong F$.\footnote{If $M$ is a point module, then 
$M(p)_{\ge 2s}(2s) \cong M(p)$ because $M(p)_{\ge j}(j) \cong M(p-j\t)$ and $2s\t=0$, .}
\item{} 
For all $j\ge 0$, the inclusion 
$$
M(p+j\tau,q-j\tau)(-j) \; \longrightarrow \;  M(p,q)
$$
 induces an isomorphism
$$
 {\rm Hom}_{\Gr(A)}(M(p,q),F)  \; \longrightarrow \;   {\rm Hom}_{\Gr(A)}(M(p+j\tau,q-j\tau),F(j)).
 $$ 
\item{} 
Let $C$ be a quotient of a line module. Suppose $C$
is equivalent to $F$. Then $F$ lies on $\ell_{pq}$
if and only if ${\rm Hom}_{\Gr(A)}(M(p \pm (k-1)\tau,q \mp (k-1)\tau),C(k-1)) \ne 0.$
\item{} 
For every $i\in\ZZ$ and every $p,q\in E$, the lines $\ell_{pq}$
and $\ell_{p+2i\t,q-2i\t}$ contain exactly the same fat points of 
multiplicity $>1$.
\end{enumerate}
\end{proposition}

 \subsubsection{Remark}
 \label{remark}
By \cite[\S5]{LS8}, 
${\rm Hom}_{\Gr(A)}(M(p+\t,q-\t)(-1),M(p,q)) \cong \CC$ so
there is a unique chain of submodules
$$
M(p,q)  \, \supseteq \, M(p+\t,q-\t)(-1)  \, \supseteq \, \cdots \, \supseteq \, M(p+j\t,q-j\t)(-j) .
$$
This is the inclusion referred to in part (2).

\begin{pf}
(1) 
Since the only 1-critical modules over $B=A/( \O_1,\O_2 )$ are its
point modules \cite{AV}, there is a non-zero central element $\Omega \in A_2$ 
such that $\O.F \ne 0$. Multiplication by $\Omega$ is an isomorphism  $F \cong F_{\ge 2}(2).$ It follows by
induction that  $F_{\ge 2s}(2s) \cong F$.

(2) The exact sequence $0 \to M(p+\tau,q-\tau)(-1) \to M(p,q)
\to M(p) \to 0$ gives rise to the exact sequence 
$$
0 \to {\rm Hom}_{\Gr(A)}(M(p),F) \to {\rm Hom}_{\Gr(A)}(M(p,q),F) \to {\rm Hom}_{\Gr(A)}(M(p+\tau,q-\tau)(-1),F).
$$ 
Since $F$ is 1-critical of multiplicity $>1$, ${\rm Hom}_{\Gr(A)}(M(p),F)=0$. 
Repeating this argument by using the inclusions remarked on in \S\ref{remark}, it follows that 
there  are injective maps 
$$
{\rm Hom}_{\Gr(A)}(M(p,q),F)  \;\longrightarrow \;  {\rm Hom}_{\Gr(A)}(M(p+j\tau,q-j\tau)(-j),F)
$$
 for all $j\ge 0.$ In particular, this is true for $j=2s$. However, by (1) the composition
\begin{align*}
{\rm Hom}_{\Gr(A)}(M(p,q),F) & \; \longrightarrow \;  {\rm Hom}_{\Gr(A)}(M(p+2s\tau,q-2s\tau)(-2s),F)
\\
& \;  \stackrel{\sim}{\longrightarrow} \;   {\rm Hom}_{\Gr(A)}(M(p,q)(-2s),F) 
 \\
&\phantom{i} \;=\;  \phantom{i.} {\rm Hom}_{\Gr(A)}(M(p,q),F(2s)) 
 \\
&\phantom{i} \;=\;  \phantom{i.}   {\rm Hom}_{\Gr(A)}(M(p,q),F(2s)_{\ge 0})
\\
&\phantom{i} \;=\;  \phantom{i.}  {\rm Hom}_{\Gr(A)}(M(p,q),F)
\end{align*}
is injective and therefore an isomorphism because all these Hom-spaces have finite dimension.  Hence for all $0\le j \le 2s$ the inclusions
$$
{\rm Hom}_{\Gr(A)}(M(p,q),F) \;  \subseteq \;  {\rm Hom}_{\Gr(A)}(M(p+j\tau,q-j\tau)(-j),F)  \; \subseteq \; {\rm Hom}_{\Gr(A)}(M(p,q)(-2s),F)
$$ 
are equalities. The result is identical for $j \ge 2s.$

(3) 
By Theorem \ref{thm4.3}, such a $C$ exists. That is, there is an inclusion $C \subseteq F.$ 
Since $C$ is a multiplicity-$k$ 1-critical quotient of a line module, it 
follows from \cite[Prop. 4.4]{SS16} that $H_C(t)=(1-t^k)(1-t)^{-2}$, and in particular
that $\dim_\CC(C_j)=k$ for all $j \ge k-1.$ Therefore, the inclusion $C \subseteq F$
actually gives an equality $C_j=F_j$ for all $j \ge k-1$ since $F$ is 
normalized of multiplicity $k$. Thus $F(k-1)_{\ge 0} \cong C(k-1)_{\ge 0}.$
By (2), $F$ lies on $\ell_{pq}$ if and only if 
${\rm Hom}_{\Gr(A)}(M(p+(k-1)\tau,q-(k-1)\tau),F(k-1)) \ne 0$; the result now follows.

(4)
By (1) and (2), ${\rm Hom}_{\Gr(A)}(M(p,q),F) \cong {\rm Hom}_{\Gr(A)}(M(p+2\t,q-2\t),F)$. Therefore 
$F$ lies on  $\ell_{pq}$ if and only if  it lies on $\ell_{p+2\t,q-2\t}.$ The result follows by an induction 
argument.
\end{pf}

Proposition \ref{prop4.4}(3) will be used frequently to  decide which lines a particular fat point lies on. 
However, to make it an effective tool we
need more information about critical factors of line modules. The 
following result gives a rather good description of this situation.

\begin{theorem}
\label{thm4.5} 
\cite[Prop. 4.4]{SS16}
Let $C$ be a 1-critical quotient of $M(p,q)$ of multiplicity  $k+1>1$.
\begin{enumerate}
\item{}
In $\QGr(A)$, there is an exact sequence 
$$
0 \to M(p-(k+1)\t,q-(k+1)\t)(-k-1) \to M(p,q) \to C \to 0.
$$ 
\item{} 
If $2(k+1)\t \ne 0$, then $p+q=\o+k\t$ for some $\o \in E[2].$
\end{enumerate}
\end{theorem}

 Theorem \ref{thm4.5} makes it clear that one must understand
${\rm Hom}_{\Gr(A)}(M(p-k\t,q-k\t)(-k),M(p,q))$ in order to determine which fat
points of multiplicity $k$ lie on $\ell_{p \pm k\t, q\mp k\t}.$
We solve this problem the next two sections.

\section{Fat point modules of intermediate multiplicity}
\label{sect.F.intermediate.e}

A fat point module $F$ has {\sf intermediate multiplicity} if $1<e(F)<s.$ 

The phrase {\sf fat point of intermediate multiplicity} will refer to the image of such a module in $\QGr(A)$. 

In \S\ref{sect.fatties.mult.s}, we will prove that the multiplicity of a fat point is at most $s$, and for most fat points it is $s$. 
In this section, Theorem \ref{thm5.4} classifies the fat points of intermediate multiplicity as objects in $\QGr(A)$.
In particular, we show there are only finitely many of  them up to isomorphism.

\subsection{}
The fat points of intermediate multiplicity are closely related to certain 
finite dimensional modules defined by Sklyanin \cite{Skly12} that were shown to be 
simple in \cite{SS16}.
For each $\o \in E[2]$ and each $k \in \NN$ there is a finite dimensional module $V(\o+k\t)$ defined in \cite[\S3]{SS16}. It is shown 
there that if $k<s$, then $V(\o+k\t)$ is a simple $A$-module of 
dimension $k+1$. The module $V(k\tau)$ may be realized as a subspace of the
ring of meromorphic functions on $\CC$, and the action of $A$ on this
realization is given explicitly in the proof of Proposition \ref{prop5.1} below (see also
\cite{SS16}).

For $k \in \NN$ and $\o \in E[2]$, define 
$$
F(\o+k\t) \; := \;  \widetilde V(\o+k\t).
$$
We will show that $F(\o+k\t)$ is a fat point of multiplicity $k+1$ if $0\le k \le s-1$ and, passing to $\QGr(A)$, that 
these are the only fat points of intermediate multiplicity, and that $F(\o+k\t)$ lies on $\ell_{pq}$ if and only if
$p+q=\o+k\t.$

\subsubsection{} 
The modules $F(\o)$ for $\o \in E[2]$ are point modules. In fact $F(\o)$ is one of the 4 special point modules $M(e_i)$ corresponding to
the four points $e_0,\ldots,e_3$ that are the vertices of the singular quadrics containing $E$.  It will be 
convenient to refer to these as fat points of intermediate multiplicity  also. 
They behave rather differently from the point modules $M(p)$ parametrized by the points $p \in E.$

\subsubsection{Notation} 
In this and later sections $\, \o$ will always denote  a point in $E[2]$  
and $k$ will always denote a non-negative integer.

\begin{proposition}
\label{prop5.1} 
Assume $0 \le k \le s-1.$
Then ${\rm Hom}_A(M(p,q),V(\o+k\t)) \ne 0$ if and only if
$p+q=\o+k\t.$
\end{proposition}
\begin{pf}
$(\Rightarrow)$ 
Since $V(\o+k\t)$ is obtained from $V(k\t)$ by twisting with an automorphism, 
as explained in \cite[\S3]{SS16}, it is enough to prove the result when $\o=0.$ 
Suppose that $p+q =k\t$, and write   $V=V(k\t).$
As in \cite[\S3]{SS16}, there is $0 \ne f \in V$ such
that $X.f=0$ for all $X \in A_1$ for which $X(\ell_{pq})=0.$
Thus
$$
(X.f)(z)  \; = \; {{X(z-\hbox{${{1}\over{2}}$}k\t)} \over {\theta_{11}(2z)}} f(z+\t) \, -\, 
{{X(-z-\hbox{${{1}\over{2}}$}k\t)} \over {\theta_{11}(2z)}} f(z-\t) 
$$
is zero for all $z \in \CC$ if $X(p)=X(q)=0.$ Choose $r\in E$ such that 
$$
-r-k\t \notin \{p,q,r\}, \qquad f(r+(\hbox{${{1} \over {2}}$}k - 1)\t) \ne 0, \quad \hbox{and} \quad 2r+k\t \ne 0.
$$
Let $X \in A_1$ be such that its divisor of zeroes on $E$ is
$(p)+(q)+(r)+(-p-q-r).$ Then 
$$
0 \;=\;  (X.f)(r+\hbox{${{1}\over{2}}$}k\t)   \;=\;    - \,
{{X(-r-k\t)} \over {\theta_{11}(2r+k\t)}} f(r+(\hbox{${{1}\over{2}}$}k-1)\t),
$$
so $X(-r-k\t)=0$ also.  Thus $p+q=k\t$ as required.

$(\Leftarrow)$ 
This was proved in \cite[\S3]{SS16}.
\end{pf}

\begin{lemma}
\label{lem5.2} 
Suppose $1 \le k \le s.$ Then $V(\o+k\t)$ is not a quotient of 
any 1-critical graded module of multiplicity $\le k$.
\end{lemma}
\begin{pf}
Suppose the result is false, and that $V(\o+k\t)$ is a quotient of a
1-critical module $C$ with $e(C)=d+1\le k.$ First we show that $C$ cannnot 
be a point module. 

Suppose to the contrary that $C$ is a point module, say $C=M(\xi)$.
Then $\xi \notin\{e_0,e_1,e_2,e_3\}$ since all simple quotients of the point modules $M(e_i)$ are 1-dimensional
and $\dim_\CC(V(\o+k\t))=k+1>1.$ Thus $\xi
\in E.$ Hence by \cite[Lem. 5.8a]{LS8}, $\xi-(k+1)\t=\xi$ whence $k=s-1$.
Since $M(\xi)$ contains a copy of $M(\xi-i\t)(-i)$ for all $i \in \NN$, it
follows that $V(\o+k\t)$ is also a quotient of $M(\xi-i\t)$ for all $0\le
i<s.$ By Proposition \ref{prop.old.2.2}(2), there is a degree 0 injection $M(\xi-i\t) \to \widetilde
V(\o+k\t)$. This gives $s$ inequivalent submodules of $\widetilde V(\o+k\t)$,
so each $M(\xi-i\t)$ is a \lq composition factor\rq\ of $\widetilde 
V(\o +k\t)$, and since $e(\widetilde V(\o+k\t))=k+1=s$, these are all the
composition factors of $\widetilde V(\o+k\t)$ as an object of $\QGr(A)$.

Now choose $p,q\in E$ such that $p+q=\o+k\t$ and $\{p,q\} \cap \{\xi-i\t
\; | \;i \in \NN\} = \varnothing.$ By Proposition \ref{prop5.1}, $V(\o+k\t)$ is a quotient of
$M(p,q)$ and hence by Proposition \ref{prop.old.2.2}(2) $\widetilde V(\o+k\t)$ contains a non-zero
submodule that is a quotient of $M(p,q)$. Since $\widetilde V(\o+k\t)$ has no
finite dimensional non-zero submodule, this image of $M(p,q)$ is of
GK-dimension 1. Hence $\widetilde V(\o+k\t)$ has a subquotient that is a
1-critical factor of $M(p,q)$. This 1-critical quotient of $M(p,q)$ is
therefore a \lq composition factor\rq\ of $\widetilde V(\o+k\t)$, and hence is
equivalent to $M(\xi-i\t)$ for some $i\in\NN.$ By [8, \S2] a 1-critical
module of multiplicity 1, is necessarily a point module, so the only 1-critical
quotients of $M(p,q)$ of multiplicity 1 are the point modules
$M(p)$ and $M(q)$. So either $M(p)$ or $M(q)$ is equivalent to (and hence
isomorphic to) some $M(\xi-i\t)$. Thus either $p=\xi-i\t$ or $q=\xi-i\t$
for some $i\in\NN.$ But this possibility is excluded by the choice of $p$
and $q$. We conclude that $C$ is not a point module.

By Theorem \ref{thm4.5},  it follows that either $2(d+1)\t=0$ or $p+q=\o' +d\t$ for some
$\o' \in E[2]$. The first possibility cannot occur since $d+1\le k<s.$
The second possibility cannot occur since, if it did then $2(k-d)\t=0$
which is impossible since $0<d<k<s.$ Hence no such $C$ can exist.
\end{pf}

\begin{lemma}
%{\bf Lemma 5.3.}
\label{lem.5.3}
Assume $0\le d\le s-1$ and  ${\rm Hom}_{\Gr(A)}(M(p',q'),M(p,q)(d)) \ne 0$.
If $p+q \notin E[2]+i\t$ for any $0\le i \le d$, then
$M(p',q') \cong M(p+j\t,q-j\t)$ for some $j \in \ZZ$.
\end{lemma}
\begin{pf}
By induction on $d$. There is nothing to do if $d=0$ so suppose 
$d\ge 1.$ Let $L$ be a submodule of $M(p,q)$ with $L \cong 
M(p',q')(-d)$.
Let $L' \subseteq M(p,q)$ be maximal subject to the condition that
$L' \supseteq L$ and $M(p,q)/L'$ is not finite dimensional. 
Then $C:=M(p,q)/L'$ is a 1-critical quotient of $M(p,q)/L.$ 
Since $e(C)\le d$ it follows from Theorem \ref{thm4.5} that either
$C \cong M(p)$ or $C \cong M(q)$. Hence by \cite[Thm. 5.5]{LS8} either
${\rm Hom}_{\Gr(A)}(M(p',q'),M(p+\t,q-\t)(d-1))$ or
${\rm Hom}_{\Gr(A)}(M(p',q'),M(p-\t,q+\t)(d-1))$ is non-zero. The hypotheses of the 
lemma apply to these two cases with a smaller $d$ so the induction argument proceeds.
\end{pf}

\begin{theorem}	
%{\bf Theorem 5.4.} 
\label{thm5.4}
$\phantom{x}$
\begin{enumerate}
\item{} 
$F(\o+k\t)$ is a fat point of multiplicity $k+1$ if  $\o \in E[2]$ and $1 \le k \le s-1$.
\item{} $F(\o+k\t)$ lies on $\ell_{pq}$ if and only if $p+q=\o+k\t.$
\item{} $\{ F(\omega+k\tau) \; | \; \omega \in E[2], \, 1 \le k\le s-1\}$ are all
the fat points of intermediate multiplicity.
\end{enumerate}
\end{theorem}
\begin{pf}
(1)
 Since $e(F(\o+k\t))=\dim( \widetilde V(\o+k\t) )= k+1$ it remains to 
prove that $F(\o+k\t)$ is critical.
If $k=0$, then $F(\o)$ is the point module $M(e_i)$ for the appropriate $e_i.$

Suppose $k>1$ and choose $p,q \in E$ such that $p+q=\o+k\t$.

Because $V(\o+k\t)$ is a quotient of $M(p,q)$ it follows from Proposition \ref{prop.old.2.2}(2) that there is a non-zero map $\varphi: M(p,q) \to F(\o+k\t)$ in $\Gr(A)$.
Let $C$ be the image of $\varphi$. To show that $F(\o+k\t)$ is critical
it suffices to prove that $C$ is critical of multiplicity $k+1$. If $C$ is not critical, or is of 
smaller multiplicity, then $F(\o+k\t)$ has a filtration by submodules such that the slices are 
either finite dimensional or of multiplicity $ \le k.$ However, $V(\o+k\t)$ would be  a quotient of one of those factors and that 
would contradict Lemma \ref{lem5.2}. 

(2) 
As remarked in (1), if $p+q=\o+k\t$, then 
${\rm Hom}(M(p,q),F(\o+k\t)) \ne 0$,  so $F(\o+k\t)$ lies on $\ell_{pq}$.
Conversely, suppose there is a non-zero map $\varphi:M(p,q) \to F(\o+k\t)$ and let $C$ be the image of $\varphi.$ 
Since $V(\o+k\t)$ is a quotient of $F(\o+k\t)$ it is a quotient 
of either $C$ or $F(\o+k\t)/C.$ Since $F(\o+k\t)$ is critical, the 
latter module is finite dimensional and, because it is a graded module, its only
composition factor is the trivial module. Hence $V(\o+k\t)$ is a 
quotient of $C$ and hence of $M(p,q).$ Thus $p+q=\o+k\t$ by Proposition \ref{prop5.1}.

(3)
Suppose  the result is false, and let $F$ be another 
 fat point of multiplicity $k+1$.  Then $F$ lies on some line by Theorem \ref{thm4.3}.
By Theorem \ref{thm4.5}  that line must be of the form $\ell_{pq}$ with $p+q=\o+k\t$ 
for some $\o \in E[2].$ But $F(\o+k\t)$ also lies on that line, so it 
follows from Theorem \ref{thm4.5} that $ \dim {\rm Hom}(M(p-(k+1)\t,q-(k+1)\t)(-k-1),M(p,q)) \ge 2.$ 
Let $L=Ae$ and $L'=Ae'$ be distinct submodules of $M(p,q)$  both 
isomorphic to $M(p-(k+1)\t,q-(k+1)\t)(-k-1).$ 

Consider the obvious map $\varphi : M(p,q) \to M(p)$. Then $\varphi$
vanishes on some linear combination of $e$ and $e'$ so after replacing
$e$ by this element, we may assume that $\varphi(e)=0$. Thus
$L \subseteq \ker(\varphi) \cong M(p+\t,q-\t)(-1),$ so
${\rm Hom}(M(p-(k+1)\t,q-(k+1)\t),M(p+\t,q-\t)[k]) \ne 0.$ It follows from Lemma \ref{lem.5.3} that
$M(p-(k+1)\t,q-(k+1)\t) \cong M(p+j\t,q-j\t)$ for some $j$. Thus
$(p-(k+1)\t)+(q-(k+1)\t) = p+q $ which forces $2(k+1)\t=0$. But this is patently
false so the result is true.
\end{pf}

\begin{corollary}
\label{cor5.5}
Assume $p+q = \o+k\t$ and $0 \le k \le s-1$. Up to scalar multiple there is a unique non-zero map
$$ 
\varphi:M(p-(k+1)\t,q-(k+1)\t)(-k-1) \to M(p,q)
$$
in $\Gr(A)$. Furthermore, ${\rm coker}(\varphi)$ is equivalent to $F(\o+k\t)$ and the latter has 
a $(k+1)$-dimensional simple quotient, namely $V(\o+k\t)$.
\end{corollary}
\begin{pf}
This follows from Theorem \ref{thm4.5}, the proof of Theorem \ref{thm5.4}.
\end{pf}

\begin{proposition}
\label{prop5.6} 
Assume $1 \le k \le s-1$.
\begin{enumerate}
\item{} 
If $p+q=\o+k\t$,  then $F(\o+k\t)$ is the 
only fat point of intermediate multiplicity that lies on $\ell_{pq}.$
\item{}  
The only fat points of intermediate multiplicity annihilated by $\O(\o+k\t)$ are $F(\o+k\t)$ and $F((\o+s\t)+(s-k-2)\t)$. 
These fat points are of multiplicity $k+1$ and $s-(k+1)$ respectively.
\item{}
Let $\o\in E[2]$. If  $p+q=\o+k\t$, then there is an exact sequence $$0\to X\to Y \to Z \to 0$$
in which 
$$X\sim F((\o+s\t)+(s-k-2)\t)(-k-1), \quad Z \sim F(\o+k\t),$$
and $Y$ is a  multiplicity-$s$ quotient of $M(p,q)$.
\end{enumerate}
\end{proposition}
\begin{pf}
 (1) 
 If $\o +k\t=\o^\prime +k^\prime\t$ with $\o,\o' \in E[2]$ 
and $0<k,k'<s-1$ then it is easy to prove that $\o=\o'$ and $k=k'.$
Hence the result follows from Theorem \ref{thm5.4}.

(2) 
Suppose $\o+k\t \ne \o' +k'\t.$
The fat points $F(\o+k\t)$ and $F(\o'+k'\t)$ are both annihilated by 
$\O(\o+k\t)$ if and only if  $\o+k\t+\o'+k'\t = -2\t.$ 
The result follows from a simple calculation.

(3) 
Write $\o'=\o+s\t$ and $k'=s-k-2$. By Theorem \ref{thm5.4}(2) there is a quotient,
$Z$ say, of $M(p,q)$ which is equivalent to $F(\o+k\t)$. In fact there is
an exact sequence $0 \to M(p',q')(-k-1) \to M(p,q) \to Z \to 0$
where $p'=p-(k+1)\t$ and $q'=q-(k+1)\t.$ Since $p'+q'
=\o'+k'\t$ there is a similar  exact sequence $0 \to
M(p'',q'')(-k-1-k'-1) \to M(p',q')(-k-1) \to X \to 0$ with
$p''=p-s\t, q''=q-s\t$ and $X \sim F(\o'+k'\t)(-k-1).$
Putting these two  exact sequences together in the obvious way gives
the result. 
\end{pf}

\begin{proposition}
\label{prop5.7}
Assume $0\le k \le s-1$. If $p+q=\o+k\t$, then $\ell_{pq}$ and $\ell_{p,q-2(k+1)\t}$
contain exactly the same multiplicity-$s$ fat points.
\end{proposition}
\begin{pf}
There is nothing to prove if $k=s-1$ so suppose $k \le s-2.$

By symmetry it suffices to show that every multiplicity-$s$ fat point  
lying on $\ell_{pq}$  lies on $\ell_{p,q-2(k+1)}.$ Let $F$ be such a
fat point and let $C$ be a 1-critical quotient of $M(p,q)$ such that $C
\sim F.$

By Proposition \ref{prop4.4}(3),  to show $F$ lies on $\ell_{p,q-2(k+1)\t}$ it suffices to prove
that
$$
{\rm Hom}_A(M(p-(s-1)\t,q+(s-2k-3)\t), C(s-1)) \ne 0.
$$
However, 
$$
 M(p-(s-k-2)\t,q+(s-k-2)\t)  \;\subseteq \; M(p,q)(s-k-2).
 $$
 Applying Corollary \ref{cor5.5} to this submodule, we have
$$
{\rm Hom}_A(M(p-(s-1)\t,q+(s-2k-3)\t), M(p,q)(s-1)) \ne 0.
$$
Since $e(F)=s$, $e(C)=s$ too;  it follows from Theorem \ref{thm4.5}  that
the map $M(p,q) \to C$ is an isomorphism in degree $s-1$. The result
follows.
\end{pf}

\section{Multiplicity-$s$ fat points}
\label{sect.fatties.mult.s}

We prove the following results in this section.
\begin{enumerate}
\item{} 
There is a central element $c$, specified in Theorem \ref{thm6.7},
such that $A[c^{-1}]_0$ is an Azumaya algebra of rank $s^2$, whence
$\Fract(\cA)$ is a division algebra of degree $s$.
\item{} 
The irreducible objects in $\QGr(A)$ consist of the
point modules, the fat points of intermediate multiplicity, and the multiplicity-$s$ fat points  (that is, $e(C)\le s$ for all 1-critical
modules $C$).
\item{} 
Every secant line contains infinitely many  multiplicity-$s$ fat points.
\end{enumerate}

  \subsection{}

Recall that $L(z) = \{\ell_{pq} \; | \; p+q =z\}$.

We  begin by analysing the line modules $M(p,q)$ such that $p-q \in 2\ZZ\t$ and the multiplicity-$s$ 
fat points that lie on them. 
Proposition \ref{prop.6.3} proves  that such an $M(p,q)$ has a 1-parameter family of 1-critical quotient modules $C^\l$, $\l\in\PP^1$.
Theorem \ref{thm6.8} shows that infinitely many  of these $C^\l$'s are pairwise inequivalent. 
This is done by showing that such a  $C^\l$ lies on a line in  $L(-p-q-2\tau)$ and that $C^\l$ and $C^\mu$ lie on
different lines  in $L(-p-q-2\tau)$ if $\l\ne\mu.$

When $p-q \not\in 2\ZZ\t$ we use the action of $A_1$ on the basis $e_{ij}$ for $M(p,q)$ that was described in \S\ref{sect.good.basis}.

\begin{lemma}
%{\bf Lemma 6.1.} 
\label{lem.6.1}
Assume $p-q \notin 2\ZZ\t.$ Let $e=\sum_{i=0}^m \alpha_i e_{i,m-i} \in M(p,q)_m$ and let $X\in A_1$.
\begin{enumerate}
\item{}  
If $z=p+q+2\tau$ and $p_k=p+(m-2k+2)\tau$, then 
\begin{align*}
\qquad 
\phantom{xx}  X.e  & \; =\;
 \beta_0 X(z-p_0)e_{0,m+1} \quad + \quad
 \\ 
& \quad \quad \sum_{k=1}^m \bigl( \beta_k X(z-p_k) \,-\, \beta_{k-1} X(p_k) \bigr)
e_{k,m+1-k} 
\; - \; \beta_m X(p_{m+1}) e_{m+1,0}
\end{align*}
where   
$$
\beta_k \; = \;  \frac{\alpha_k}{\theta_{11}(p-q-2(2i-m)\tau}\, .
$$
\item{} 
If $f=\l_0e_{0,m+1}\, + \, \l_1 e_{m+1,0}$,  % \in M(p,q)_{m+1}$
 then $X.e \in \CC f$ if and only if 
\begin{equation}
\label{eq.dag}  
\begin{pmatrix}
\lambda_1 X(z-p_0) & 0 &   \ldots&& 0&  \lambda_0 X(p_{m+1}) \cr
-X(p_1) & X(z-p_1) &  \ldots &  & 0 & 0 \cr
0 & -X(p_2) &   \ldots &  & 0 & 0 \cr
\vdots   & \vdots  &   & & \vdots & \vdots \cr
\vdots  & \vdots & &   & \vdots & \vdots \cr
0 & 0&  \ldots &  & -X(p_m) & X(z-p_m) \cr
\end{pmatrix} 
\begin{pmatrix}
\beta_0 \cr
\beta_1 \cr
\vdots \cr
\vdots \cr
\vdots \cr
\beta_m \cr
\end{pmatrix}
=0.
\end{equation}
\end{enumerate}
\end{lemma}
\begin{pf}
(1) 
The expression for $X.e$ is a straightforward calculation.

(2)
It follows from (1) that $X.e$ is in $\CC f$ if and only if $\b_k X(z-p_k) - \b_{k-1}
X(p_k)=0$ for all $1 \le k \le m$ and $(\b_0 X(z-p_0), -\b_m X(p_{m+1}))
= \mu  (\l_0, \l_1)$ for some $\mu \in \CC$.  Since $(\l_0,\l_1) \ne (0,0)$, the last condition holds  if and only if 
$$
\det
\begin{pmatrix} 
\lambda_0 & \lambda_1   \cr
\b_0 X(z-p_0) & -\b_m X(p_{m+1})  \cr
\end{pmatrix}
 =0;
$$
this happens if and only if the first row of the matrix annihilates $(\b_0,\ldots,\b_m)^t$.
The other conditions hold if and only if the other rows of  the matrix annihilate $(\b_0,\ldots,\b_m)^t$.
\end{pf}

\begin{proposition}
%{\bf Proposition 6.2.}
\label{prop.6.2}
 If $p-q \notin 2\ZZ\tau$, then 
$$
\dim_\CC \big( {\rm Hom}(M(p-s\tau,q-s\tau)(-s),M(p,q))\big) \;=\; 2.
$$
This Hom-space is spanned by the two $A$-module homomorphisms that send the generator $e_{00}$
of $M(p-s\t,q-s\t)$ to $e_{0s}$ and $e_{s0}.$
\end{proposition}
\begin{pf}
Since $p-q \notin 2\ZZ\tau$  we will use the basis $\{e_{ij}\}$ described just before Lemma \ref{lem.6.1}. Since $2s\tau=0$,  
 $$
 Ae_{0s} \; \cong \; M(p+s\tau,q-s\tau)(-s) \;=\; M(p-s\tau,q+s\tau)(-s) \;  \cong \;  Ae_{s0}
 $$
 so the dimension of the Hom-space is at
least 2. It remains to show that linear combinations of $e_{0s}$ and 
$e_{s0}$ are the only elements in $M(p,q)_s$ that generate a submodule isomorphic to
$M(p-s\tau,q-s\tau)(-s)$.  

Let $u,v \in A_1$ be such that $\ell_{p-s\tau,q-s\tau}$ is the line $\{u=v=0\}$.

Take $n=s$ and $e =\sum_{i=0}^s \alpha_i e_{i,s-i} $ as in the previous lemma, and suppose  $X.e=0$ for all $X\in \CC u \oplus \CC v$.   Certainly $X.e_{0s}=X.e_{s0}=0$ so we may assume
that  $e =\sum_{i=1}^{s-1} \alpha_i e_{i,s-i} $ and we must show that
$\alpha_i=0$ for all $i.$ Suppose $\alpha_0=\alpha_1=\ldots=
\alpha_{i-1}=0$ for some $1\le i \le s-1$. 
The term $X.\alpha_i e_{i,s-i}$ has a component in
$\CC e_{i,s+1-i}$ which cannot cancel with any other term arising from
$X.e$; hence that component is zero. Thus $\alpha_i X(q+(2i-s)\tau)=0$ 
and this must be true for all $X\in \CC u \oplus \CC v.$ However, since
$q+(2i-s)\tau \notin \{p-s\tau,q-s\tau\}$ some $X \in \CC u \oplus \CC v$ does
not vanish at it. It follows that $\alpha_i=0$ and the induction 
argument proceeds. Eventually we conclude that $\alpha_i=0$ for all $i$, as required.
\end{pf}

\begin{proposition}
%{\bf Proposition 6.3.} 
\label{prop.6.3}
Assume $p-q \notin 2\ZZ\tau$.
\begin{enumerate}
\item{} 
For all except finitely many $\lambda=(\l_0,\l_1) \in \PP^1$,  the module
\begin{equation}
\label{C-lambda}
C^\lambda \; :=\; \frac{M(p,q)}{A(\l_0e_{0s} + \l_1 e_{s0})}
\end{equation}
is 1-critical of multiplicity $s$. 
\item{} 
 Every multiplicity-$s$ 1-critical quotient of $M(p,q)$  is equal to some $C^\lambda$.
\end{enumerate}
\end{proposition}
\begin{pf}
(1) 
Let $\l=(\l_0,\l_1) \in \PP^1$ and define $f_\l=\l_0 e_{0s} +\l_1 e_{s0}.$
Let $C^\l$ be the module in (\ref{C-lambda}). 
By Proposition \ref{prop.6.2}, $Af_\lambda \cong M(p-s\tau,q-s\tau)(-s)$  so the Hilbert series  of
$C^\lambda$ is $(1-t^s)(1-t)^{-2}$ and $e(C^\lambda)=s.$

Suppose, first, that $p+q \notin E[2] + \ZZ\tau.$ 

If $C^\lambda$ is not critical, then there is a graded submodule $L$ 
of $M(p,q)$, properly containing $Af_\lambda$  such that $N:= M(p,q)/L$ 
is 1-critical (pick $L \supseteq Af_\l$ maximal such that $M(p,q)/L$ is
infinite dimensional). By \cite[Prop. 4.4]{SS16} either $N$ is a point module, or
$H_N(t)=(1-t^d)(1-t)^{-2}$ for some $d$ such that  $2d\tau=0$
or $p+q \in E[2] + (d-1)\tau.$ 

Suppose $N$ is not a point module. Since $N$ is a proper quotient 
of $C^\lambda$, consideration of Hilbert series shows that $d<s,$ whence
 $2d\tau \ne 0.$ Therefore $p+q \in E[2] +(d-1)\tau.$
However, this possibility is excluded by hypothesis so we conclude that
$N$ is a point module.

By [8, \S 5] a line module has at most three submodules such that the 
quotient is a point module. In fact, since $p+q \notin E[2]$, 
$M(p,q)$ can have at most two submodules such that the quotient is a 
point module. Those submodules are isomorphic to $M(p+\tau,q-\tau)(-1)$ and
$M(p-\tau, q+\tau)(-1)$. Thus there are at most 2 possibilities 
for $L$. Moreover, since $L \cong M(p^\prime,q^\prime)(-1)$ is a shifted 
line module, and $p^\prime + q^\prime=p+q$, we may apply a similar 
argument to $L \supseteq Af_\lambda.$ In particular, since $e(L/Af_\lambda )<s$, $L/Af_\lambda$ 
 has a 1-critical quotient that is necessarily a shifted point module. 
	
Continuing in this way, we see there are only finitely many
$\lambda$ such that $C^\lambda$ is not critical, and that each of these
has a 1-critical \lq composition\rq \ series with all its composition factors being shifted point modules.

Suppose, now,  that $p+q=\o+(k-1)\t \in E[2]+\ZZ\t.$

It is now possible for $N$ to be not only a point module, but also a critical
module of multiplicity $k$. However, by Corollary \ref{cor5.5}, if $N$ is of multiplicity
$k$, then $L$ is uniquely determined. The proof then proceeds as before.

(2) By Theorem \ref{thm4.5}(1), if $C$ is a 1-critical graded quotient of $M(p,q)$ of 
multiplicity $s,$ there is an exact sequence $ 0 \to M(p-s\t,q-s\t)(-s) \to M(p,q) \to C \to 0 $.
After Proposition \ref{prop.6.2} it follows that the $C^\l$'s  are the only possible 
1-critical quotients of $M(p,q)$ of multiplicity $s$. 
\end{pf}

{\bf Remarks. 1.} Proposition \ref{prop.6.3} can be improved. 
If $p+q \notin E[2]+\ZZ\t$ only two   $C^\l$'s  fail 
to be critical, namely $M(p,q)/Ae_{0s}$ and $M(p,q)/Ae_{s0}$.
If $p+q \in E[2]+\ZZ\t$ only three of the $C^\l$ will fail 
to be critical, namely $M(p,q)/Ae_{0s}$ and $M(p,q)/Ae_{s0}$ and one other
that is a non-split extension in $\QGr(A)$ of a fat point module of intermediate multiplicity by 
a shifted fat point module of intermediate multiplicity.

{\bf 2.} To get a better feel for Proposition \ref{prop.6.3} recall the geometric 
situation. Certainly , every secant line $\ell_{pq}$ contains a 1-parameter family of
points; two of these points lie on $E$, and if $p+q \in E[2]$ then 
one of the points on $\ell_{pq}$ is the singular point $e_i$ of $Q(\o_i)$
for a suitable $i\in\{0,1,2,3\}.$ Thus Proposition \ref{prop.6.3} is a reasonable 
analogue of the geometric situation.

 \subsection{}

The next three results address the question of when two lines can contain a
common multiplicity-$s$ fat point. This problem is dealt with more
completely in \S\ref{sect.fat.pts.s}. First, observe that if $F$ is a fat point module,
then $F$ is not an $A/( \O_1,\O_2 )-$module, so there is a
unique $\O(z)$ (up to scalar multiple) that annihilates $F$. Hence if $\ell\in
L(z)$ and $\ell^\prime \in L(z^\prime)$ contain a common fat point then
either $z^\prime=z$ or $z^\prime=-z-2\tau.$ These two possibilities are 
addressed in Propositions \ref{prop.6.4} and  \ref{prop.6.5}. However, by Proposition \ref{prop4.4}(4),
$\ell_{pq}$ and $\ell_{p+2i\t,q-2i\t}$ contain exactly the same fat points
for all $i \in \ZZ.$
\medskip

{\bf Notation.} 
The next few results use the notation in the proof of Proposition \ref{prop.6.3}:  
for $\l= (\l_0,\l_1) \in \PP^1$ we define $f_\l:=\l_0 e_{0s} + \l_1 e_{s0}$ and $C^\l:=M(p,q)/Af_\l$.

\begin{proposition}
%{\bf Proposition 6.4.}
\label{prop.6.4}
Assume $p-q \notin 2\ZZ\t.$ 
Let $\{ C^\l \}$ be the multiplicity-$s$ 1-critical modules lying on $\ell_{pq}.$
Let $\ell_{p' q'} \in L(p+q)$ and suppose that $\ell_{p' q'} \ne 
\ell_{p+2i\t,q-2i\t}$ for any $i \in \ZZ.$ 
\begin{enumerate}
\item{} 
If $p+q=\o-\t$ for some $\o \in E[2]$, then there is a unique $C^\l$ lying on $\ell_{p' q'}$.
\item{} 
If $p+q \notin E[2]-\t$, then no $C^\l$ lies on $\ell_{p' q'}.$
\end{enumerate}
\end{proposition}
\begin{pf}
Choose $u,v \in A_1$ such that the line $\{u=v=0\}$ is $\ell_{p' +(s-1)\t, q'
-(s-1)\t}.$ By Proposition \ref{prop4.4}(3), $C^\l$ lies on $\ell_{p' q'}$ if and only if
$${\rm Hom}(M(p'+(s-1)\t, q' - (s-1)\t),C^\l(s-1)) \; \ne \; 0.$$ By Proposition \ref{prop.6.2},
this ${\rm Hom}$-space is non-zero if and only if there is a non-zero $e$ in exists $M(p,q)_{s-1}$ such that $u.e,v.e \in \CC f_\l.$ 
By Lemma \ref{lem.6.1}, this is equivalent
to the existence of a common solution $0\ne(\b_0,\ldots , \b_{s-1})^t$ to
the two matrix equations (\ref{eq.dag}) for $X=u$ and $X=v.$
In (\ref{eq.dag}),  $z=p+q+2\t,\, m=s-1,\, p_i =p+(s-2i+1)\t$, and $p_{m+1}=p_0.$
The determinant of the matrix in (\ref{eq.dag}) is
$$
\l_1 X(z-p_0)\ldots X(z-p_{s-1}) \; + \; \l_0 X(p_0) \ldots X(p_{s-1}).
$$

The hypothesis on $\ell_{p' q'}$implies that neither $p'+(s-1)\tau$
nor $q'-(s-1)\t$ is   in $\{z-p_i, p_i\; | \; i \in \ZZ \}$ so we can, and do, 
choose $u$ and $v$ such that $u(z-p_i) \ne 0$ and $u(p_i)\ne 0$ for all $i$,
and $v(p_0)=0$. Hence the
determinant obtained with $X=u$ is a non-zero form having a unique zero
$\l=(\l_0,\l_1)$, and this zero has the property that $\l_0\l_1 \ne 0.$
Furthermore, there is a unique solution $(\b_0,\ldots,\b_{s-1})^t$ to the
matrix equation with $X=u$, and $\b_i \ne 0$ for all $i$. 

Now consider (\ref{eq.dag}) with $X=v.$ Looking at the first row of the matrix,
the fact that $v(p_0)=0$ means that we can only have a common solution $(\b_0,\ldots,\b_{s-1})^t$ when  $v(z-p_0)=0$ also. However, since $p-q
\notin 2\ZZ\t$, it follows that $z-p_0 \notin \{p'+(s-1)\t, q' -(s-1)\t,
p_0\}$, so the only way $v$ could vanish at $z-p_0$ is if $z-p_0$
is the fourth zero of $v$.  Since the sum of the 4 zeroes of $v$ is
zero, its other zero is $-(2p+q+(s+1)\t).$ 
Therefore, if $z-p_0 \ne  -(2p+q+(s+1)\t),$ the matrix equations for
$X=u$ and $X=v$ cannot have a common solution, and we conclude that no
$C^\l$ can lie on $\ell_{p' q'}.$

A calculation shows that $z-p_0 =  -(2p+q+(s+1)\t)$ if and only if
$p+q \in E[2]-\t.$ 
Thus,  if $p+q \notin E[2] -\t$, then no $C^\l$ lies on $\ell_{p' q'},$ whence (2) is true.

(1)
Now suppose that $p+q=\o-\t$ with $\o \in E[2].$ Consider the rational function
${{u}\over{v}}$ on $E$. It has two zeroes,  at 
$r$ and $z-r$ say, and two poles, at $p_0$ and $z-p_0$. 
Since The rational functions ${{u}\over{v}}(x)$ and ${{u}\over{v}}(z-x)$ have the
same zeroes and poles so they are scalar multiples of each other; evaluating the functions at an  
$x$ such that $2x=z$, we see that this scalar multiple must be 1.
Therefore $\bigl( {{u} \over {v}} \bigr) (\xi) =
\bigl( {{u} \over {v}} \bigr) (z-\xi)$ for all $\xi \in E$. 
It follows that $(\b_0,\ldots,\b_{s-1})^t$ is also a solution to the matrix
equation (\ref{eq.dag}) with $X=v.$ This proves that
$${\rm Hom}_{\Gr(A)}\big(M(p'+(s-1)\t,q'-(s-1)\t),C^\l(s-1)\big)\; \ne \; 0$$ for the unique $\l$ 
described above. Thus to prove (1), it remains to check that this $C^\l$ is
critical. Since $p+q =\o-\t = (\o+s\t)+(s-1)\t$ with $\o+s\t \in E[2]$ it
follows from \S4 that $\ell_{pq}$ has no fat points of intermediate
multiplicity (since $e(F(\o-\t))=s$). Therefore the only $C^\mu$ that are not
critical are $C^{(0,1)}$ and $C^{(1,0)}.$ But $\l_0\l_1 \ne 0$ so this
particular $C^\l$ is critical.
\end{pf}

\begin{proposition}
%{\bf Proposition 6.5.} 
\label{prop.6.5}
Assume $p-q \notin 2\ZZ\tau$ and 
let $\ell_{ p^\prime   q^\prime} \in L(-p-q -2\t)$. 
Let $\{C^\l\}$ be the multiplicity-$s$ 1-critical modules that are quotients of  $M(p,q).$
\begin{enumerate}
\item{} If $p+q \notin E[2]+\ZZ\t$, then there is exactly one
$\lambda$ such that $C^\l$ lies on $\ell_{p^\prime q^\prime}$, except when
$\{p',q'\} \cap \{p+2j\t,q+2j\t \; | \; j \in \ZZ \} \ne \varnothing$.
\item{} If $p+q=\o+k\t$ with $\o \in E[2]$ and $0\le k \le s-2$ then
\begin{enumerate}
\item{} if $\ell_{p^\prime q^\prime}= \ell_{p-2(i-1)\t, q+2j\t}$ 
for some $i,j \in \ZZ$, then every $C^\l$ lies on $\ell_{p^\prime q^\prime}$;
\item{} if $\ell_{p^\prime q^\prime}\ne \ell_{p-2(i-1)\t, q+2j\t}$ 
for any $i,j \in \ZZ$, then no $C^\l$ lies on $\ell_{p^\prime q^\prime}$.
\end{enumerate}
\end{enumerate}
\end{proposition}

\smallskip

{\bf Remark.} 
The case $p+q=\o+(s-1)\t$ is not covered by Proposition \ref{prop.6.5}(2), but in
that case $L(p+q)=L(-p-q-2\t)$, so it is covered by Proposition \ref{prop.6.4}(1).
Part (2b) of Proposition \ref{prop.6.5} was proved in Proposition \ref{prop5.7} for all $p,q \in E$. 

\smallskip

\begin{pf}
Let $u,v \in A_1$ be such that the line $\{u=v=0\}$ is $\ell_{p^\prime+(s-1)\t, 
q^\prime-(s-1)\t} $. As in the proof of Proposition \ref{prop.6.4},
 $C^\l$ lies on $\ell_{p^\prime q^\prime}$ if 
and only if the two matrix equations (\ref{eq.dag}) for $X=u$ and $X=v$ have a simultaneous solution.
In (\ref{eq.dag}), $z=p+q+2\t,\, m=s-1,\, p_i =p+(s-2i+1)\t$ and 
$p_{m+1}=p_0.$
The determinant of the matrix in (\ref{eq.dag}) is
$$
\l_1 X(z-p_0)\ldots X(z-p_{s-1}) \; + \;
\l_0 X(p_0) \ldots X(p_{s-1}).
$$

Now choose $u$ and $v$ such that both have no zeroes in the set
$$
\{ p_k,z-p_k \; | \;0\le k \le s-1\} \,\,  \, - \, \,\,   \{
p^\prime+(s-1)\t, q^\prime-(s-1)\t \}.
$$

\underbar{Claim:} $(\b_0,\ldots,\b_{s-1})^t$ is a solution to (\ref{eq.dag}) for
$X=u$ if and only if it is a solution for $X=v.$ \underbar{Proof}: The
rational function ${{u}\over{v}}$ on $E$ has two zeroes, the
sum of which is $z$. As in the proof of Proposition \ref{prop.6.4}(1) it
follows that $\bigl( {{u} \over {v}} \bigr) (\xi) =
\bigl( {{u} \over {v}} \bigr) (z-\xi)$ for all $\xi \in E$. Thus
$$ 
{\rm rank}
\begin{pmatrix} 
-u(p_k) & u(z-p_k)  \cr
-v(p_k) & v(z-p_k)  \cr
\end{pmatrix}
\le 1
$$
for all $k.$ If the top row of this matrix is $(0 \; 0)$ then the careful
choice of $u$ and $v$ ensures that $\{p_k,z-p_k\} \subseteq \{p'+(s-1)\t,
q'+(s-1)\t\}$ so the bottom row of the matrix is also $(0 \; 0)$.
Therefore the condition rank $\le 1$ implies that the two rows are non-zero
multiples of each other. Similarly, the rows of 
$$
\begin{pmatrix} 
\lambda_1 u(z-p_0) & \lambda_0 u(p_0)  \cr
\lambda_1 v(z-p_0) & \lambda_0 v(p_0)  \cr
\end{pmatrix}
$$
are non-zero multiples of each another. The truth of the claim follows.
$\square$

(1)
 It follows from the careful choice of $u$ that 
the determinant of the matrix in (\ref{eq.dag}) for $X=u$ is a non-zero form: the only way it could not
be is if there were some $i,j$ such that $\{p_i,z-p_j\} \subseteq
\{p'+(s-1)\t,q'+(s-1)\t\}$ but this is excluded by the hypotheses on
$p+q$ and $p-q$. Hence there is a unique $\l$ for which the determinant
vanishes, and therefore a unique $\l$ such that
$${\rm Hom}_{\Gr(A)}\big(M(p'+(s-1)\t,q'-(s-1)\t),C^\l(s-1)\big)\; \ne\; 0.$$
Since $\ell_{pq}$ contains no fat points of intermediate multiplicity,
there are only two $C^\l$ that fail to be critical, namely $C^{(0,1)}$ and
$C^{(1,0)}.$ If $\l=(0,1)$ then  $u(z-p_i)=0$ for some $i$, so
$z-p_i \in \{p'+(s-1)\t,q'-(s-1)\t\}$ from which it follows that
$\{p',q'\} \cap \{p+2j\t,q+2j\t \; | \; j \in \ZZ \} \ne \varnothing$.
The case $\l=(1,0)$ is similar.

(2a)
The hypothesis on $\ell_{p' q'}$ may be rephrased as
$\ell_{p' +(s-1)\t, q'-(s-1)\t} = \ell_{p_i,z-p_j}$  whence
$u(p_i)=u(z-p_j)=0$. Therefore the determinant of the matrix in (\ref{eq.dag}) at $X=u$ is identically zero
for all $\l$. 

(2b)
The hypothesis on $\ell_{p' q'}$ ensures that for all $i$ and
$j$,  $\{p_i,z-p_j\}$ is not contained in $\{p'+(s-1)\t,q'+(s-1)\t\}$ 
so, as in case (1), there is a unique $\l$ such that
${\rm Hom}_{\Gr(A)}(M(p'+(s-1)\t,q'-(s-1)\t),C^\l(s-1)) \ne 0.$ However, in contrast
to case (1), $\ell_{pq}$  contains fat points of intermediate multiplicity,
and we will show that the $C^\l$ we have just found is not critical. 
To do this, we describe a particular $C^\mu$ that is not critical, and
show that ${\rm Hom}_{\Gr(A)}(M(p'+(s-1)\t,q'-(s-1)\t),C^\mu(s-1)) \ne 0,$  and it
then follows from the uniqueness of $\l$ that $\l=\mu.$

Set $F:=F(\o-(k+2)\t)$. By Proposition \ref{prop5.6}(3) and its proof, there is a $\mu \in \CC$ fitting into an exact sequence 
$$
0\to X \to C^\mu \to Z \to 0
$$ 
in which $X \sim F(-k-1).$ By Theorem \ref{thm5.4}(2), $F$ lies on every $\ell_{p_1 q_1} \in
L(-p-q-2\t)$. It therefore follows from Proposition \ref{prop4.4}(2) that ${\rm Hom}_{\Gr(A)}(M(p_1,q_1),F(j))
\ne 0$ for all such $\ell_{p_1 q_1}$ and all $j>0.$
Since 
$$C^\mu(s-1) \; \supseteq \; X(s-1)\, \sim \, F(s-k-2)$$
it follows that ${\rm Hom}_{\Gr(A)}(M(p_1,q_1), C^\mu(s-1)) \ne 0$ for all $\ell_{p_1 q_1}
\in L(-p-q-2\t).$ In particular ${\rm Hom}_{\Gr(A)}(M(p'+(s-1)\t, q'
-(s-1)\t),C^\mu(s-1)) \ne 0.$
\end{pf}

\begin{proposition}
%{\bf Proposition 6.6.} 
\label{prop.6.6}
Assume $p-q \notin 2\ZZ\tau$ and let $F$ be a multiplicity-$s$ fat point module lying on $\ell_{pq}.$ 
\begin{enumerate}
\item{}
$F$ lies on some line in $L(-p-q-2\t);$
\item{} 
If $p+q \notin E[2]+\ZZ \tau$,  then $F$ lies on at most $s$ lines in $L(-p-q-2\t)$.
\end{enumerate}
\end{proposition}
\begin{pf}
We must show the following (in the notation of Proposition \ref{prop.6.3}):
if $C^\l$ is critical, then $ {\rm Hom}_{\Gr(A)}(M(p^\prime, q^\prime),C^\lambda(s-1)) \ne 0$ for some
$\ell_{p^\prime q^\prime} \in L(-p-q-2\tau)$.
 The requirement that
$C^\l$ be critical implies that $\l \ne (0,1),(1,0).$

Let $\ell_{p^\prime q^\prime}$ be an arbitary line in
$L(-p-q-2\tau).$ Let  $u,v \in A_1$ be such that $\ell_{p^\prime q^\prime}=\{u=v=0\}.$ We must show that there is a 
choice of $p^\prime$ and $q^\prime$ for which there is a non-zero $e \in  M(p,q)_{s-1}$ 
such that both $u.e$ and $v.e$ are in 
$\CC f_\lambda.$ Let $z$, $p_j$, and $z-p_j$, be as in Proposition \ref{prop.6.3}.
There is an element $e\in M(p,q)_{s-1}$ with the required property if and
only if the matrix equation (\ref{eq.dag}) in Lemma \ref{lem.6.1} has a common non-zero solution
$(\beta_0, \ldots, \beta_{s-1})^t$ for both $X=u$ and $X=v.$

For  $0\le j \le s-1$ write $\ell_j=\ell_{p_j,z-p_j}$, fix
linear isomorphisms $\varphi_j:\PP^1 \to \ell_j$ with $\varphi_j(1,0)=p_j$
and $\varphi_j(0,1)=z-p_j$, and define maps $\Phi_j:L(-z) \to \PP^1$ by
$\Phi_j(\ell_{p^\prime q^\prime})=\varphi^{-1}_j \bigl(\ell_j \cap
\ell_{p^\prime q^\prime} \bigr)$. Each $\Phi_j$ is well-defined: first,
the points $p^\prime, \, q^\prime, \, p_j, \, z-p_j \,$ are coplanar which
ensures that $\ell_j$ and $\ell_{p' q'}$ do intersect; second, 
since $p+q \notin E[2] -\t$ the lines $\ell_j \in L(z)$ and
$\ell_{p' q'} \in L(-z)$ are distinct. 
Furthermore, $\Phi_j$ is a morphism of varieties. Recall that
$L(-z) \cong \PP^1$ is the fiber over $-z$ of the addition map $+:E\times E 
/\ZZ_2
\to E$. If $\Delta$ is the diagonal, then the map to the Grassmannian,
$(\PP^3 \times \PP^3) \, - \, \Delta \to \GG(2,4)$, that sends a pair of 
points to the line they span is a morphism. The map
$$ 
\GG(2,4) \; \supseteq \;    \{ \hbox {lines } \ell^\prime \hbox { in } \PP^3 \; | \;\ell^\prime  \cap \ell_j \ne \varnothing \} \, - \, \{\ell_j\} \to \ell_j 
\; \cong \; \PP^1
$$ 
given by $\ell^\prime \to \ell^\prime \cap \ell_j,$ is also a morphism. The map
$\Phi_j$ is  the restriction  to $L(-z)$ of the composition of these morphisms.
Clearly $\Phi_j$ is injective, so it is given by a linear form $L(-z)
\cong \PP^1 \to \PP^1.$

Now define $\Phi:L(-z) \to \PP^1$ as follows. Let $\ell \in L(-z).$
For each $j$ write  $\Phi_j(\ell)=(\gamma_j, -\delta_j)$ and define 

\centerline{ $
\Phi(\ell)=(\lambda_1 \gamma_0 \ldots \gamma_{s-1}, \lambda_0 \delta_0 
\ldots \delta_{s-1}).$}

\noindent
Since $z \notin E[2] + \ZZ\tau$, $\{p^\prime, q^\prime\} \ne \{p_i, z-p_j\}$ for all $i$ and $j$. 
This ensures that $\Phi(\ell)\ne (0,0).$  Thus $\Phi$ is a degree $s$ form on $L(-z).$  It is
clear that $\Phi$ is not constant, and therefore surjective, since 
$\Phi(\ell_{p_j,-z-p_j})= (1,0)$ and
$\Phi(\ell_{z-p_j,-2z+p_j})=(0,1).$ In particular, $(1,-1)$ is in the 
image. Now choose $\ell_{p^\prime q^\prime} \in L(-z)$ such that $\Phi(\ell_{p^\prime q^\prime})
=(1,-1);$ there are at most $s$ such lines. 

Since  $\{p^\prime, q^\prime\} \ne \{p_j, z-p_j\}$,  we can, and do,  choose $u\in A_1$ such that $u(p^\prime)=u(q^\prime)=0$ and 
$u(p_j) \ne 0 \ne u(z-p_j)$ for all $0\le j \le s-1$. Pick representatives
of $p^\prime, \, q^\prime, \, p_j, \, z-p_j$ in $A_1^*.$ These are linearly 
dependent and if $\Phi_j(\ell_{p^\prime q^\prime})=(\gamma_j, -\delta_j)$
then $\gamma_j u(p_j) -\delta_j u(z-p_j) =0.$ Hence, as points on $\PP^1$, $(u(p_j),u(z-p_j))=(\delta_j,\gamma_j).$ Therefore
\begin{align*}
& \big(\lambda_1 u(z-p_0)\ldots u(z-p_{s-1}),\, \lambda_0 u(p_0)\ldots u(p_{s-1})\big)
\\
& \phantom{xx} \;= \; \big(\lambda_1 \gamma_0 \ldots \gamma_{s-1}, \lambda_0 \delta_0 \ldots \delta_{s-1}\big) 
\\ 
& \phantom{xx} \;= \; (1,-1),
\end{align*}
 so the determinant of the matrix in (\ref{eq.dag}) with $X=u$ is zero. Let $(\beta_0,\ldots,\beta_{s-1})^t$ be a non-zero solution to
the matrix equation for $X=u$. As in the proof of Proposition \ref{prop.6.5}, $(\beta_0,\ldots,\beta_{s-1})^t$ is also a
solution for $X=v.$
\end{pf}

\begin{theorem}
\label{thm6.7}
Let 
$$
c\; =\; \prod_{\omega\in E[2]} \prod_{k=1}^{s-1} \Omega(\omega+k\tau).
$$
\begin{enumerate}
\item{} $A[c^{-1}]_0$ is an Azumaya algebra of rank $s^2$ over its center.
\item{} ${\rm Fract}(A[c^{-1}]_0)$ is a central simple division algebra 
of degree $s$.
\item{} If $F$ is a fat point module such that $cF \ne 0$, then $e(F) = s.$
\item{} If $F$ is a fat point module, then $e(F) \le s.$
\end{enumerate}
\end{theorem}
\begin{pf}
Write $\Lambda_0=A[c^{-1}]_0.$
 By \cite[Prop. 7.5]{ATV1}, the category of finite dimensional 
$\Lambda_0$-modules is equivalent to the category of normalized $A$-modules
that are $c$-torsion-free (see \S\ref{sect.normalized}). 
If $F$ is a normalized $A$-module,  then the the natural $A$-module homomorphism 
$F \to A[c^{-1}]\otimes_A F$ is an isomorphism onto the degree $\ge 0$ part, and the $\Lambda_0$-module 
corresponding to $F$ may be identified with $F_0$. Thus the dimension of a 
simple $\Lambda_0$-module equals the multiplicity of the corresponding 
normalized $A$-module. Furthermore, every simple $\Lambda_0$-module corresponds to a 1-critical $A$-module, i.e., 
to a fat point 
module.

We now present the proof in a sequence of steps.

Step 1. 
If $p-q \notin 2\ZZ\t$ and $p+q \notin E[2]+\ZZ\t$, then $\ell_{pq}$ contains 
infinitely many distinct multiplicity-$s$ fat points.
(Theorem \ref{thm6.8} shows this holds for {\it all} $p,q \in E$.)
\underbar{Proof}:  
Let $z=-p-q-2\tau.$ 
Let $\{C^\l\}$ be the multiplicity-$s$ 1-critical modules lying on
$\ell_{pq}.$
By Proposition \ref{prop.6.5}(1), the set of lines in $L(z)$ that contain only one 1-critical $C^\l$ form 
an infinite subset, $L(z)^\prime \subseteq L(z)$ say. By Proposition \ref{prop.6.6}(2),
each critical $C^\l$ lies on $\le s$ lines in $L(z)$. 
If two $C^\lambda$'s arising in this way were equivalent, then they would lie on exactly the same lines. 
Hence there are infinitely many inequivalent critical $C^\l$'s lying on $\ell_{pq}.$
$\square$

Step 2.
$\{0\}=\bigcap {\rm Ann}(F)$ where the intersection is taken over all fat point
modules $F$ such that $c.F \ne 0$ and $e(F)=s$. 
\underbar{Proof}:
The proof of Theorem \ref{thm3.6} showed that $\bigcap \Ann(M(p,q))=\{0\}$ where the intersection 
is taken over all $p$ and $q$ for which
$p-q \notin 2\ZZ\t$ and $p+q \notin E[2]+\ZZ\t$. For such $p$ and $q$, every fat point lying
on  $\ell_{pq}$ is  $c$-torsion-free so  it suffices to show that $\Ann(M(p,q)) = \bigcap \Ann(F)$ where the intersection is taken over
all multiplicity-$s$ fat points lying on $\ell_{pq}.$ By Step 1, there are
infinitely many such $F$ so the desired follows from the argument at the
end of the proof of Proposition \ref{lem3.5}.

Step 3.
In $\Lambda_0$, $\{0\}=\bigcap \Ann(S)$ where the intersection is taken over all
simple modules of dimension $s$. \underbar{Proof}: Suppose $a \in A$ is such that $ac^{-m} \in \Lambda_0$ 
annihilates all such $S$. Then by the equivalence of categories mentioned at the start of the proof,
 $a.F_0=0$ for all  multiplicity-$s$ fat points $F$ that are $c$-torsion-free. Hence by Step 2, $a=0$. 

We now prove the four statements in the theorem.

(1) 
By Step 3, $\Lambda_0$  satisfies a polynomial identity of degree $2s$. Thus, every simple $\Lambda_0$-module
has dimension $\le s.$ Suppose $\Lambda_0$ has
a simple module of dimension $<s.$ The equivalence of categories implies there
is a $c$-torsion-free critical $A$-module of multiplicity $<s.$ But
this gives a fat point of intermediate multiplicity; by Theorem \ref{thm5.4}(3), such a
fat point module is annihilated by $c$. This contradiction shows that every simple 
$\Lambda_0$-module has dimension $s$. Hence $\Lambda_0$ is an Azumaya
algebra by Artin's Theorem \cite[Thm. 8.3, Cor. 8.4]{A69}.

(2) 
This is an immediate consequence of (1) since $A$ is a domain.

(3)
By the above every simple $\Lambda_0$-module has dimension $s$, so by the 
equivalence of categories, every $c$-torsion-free fat point module has multiplicity $s$.

(4)
Let $F$ be an arbitrary fat point module of multiplicity $>1$. Let $B=A/(\O_1,\O_2)$. By \cite{SS15}, 
$B$ is isomorphic to the twisted homogeneous coordinate ring $B(E,\cL,\tau)$ so, by \cite{AV}, 
the category $\QGr(B)$ is equivalent to $\Qcoh(E)$ and under this equivalence 
the skyscraper sheaves correspond to point modules;
therefore $F$ is not a $B$-module; thus  $\O.F \ne 0$ for some $\O \in \CC\O_1 \oplus \CC\O_2$. 
Notice that $\Fract(A[\O^{-1}]_0) = \Fract (A[c^{-1}]_0)$ since both 
equal 
$$
\{ab^{-1} \, | \, a,b\in A \hbox{ are homogeneous of the same degree}\}.
$$
Hence $A[\O^{-1}]_0$ also satisfies a polynomial identity of
degree $2s.$ So every simple $A[\O^{-1}]_0$-module has dimension $\le s$. 

There is a similar equivalence of categories between $\O$-torsion-free 
$A$-modules and $A[\O^{-1}]_0$-modules. Under this equivalence, $F$ corresponds to a
simple $A[\O^{-1}]$-module of dimension equal to $e(F)$. Hence by the previous
paragraph $e(F)\le s.$
\end{pf}

{\bf Remarks. 1.}
In the notation introduced at the end of \S3, 
$$
\L_0 \; =\; A_{(c)} \; =\; \G(S_{(c)},\cA)
$$
 and $D=\Fract(\cA)$.
 
{\bf 2.}
Thus the irreducible objects in $\QGr(A)$ are the point modules
$M(p)$, the $4(s-1)$  fat point modules of intermediate multiplicity,  namely the $F(\o+k\t)$'s for $\o \in E[2]$ and $0\le k < s-1)$, 
and the remaining ones, these being the multiplicity-$s$ fat point modules. The last ones form the generic case.

\begin{theorem}
\label{thm6.8} 
Let $p,q \in E$. % and let $\ell_{pq}$ be the corresponding secant line.
\begin{enumerate}
\item{} 
$\ell_{pq}$ contains infinitely many distinct multiplicity-$s$ fat points.
\item{} Every multiplicity-$s$ fat point lying on $\ell_{pq}$ has a 1-critical representative $C$ that fits into an exact sequence
$$
0 \to M(p-s\t,q-s\t)(-s) \to M(p,q) \to C \to 0.
$$
\item{} 
$\dim_\CC\!\big({\rm Hom}_{\Gr(A)}(M(p-s\t,q-s\t),M(p,q)(s))\big)=2.$
\end{enumerate}
\end{theorem}
\begin{pf}
(1) 
A partial version of this was proved in Theorem \ref{thm6.7}.

By (2.9) applied to $M(p,q)$ we see that $M(p,q)$ has infinitely many  shift-inequivalent
1-critical quotients. A line module can
have at most 3 point modules as quotients by [8, \S5], and $A$ has
only a finite number of fat points of intermediate multiplicity.
Therefore amongst these quotients of $M(p,q)$ there must be
infinitely many having multiplicity $s$. 

(2)
This is a consequence of Theorems  \ref{thm4.3} and \ref{thm4.5}.

(3)
(This was proved in Proposition \ref{prop.6.2} when $p-q \notin 2\ZZ\t.$) It follows
at once from (1) and (2) that the dimension of this ${\rm Hom}$-space is at
least 2. If it were strictly larger, then $\dim\,{\rm Hom}(M(p-s\t,q-s\t),
M(p+\t,q-\t)(s-1)) \ge 2$ because $M(p,q)$ contains a copy of
$M(p+\t,q-\t)(-1)$ that is of codimension 1 in degree $s$. 
The cokernels of these maps are quotients of $M(p+\t,q-\t)$ of 
multiplicity $s-1.$ As in the proof of Proposition \ref{prop.6.3}, it follows that infinitely 
many of these quotients are critical, i.e., fat points of
multiplicity $s-1$. By Theorem \ref{thm4.5},  it follows that $M(p-s\t,q-s\t) \cong
M(p+\t-(s-1)\t,q-\t-(s-1)\t)$. Thus
$\ell_{p-s\t,q-s\t}=\ell_{p-(s-2)\t,q-s\t}.$ But this is impossible.
\end{pf}

\section{  %Section 7. 
Classification of multiplicity-$s$ fat points}
\label{sect.fat.pts.s}

As always, $s$ is the order of $2\tau$, i.e., the smallest positive integer such that $s \tau \in E[2]$.
 
\subsection{}
In this section we  show there is a 3-dimensional variety parametrizing the multiplicity-$s$ fat points in $\Projnc(A)$.  
Together with the results in \S5 this will classify the ``generic'' fat points. 
The fat points remaining unclassified are those of multiplicity $s$ that lie on
the lines $\{\ell_{pq} \; | \; p+q \in E[2]+\ZZ \tau\}.$

By Theorem \ref{thm4.3}, fat points lie on secant lines (a given fat point can lie on infinitely
many different secant lines) and, conversely, each secant line  contains infinitely many
different fat points. In order to clarify the situation we need to
introduce a notion of equivalence between lines. 

Several important facts about line modules are summarized in \S\ref{sect.prelims}. 

%However, we first recall some results from \cite{LS8}.

%It is shown in \cite[\S3]{LS8} that the secant lines of $E$ occur in families 
%as follows. For each $z \in E$ we refer to $L(z):=\{ \ell_{pq} \; | \;p+q=z \}$, as a {\it family of lines}.  The
% union of the lines in $L(z)$ forms a quadric $Q(z)$ say, that contains $E$. There is a pencil of quadrics containing $E$, and each is of the form $Q(z)$ for some $z \in E.$ Furthermore, $Q(z)=Q(-z)$ and these are the 
%only equalities among the quadrics. Hence there are two families of lines
%on $Q(z)$ namely $L(z)$ and $L(-z)$: these are the two rulings on $Q(z) \cong \PP^1 \times \PP^1$ when $Q(z)$ is smooth. If $\o \in E[2]$ then $Q(\o)$ is singular, and contains only one family of lines; these are 
%the only singular quadrics in the pencil.

%There is a unique degree 2 central element $\O(z)$ that annihilates the line modules $\{M(p,q) \; | \; p+q=z\}$.  In analogy with the 4 singular quadrics there are 4 central elements  that annihilate only one family of line modules, namely
%$\O(\o-\t)$ for $\o \in E[2].$ 

\subsubsection{Definition}
Two secant lines $\ell$ and $\ell^\prime$ are {\sf equivalent}, denoted 
$\ell \sim \ell^\prime,$ if there are infinitely many different fat  points that lie on both of them.  
The equivalence class of a line $\ell$ is denoted by $[\ell]$.

\subsubsection{}
When $Q(z)$ is smooth, i.e., when $z\notin E[2]$, each point on it is the intersection of a unique pair of 
lines one from $L(z)$ and one from $L(-z)$.  Our classification of multiplicity-$s$ fat points follows this classical model:  
we will characterize a fat point in terms of the lines in $\Projnc(A)$ that it lies on.  Although this strategy must be refined it is helpful:
Theorem \ref{thm7.2} shows that if $z \notin E[2] +\ZZ\t$, then a multiplicity-$s$  fat point annihilated by 
$\O(z)$ lies on exactly two equivalence
classes of lines, one from the family $L(z)$ and one from the family
$L(-z-2\t)$. Conversely, if $z \notin E[2] +\ZZ\tau$, then most pairs
of lines $\ell\in L(z)$ and $\ell^\prime \in L(-z-2\tau)$ have a unique
common multiplicity-$s$ fat point.

\begin{proposition}
\label{prop7.1} 
Let $p,q \in E$ be such that $p-q \notin 2\ZZ\t$. Then 
\begin{enumerate}
\item{} 
$\ell_{pq}$ and $\ell_{p+2i\t,q-2i\t}$ are
equivalent for all $i \in \ZZ$;
\item{} if $p+q\notin E[2] +\ZZ\t$, then $[\ell_{pq}] = \{\ell_{p+2i\t,q-2i\t}
\; | \; i \in \ZZ\}$;
\item{} 
if $p+q \in E[2]+k\t$, then $[\ell_{pq}]= 
\{\ell_{p+2i\t,q-2i\t}, \, \ell_{p+2i\t,q-2(k+i+1)\t} \; | \; i \in \ZZ \}$.
\end{enumerate}
\end{proposition}
{\bf Remark.} 
If $p+q\in E[2] -\t$, then Proposition \ref{prop7.1}(3) will show that
$$
[\ell_{pq}] \;=\;  \{\ell_{p+2i\t,q-2i\t} \; | \; i \in \ZZ\}.
$$

\begin{pf}
(1) 
(The proof of (1) works without any assumption on $p-q$.)
Since $A$ satisfies a polynomial identity,  (2.8) implies that 
$A/ \! \Ann(M(p,q))$ has infinitely many inequivalent 1-critical
quotients. In other words, infinitely many different fat points lie on $\ell_{pq}.$ 
By Proposition \ref{prop4.4}(4), these fat points also lie on $\ell_{p+2\t,q-2\t}$ so $\ell_{pq}$ and $\ell_{p+2\t,q-2\t}$ 
are equivalent. The result now follows by induction.

(2) 
Write $z=p+q$ and let $\{C^\l\}_{\l \in \PP^1}$ be the multiplicity-$s$ quotients of
$M(p,q)$ constructed in Proposition \ref{prop.6.3}. 
Suppose  $\ell_{p' q'} \sim \ell_{pq}$. Then $M(p,q)$ and
$M(p',q')$ have the same annihilator since they have infinitely many
common fat points. Therefore $\O(p+q) =\O(p' + q')$, whence either
$p' + q'=z$ or  $p' + q'=-z-2\t.$  In particular, every line
equivalent to $\ell_{pq}$ lies in $L(z) \cup L(-z-2\t).$

Suppose $p^\prime + q'=z.$
It was already proved in Proposition \ref{prop.6.4} that if $\ell_{p' q'}$ is not of the
form $\ell_{p+2i\t,q-2i\t}$ then there is at most one multiplicity-$s$ fat point that lies on both $\ell_{pq}$ and $\ell_{p' q'}.$
Therefore (without any hypothesis on $p+q$)
the only lines in $L(z)$ equivalent to $\ell_{pq}$ are the lines
$\ell_{p+2i\t,q-2i\t}.$

Suppose  $p' + q' = -z-2\t.$ Then $\ell_{p' q'}$ is equal to $\ell_{p+2i\t,q+2j\t}$ for some $i,j \in \ZZ$ because if it were
not, then Proposition \ref{prop.6.5}(2) would show there is a unique $C^\l$ lying on $\ell_{p' q'}$ thereby contradicting the hypothesis that 
$\ell_{pq} \sim \ell_{p' q'}.$  It follows that $z= \o +k\tau$ for some $\o \in E[2]$
and $k=i+j-1.$ 
Hence if $p+q \notin E[2]+k\tau$ then the previous paragraph applies, and
(2) follows.

(3)
Suppose $z=\o+k\tau.$
By Proposition \ref{prop5.7},  $\ell_{pq} 
\sim \ell_{p,q-2(k+1)\t}$ and hence, by (1), $\ell_{pq} \sim \ell_{p+2i\tau,
q-2(k+i+1)\tau}.$ The previous paragraph shows these are the only
lines in $L(-z-2\t)$ that are equivalent to $\ell_{pq}.$
\end{pf}

\begin{theorem}
\label{thm7.2}
%{\bf Theorem 7.2.}
Let $F$ be a multiplicity-$s$ fat point module. If $\Omega(z).F=0$, then  $F$ lies on a line in $L(z)$ and on 
a line in $L(-z-2\t).$
\end{theorem}
\begin{pf}
By Theorem \ref{thm4.3}, 
 $F$ lies on some secant line; there is no loss of generality in assuming that line belongs to $L(z)$
so we will do that. The result follows from Proposition \ref{prop.6.6} if this line is of the 
form $\ell_{pq}$ with $p-q \notin 2\ZZ\t$.  If $F$ does not lie on a 
line of this form, choose $\ell_{pq} \in L(z)$ not containing $F$.

Write $\O=\O(z).$
Let $\ell^\perp \subseteq A_1$ be the subspace vanishing on $\ell_{pq}$. As
in Theorem \ref{thm4.3}, for some $u \in \ell^\perp-\{0\}$ there is a non-zero map $\varphi: A/Au+A\O \to F$ and  a diagram
$$
\xymatrix{
0 \ar[r] & M(\ell')(-1) \ar[r]^\phi &   A/Au+A \Omega  \ar[d]^\varphi  \ar[r]^\psi &  M(\ell_{pq}) \ar[r] & 0
\\
&&  F
}
$$
in which the row is exact. As in  \cite[Thm. 4.2]{SS16},  $\ell'=\ell_{p' q'}$ 
where 
$$
(u)_0\;=\; (p)+(q)+(p'+\t) +(q'+\t).
$$
Since $F$ does not lie on $\ell_{pq}$ the composition $\varphi\phi:M(\ell')(-1) \to F$ is non-zero. 
Hence, by Proposition \ref{prop4.4}(3), $F$ lies on $\ell_{p'-\t, q'+\t}.$ 
The result follows since $p'+q'=-p-q-2\t=-z-2\t.$
\end{pf}

\subsubsection{Remark}
 If there were a simple proof that no multiplicity-$s$ fat points lie on every line in $L(z)$, then Theorem \ref{thm7.2}  
would follow easily without requiring the  careful  analysis in Proposition \ref{prop.6.6}.

\begin{proposition}
\label{prop7.3} 
If $z\notin E[2]+\ZZ\t$, then no line in $L(z)$ is equivalent to a line in $L(-z-2\t).$\footnote{By Proposition \ref{prop7.1}(3), the hypothesis that  $z\not\in E[2]+\ZZ\t$ is necessary.}
\end{proposition}
\begin{pf}
Suppose to the contrary that $\ell \in L(z)$ and $\ell' \in L(-z-2\t)$ are equivalent.

Then $\Ann(M(\ell))=\Ann (M(\ell'))$. By Proposition \ref{prop1.1}, 
there is an injective map $\varphi:M(\ell')(-k) \to M(\ell)$ for some $k\ge 0$. 
Without loss of generality, we will assume that $\ell$ and $\ell'$ are such that $k$ is as small as possible. 
Since $z \ne -z-2\t$,  $\ell \ne \ell'$. Hence $k \ge 1.$

Since $\coker(\varphi) \ne 0$, its GK-dimension is 1 so it has a cyclic 1-critical quotient, $C$ say. 
Hence there is an exact
sequence 
$0 \to M(\ell'')(-e) \to M(\ell) \to C \to 0$ where
$e=e(C)$ and ${\rm im}(\varphi) \subseteq M(\ell'').$  Thus $\varphi$
gives an injective map $M(\ell')(-(k-e)) \to M(\ell'').$ In
particular, $\ell'$ and $\ell''$ are equivalent.

There are 3 possibilities for $C$, and accordingly 3 possibilities for $\ell'':$
\begin{enumerate}
\item{} $C \cong M(p)$ for some $p \in E;$
\item{} $C \cong F(\o+(e-1)\t)$ for some $\o \in E[2];$
\item{} $C$ is a multiplicity-$s$ fat point.
\end{enumerate}
By Theorem \ref{thm4.5},  (2) does not occur  since $z \notin E[2]+\ZZ\t.$ 
If (1) or (3) occurs, then   $\ell'' \in L(z)$; combined with the conclusion of the previous paragraph, this contradicts 
the minimality of $k$. Thus, $\ell$ and $\ell'$ are not equivalent. The truth of the proposition  follows.
\end{pf}

\begin{theorem}
\label{thm7.4} 
Assume $z \notin E[2] +\ZZ\t$ and let  $\ell \in L(z)$.
\begin{enumerate}
\item{}
 If $\ell' \in L(-z-2\t)$,  there is at most one fat point that lies on both $\ell$ and $\ell'.$ 
\item{} 
Only finitely many $\ell' \in L(-z-2\t)$ fail to have a fat point lying on $\ell$.
\end{enumerate}
\end{theorem}
\begin{pf}
 Write $\ell=\ell_{pq}.$
By Theorem \ref{thm6.8}(2), the fat points lying on $\ell$ are obtained as follows. Let 
$$
\psi_1,\psi_2: M(p-s\t,q-s\t)(-s) \; \longrightarrow M(p,q)
$$
 be linearly 
independent. For each $\l=(\l_1,\l_2) \in \PP^1$ let $\psi_\l=\l_1\psi_1 
+ \l_2\psi_2$, and define $C^\l=M(p,q)/{\rm im}(\psi_\l)$. The argument in 
Proposition \ref{prop.6.3} shows that $C^\l$ is critical for all except a finite number of 
$\l.$ Hence the multiplicity-$s$ fat points lying on $\ell$ are 
parametrized by $\PP^1$ minus a finite number of points; write $\PP_\ell$ 
for this projective line.  We can identify $L(-z-2\t)$
with $\PP^1$ also. Define 
$$
Z \; :=\; \{ (\l, \ell') \; | \;C^\l \hbox{ lies on } \ell' \} \; \subseteq \; \PP_\ell \times L(-z-2\t) \; \cong \; \PP^1 \times \PP^1.
$$

To see that $Z$ is a closed subvariety we appeal to a rather general
principle. Let $\varphi:V_1 \times V_2 \to V_3$ be a linear map of
finite dimensional vector spaces, fix $n_1,n_2,d \in \NN$  and define
$Z=\{(U_1,U_2)\in \GG(n_1,V_1) \times \GG(n_2,V_2) \,\; | \;\,
\dim\,\varphi (U_1 \times U_2) \le d\}$ where $\GG(n_i,V_i)$ denotes the
Grassmannian of $n_i$-dimensional subspaces of $V_i$. Then $Z$ is a
closed subvariety.

The projection $\pi_2:Z \to L(-z-2\t)$ to the second factor  is a proper
morphism. If it were constant, then $\ell$ would be equivalent to
a line in $L(-z-2\t)$ which would contradict Proposition \ref{prop7.3}. Since $\pi_2$ is
proper it   is surjective. Hence there is a $C^\l$ that lies on both $\ell'$ and $\ell$.

To show there is at most one $C^\l$, and hence at most one fat point, lying on both $\ell$ and $\ell'$, 
we must show that $\pi_2$ has degree 1.
 It suffices to show that the fiber of $\pi_2$ over a general point of $L(-z-2\t)$ consists of only one point. 
The lines $\{ \ell_{p'q'}\in L(-z-2\t) \; \vert \; p'-q'
\notin 2\ZZ\t \}$ form a Zariski dense open subset of $L(-z-2\t)$ and
for these lines  Proposition \ref{prop.6.5}(1) shows there is only one $\l
\in \PP^1$ such that $C^\l$ lies on $\ell_{p'q'}.$
\end{pf}

\begin{corollary}
\label{cor7.5} 
If $F$ and $F'$ are multiplicity-$s$  fat points lying on exactly the same lines, then $F \sim F'$.
\end{corollary}
\begin{pf}
The module $F$ is annihilated by $\O(z)$ for some $z \in E$.
 By Theorem \ref{thm7.2},  $F$ lies on a line $\ell \in L(z)$ and on a line $\ell' \in L(-z-2\t)$. By hypothesis, so does $F'$. Hence,
by Theorem \ref{thm7.4}, $F$ and $F'$ are equivalent.
\end{pf}

\begin{theorem}
\label{thm7.6} 
Let $p,q \in E$.
\begin{enumerate}
\item{} 
If $p+q\notin E[2] +\ZZ\t$, then  $[\ell_{pq}] \,=\, \{\ell_{p+2i\t,q-2i\t}  \; | \;  i \in \ZZ\}$.
\item{} 
If $p+q=\o+k\t$ for some $\omega \in E[2]$,  then 
$$
[\ell_{pq}] \;=\;  \{\ell_{p+2i\t,q-2i\t}, \, \ell_{p+2i\t,q-2(k+i+1)\t}  \; | \;  i \in \ZZ \}.
$$
\end{enumerate}
\end{theorem}
\begin{pf}
By Proposition \ref{prop7.1}(1), $\ell_{pq} \sim \ell_{p+2i\t,q-2i\t}$ for all $i \in \ZZ$ (if $p-q \notin 2\ZZ\tau$).

Proposition \ref{prop5.7} and Theorem \ref{thm6.8}(1) show that $\ell_{pq} \sim \ell_{p,q-2(k+1)\t}$ whenever
$p+q=\o+k\t.$ Hence to complete the proof it suffices to prove the
following (stronger) statement.

\underbar{Claim:}
Assume $p+q \notin E[2]-\t.$ 
If $\ell_{p_1 q_1 },\ell_{p_2 q_2} \in L(p+q)$ have a common multiplicity-$s$ fat point, 
then $\ell_{p_2 q_2}=\ell_{p_1+2i\t,q_1-2i\t}$ for some  $i \in \ZZ.$

\underbar{Proof}:
Suppose $\ell_{p_1 q_1}$ and $\ell_{p_2 q_2}$ have
a common multiplicity-$s$ fat point. By Theorem \ref{thm7.2},  that fat point also lies on a
line $\ell \in L(-z-2\t).$ We adopt the notation in  Theorem \ref{thm7.4} with
the roles of $z$ and $-z-2\t$ interchanged. As in Theorem \ref{thm7.4}, let $Z
\subseteq\PP_\ell \times L(z)$ be the closed subvariety of pairs
$(\l,\ell')$ such that $C^\l$ lies on $\ell'$ ($C^\l$ is a
multiplicity-$s$ quotient of $M(\ell)$ that is critical for all
except  finitely many $\l$). 

Let $E'$ denote the isogenous elliptic curve $E/\langle 2\t \rangle$, and define an equivalence relation $\sim$ on $E'$ by $p
\sim z-p$. The equivalence classes are orbits under an action of
the group $\ZZ_2$, so we can form the quotient variety
 $E'/\!\!\sim \, \cong\PP^1.$ If $p \in E$, write $[p]$ for
the equivalence class of $p$ in $E'/\!\!\sim.$ 
We also define an equivalence relation on $E$ by $p \sim z-p$. We identify
$E/\!\sim$ with $L(z)$ by making $p$ correspond to the line
$\ell_{p,z-p}.$ There is a commutative
diagram
$$
\xymatrix{
& E \ar[d]   \ar[r]^{\alpha}   & E/\!\sim  \ar[d]^{\c} & =L(z)  \cong \PP^1 
\\
& E'  \ar[r]_{\a'} & E'/\!\sim & \cong \PP^1  
\\
}
$$
where the left-most vertical arrow is the quotient map and $\g$ is the unique surjection making the diagram commute. 
We  therefore obtain a morphism  
$$
\hbox{Id} \times \c : \PP_\ell \times L(z) \longrightarrow \PP_\ell \times  (E'/\!\sim).
$$
Define $W=(\hbox{\rm Id} \times \c)(Z)$, and let $\theta_1:W \to
\PP_\ell$ and $\theta_2:W \to E'/\!\!\sim$ be the projections. Thus
$W$ is the closed subvariety consisting of pairs $(\mu,[r])$
such that $C^\mu$ (which lies on $\ell$) lies on $\ell_{r, z-r}$; 
we are implicitly using Proposition \ref{prop7.1}(1)  where we showed that $\ell_{p' q'} \sim
\ell_{p' +2i\t, q'-2i\t}$ for all $i \in \ZZ.$ 

Thus, the hypothesis is that there is a $\l \in \PP_\ell$ such that $(\l,[p_1])$ and $(\l,[p_2])$ belong to $W$, and
we wish to prove that $[p_1]=[p_2]$. It suffices to show that ${\rm deg}(\theta_1)=1.$
Clearly 
$$
U\; := \; \{(\mu,[r]) \in W \; | \;r-(z-r) \notin 2\ZZ\t \} \; =\; \theta_2^{-1}(E'/\!\sim  \,- \,[z+E[2]'])
$$ 
is a dense open subset  of $W$. 
Now fix some $(\mu,[r]) \in U$; then $C^\mu$ lies on
$\ell_{r,z-r}$ and $r -(z-r) \notin 2\ZZ\t$ so Proposition \ref{prop.6.4}(2) applies to
$\ell_{r,z-r}$.  However, that result shows that $\theta_2\theta_1^{-1}(\mu)=[r]$ is 
a singleton. But ${\rm deg}(\theta_2)=1$ by Theorem \ref{thm7.4}, so $\theta_1^{-1}(\mu)$ is a 
singleton.  Therefore ${\rm deg}(\theta_1)=1.$
\end{pf}

\begin{corollary}
%{\bf Corollary 7.7.} 
\label{cor.7.7}
Suppose the order of $\tau$, $n$, is odd. If $F$ is a fat point such
that $cF\ne 0$, then $F \sim F(1).$ Consequently, there is a bijection
between the sets:
\begin{enumerate}
\item{} equivalence classes of simple $A$-modules $S$, such that $c.S
\ne 0;$
\item{} fat points $F$ such that $c.F \ne 0;$
\item{} graded prime ideals $\fp$ such that $c \notin \fp$ and $d(A/\fp)=1.$
\end{enumerate}
\end{corollary}
\begin{pf}
If $p,q \in E$, then $F$ lies on $\ell_{pq} \Longleftrightarrow F(1)$
lies on $\ell_{p+\t,q-\t}$ by Proposition \ref{prop4.4}(2) $\Longleftrightarrow F(1)$
lies on $\ell_{pq}$ by Proposition \ref{prop7.1}. Therefore $F$ and $F(1)$ lie on exactly 
the same lines, so $F \sim F(1)$ by Theorem \ref{thm7.4}. The rest of the result follows
from (2.6).
\end{pf}

{\bf Remarks. 1.} 
If $n$ is even it is no longer true that $F \sim F(1).$
Suppose $n$ is even. Let $F$ be a fat point such that $cF \ne 0.$ If
$F \sim F(1)$, then $\ell_{pq}$ and $\ell_{p+\t,q-\t}$ have a fat point in
common which contradicts Theorem \ref{thm7.6}.

{\bf 2.}
If $\ell_{pq}$ and $\ell_{p^\prime q^\prime}$ are equivalent lines, then 
$$
\Ann(M(p,q)) \;=\;  \Ann(M(p^\prime,q^\prime))
$$
by Corollary \ref{cor.1.5}. The converse is true if $n$ is odd, and false if $n$ is even. 

First, suppose $n$ is even and choose $p,q$ such that $p+q \notin E[2] +\ZZ\tau.$
The modules $M(p+\tau,q-\tau)(-1)$ and $M(p,q)$ have the same annihilator because the first embeds in the second
and $\O M(p,q)$ is isomorphic to a non-zero submodule of the first for some $\O\in \CC\O_1+\CC\O_2$. 
However, by Theorem \ref{thm7.6}, $\ell_{pq}$ and $\ell_{p+\t,q-\t}$ are not equivalent.

Now suppose $n$ is odd, and that $\Ann(M(p,q))=\Ann (M(p^\prime,q^\prime))$. 
By Proposition \ref{prop1.1}, there is an
injection $M(p^\prime,q^\prime)(m)\to M(p,q)$ for some $m$ so, by Lemma \ref{lem3.2},
$\{p^\prime,q^\prime\}=\{p+(d-2k)\tau,q-d\tau\}$ for some $d,k \in \ZZ$
and either $2k\tau=0$ or (if $2k\tau \ne 0$) $p+q=\o+(k-1)\tau$ for 
some $\o \in E[2]$. If $2k\tau=0$ then it follows from Theorem \ref{thm7.6} that the lines
are equivalent (because $n$ is odd it doesn't matter if $d$ is odd).
If $2k\tau \ne 0$ then by  Theorem \ref{thm7.6}(2) the lines are equivalent.

\begin{theorem}
%{\bf Theorem 7.8.} 
\label{thm.7.8}
If  $\tau$ has order $n<\infty$, then
\begin{enumerate}
\item{} 
$A[c^{-1}]$ is an Azumaya algebra of rank $n^2$ over its center;
\item{} 
$A$ satisfies a polynomial identity of degree $2n$ and none of lower degree.
\end{enumerate}
\end{theorem}
\begin{pf}
Clearly, it suffices to prove (1). To do this we will show
that every simple $A[c^{-1}]$-module has dimension $n$ and then invoke
Artin's Theorem \cite[Thm. 8.3, Cor. 8.4]{A69}. Thus it suffices to show that, if $F$ is a fat
point module not annihilated by $c$, then $F$ has an $n$-dimensional
simple quotient. Let $F$ be a fat point module such that $c.F \ne
0.$ Thus $e(F)=s.$

First suppose $n$ is odd. By Corollary \ref{cor.7.7}, $F \sim F(1)$ so, by Proposition \ref{prop.1.2},
$F$ is the only 1-critical $A/\!\Ann(F)$-module up to equivalence.
Hence (2.10) applies. This says that every simple $A/\!\Ann(F)$-module
has dimension $e(F)=s$. But $n$ is odd, so $s=n$.

Now suppose  $n$ is even. Let $S$ be a non-trivial 
simple quotient of $F$. Then $S$ is also a quotient of $F(-1)$.
By Proposition \ref{prop.old.2.2}, both $F$ and $F(-1)$ embed in $\widetilde S$ via degree 0 maps.
By Remark 1 above, $F$ and $F(-1)$ are not equivalent.
Therefore, their images in $\widetilde S$ intersect in 0. It follows
that $\dim(\widetilde S_m) \ge 2e(F)=n$ for $m \ge 1.$ Thus $\dim(S) \ge
n.$ However, there is a degree-two central element $\Omega$ such that
$\Omega.F \ne 0$. The dimension of $F/(1-\Omega)F$ is $n$ so it has a non-trivial simple quotient of dimension $\le n.$ 
By Theorem \ref{2.5}, all non-trivial simple quotients of $F$ have the same dimension so $\dim(S)=n.$
\end{pf}

Our next goal is to parametrize the fat points that are not annihilated
by the central element $c$ in Theorem \ref{thm6.7}. All these fat points have multiplicity $s$.

\begin{corollary}
%{\bf Corollary 7.9.} 
\label{cor7.9}
Assume $z \notin E[2]+\ZZ\t.$ The
fat points annihilated by $\O(z)$ are parametrized by a dense 
subset of $\PP^1 \times \PP^1$ and they all have multiplicity $s$.
\end{corollary}
\begin{pf}
Let $F$ be such a fat point. Its multiplicity is $s$
by \S5. By Theorem \ref{thm7.2}, $F$ lies on a line $\ell_{pq} \in L(z)$ and on a line
$\ell_{p' q'} \in L(-z-2\t).$ 

Set $E'=E/\langle 2\t \rangle$. Define equivalence relations
$\sim$ and $\approx$ on $E'$ by $p \sim z-p$ and $p \approx
-z-2\t-p$. The quotient varieties $E'/\!\sim$ and $E'/\!\approx$  are isomorphic to $\PP^1.$ Although
the lines $\ell_{pq}$ and $\ell_{p' q'}$ that contain $F$ are
not uniquely determined, their equivalence classes are by Theorem \ref{thm.7.8}.
Hence $F$ determines a point $([p],[p']) \in (E'/\!\sim)
\times (E'/\!\approx) \cong \PP^1 \times \PP^1$. Conversely, by
Theorem \ref{thm7.4}, for a dense set of $([p],[p'])$ the lines $\ell_{pq}$ and
$\ell_{p' q'}$ have a common fat point.
\end{pf}

{\bf Remarks.} 
1. 
The parametrization of the fat points in Corollary \ref{cor7.9} may be
described more elegantly by using the geometry of the isogenous elliptic curve $E'=E/\langle 2\tau \rangle$ as a quartic space curve.  
To do this we first embed $E/\langle 2\tau \rangle$ in $\PP^3$ via the complete linear system $|4(0)|$ and identify $E'$ with its
image. Thus, four points of $E'$ sum to 0 if and only if there is a hyperplane in $\PP^3$ whose scheme-theoretic intersection with
$E'$ is the sum (in the sense of divisors) of those four points. 

We label the quadrics containing $E'$ in the following way: $Q(\overline{z})$, $\overline{z} \in E'$, is the union of the secant lines 
$\ell$  to $E'$ such that the sum of the points in the scheme-theoretic intersection $\ell \cap E'$ is $\overline{z}$. Thus, if
$\overline{z}$ is the image in $E'$ of a point $z \in E$, then  $Q(\overline{z})$ is smooth if and only if $z \notin E[2]+\ZZ\t$. 
For such $z$, the two rulings on $Q(\overline{z})$ are $L(\overline{z})$ and $L(-\overline{z})$. 
More importantly for us, there is a bijection between the lines in $L(\overline{z})$ and the equivalence classes
of lines in $L(z).$ In fact, the map $[\ell_{pq}] \to \ell_{\overline{p} \overline{q}}$ is an isomorphism of varieties $L(z)/\!\!\sim \, \to
L(\overline{z})$. There is a similar isomorphism $L(-z-2\t)/\!\!\sim  \, \to
L(\overline{-z-2\t})=L(\overline{-z})$ to the lines in the other ruling on $Q(\overline{z}).$

Now define a map 
$$ 
\{\hbox{multiplicity-$s$ fat points annihilated by $\O(z)$}\} \; \longrightarrow \; \{ \hbox{points on  $Q(\overline{z})$}   \}
$$
by sending a fat point $F$ to 
$\ell_{\overline{p_1} \overline{q_1}} \cap
\ell_{\overline{p_2} \overline{q_2}}$
whenever $F$ lies on $\ell_{p_1 q_1} \in L(z)$ and $\ell_{p_2,q_2} \in
L(-z-2\t)$. By Corollary \ref{cor7.5}, this is injective, and as in the proof of Corollary \ref{cor7.9} the
image is dense in $Q(\overline{z})$;  the image of this map is probably $Q(\overline{z})-E'.$

2.
When $z \in E[2]+\ZZ\t$ there should be a similar point of view, but the
difficulty is that the quadric $Q(\overline{z})$ is then singular so  the points on it are not given by intersections of lines.

\subsection{Some 3-folds}
\label{sect.3-folds}
We now define a 3-fold that parametrizes most multiplicity-$s$ fat points. Secant lines are parametrized by degree-two effective divisors on $E$, i.e., by points on the second symmetric power $S^2E$.
Now consider the surface 
\begin{equation}
\label{3-fold.S}
S \;  := \; \frac{E \times E}{\langle (2\t,-2\t)\rangle} \Bigg/   \ZZ_2
\end{equation}
obtained by first forming the quotient of $E\times E$ by the cyclic subgroup
of order $s$ generated by $(2\t,-2\t)$, and then quotienting out by
the $\ZZ_2$-action that swaps the coordinates. 

The surface $S$ parametrizes the equivalence classes of lines under the equivalence relation
generated by $\ell_{pq} \approx \ell_{p+2\t,q-2\t}$ (this is not
the same as the equivalence $\sim$ that we are really interested in but , by Theorem \ref{thm.7.8}, $\sim$ and $\approx$
are closely related; indeed, for most lines the two equivalence relations coincide). 

Now define for each $x \in E$,
$$
(S \times S)_{x}  \; := \;  \{ ((p,q),(p',q')) \in S \times S  \; | \; p+q+p'+q'   \in  x+\langle 2\t \rangle  \}.
$$
The 3-folds $(S \times S)_{x}$ are the fibers of the right-most vertical morphism in the commutative diagram
$$
\xymatrix{
E^2 \times E^2 \ar[d]_+ \ar[r] & S \times S \ar[d] 
\\
E \ar[r] & E/\langle 2\t \rangle
}
$$
Let
\begin{equation}
\label{3-fold.T}
T\; := \; (S\times S)_{-2\t}/\ZZ_2
\end{equation}
be the quotient  by the $\ZZ_2$ action that sends $((p,q),(p',q'))$ to $((p',q'),(p,q))$.

\begin{theorem}
\label{thm7.10}
 Define 
$$
c \; :=\;        \prod_{\omega \in E[2]}      \prod_{k=0}^{s-2}  \Omega(\omega+k\tau).
$$
The $c$-torsion-free fat points of multiplicity $s$ are parametrized by a Zariski dense open subset of
the 3-fold $T$ in (\ref{3-fold.T}).
\end{theorem}
\begin{pf}
We keep the notation defined just prior to the statement of the Theorem. 

Let
$$
U'  \; :=\;  \{((p,q),(p',q')) \in T \; | \;p+q \notin E[2]+\ZZ\t  \subseteq T
$$ 
Clearly $U'$ is a Zariski dense open subset of $T$.

Let $F$ be a multiplicity-$s$ fat point  such that
$cF\ne 0$. Then $\O(z).F=0$ for some $z \notin E[2]+\ZZ\t.$
By Theorem \ref{thm7.2},   $F$ lies on lines $\ell_{pq} \in L(z)$ and $\ell_{p' q'} 
\in L(-z-2\t)$ such that $p+q+p'+q'=-2\t$. By Theorem \ref{thm.7.8}, the pairs $(p,q)$ for which $F$ lies on $\ell_{pq} $ all map to 
the same point of the surface $S$.
The same remark applies to $(p',q')$. Since $\O(z)=\O(-z-2\t)$, this procedure associates to $F$ a unique 
point $((p,q),(p',q')) \in T;$ the fact that $F$ is
$c$-torsion-free ensures that this point is actually in $U'$. 
Conversely, Theorem \ref{thm7.4} shows that
given $((p,q),(p',q')) \in U'$ then the lines $\ell_{pq}$
and $\ell_{p' q'}$ have at most one common fat point; in fact
Theorem \ref{thm7.4} shows that if we
remove a suitable hypersurface from $U'$, we obtain a Zariski open
subset $U \subseteq U'$ such that for each $((p,q),(p',q')) \in
U$ the lines $\ell_{pq}$ and $\ell_{p' q'}$ contain a
unique common multiplicity-$s$ fat point. Hence the points of $U$
are in bijection with the multiplicity-$s$ fat points that are $c$-torsion-free.
\end{pf}

\begin{theorem}
\label{thm7.11}
The variety $T$ is rational. Thus, the $c$-torsion-free
fat points are parametrized by a rational variety.
\end{theorem}
\begin{pf}
 First we observe that $T$ is isomorphic to $Y:= (S\times S)_0/\ZZ_2$ 
where the $\ZZ_2$ action swaps the two coordinates. The isomorphism $\psi:T \to 
Y$ is given as follows: fix $\eta \in E$ such that $2\eta=\t,$ and
if $t=((p,q),(p',q')) \in T$ then $\psi(t)=((p+\eta,q+\eta),(p'+\eta,q' +
\eta)) \in Y.$

Define $E'=E/\langle 2\t \rangle$ and embed $E'$ in $\PP^3$ as in the remarks after Corollary \ref{cor7.9}.
The labelling of the quadrics containing $E'$ is such that $Q(\xi)=Q(-\xi)$ and $Q(\xi)$ is the union of the lines
$\ell_{p,\xi-p}$ for $p \in E'$. Fix a hyperplane $H \subseteq \PP^3$ that does
not pass through $0 \in E'$. Let $\pi :\PP^3 \to H$ be the projection with
center $0.$ The restriction of $\pi$ to each $Q(\xi)$ is a birational isomorphism $\pi:Q(\xi) \to H$.

We will now define a rational map $f:Y \to H \times \PP^1 \cong \PP^2 \times \PP^1$. 

Let $((p,q),(p',q')) \in Y.$ Write $\overline \ell_{pq}$
for the secant line to $E'$ that passes through the images of $p$ and $q$
in $E'.$ Then $\overline \ell_{pq}$ depends only on the image of $(p,q)$ in
$S$. Define $\overline \ell_{p' q'}$ similarly. Since $p+q=-p'- q'$ in $E/\langle 2\tau \rangle$,
 $\overline \ell_{pq}$ and $\overline \ell_{p' q'}$ lie on a common $Q(\xi)$. 
 Suppose $p+q \notin E[2]+\langle 2\t \rangle.$ Then $Q(\xi)$ is smooth, and $\overline
\ell_{pq}$ and $\overline \ell_{p' q'}$ belong to different rulings so
intersect at some point $y \in Q(\xi).$ Define $f((p,q),(p',q')) =
(\pi(y), \pm(p+q)).$ This is well-defined since $y$ depends only on
$((p,q),(p',q'))$ as a point of $Y$, and we are viewing $\pm(p+q)$ as an
element of $E'/\pm \cong \PP^1.$

The injectivity of $f$ on its domain of definition follows from the fact
that the restriction of $\pi$ to each $Q(\xi)$ is injective. Furthermore, it
is obvious that the image of $f$ is dense, so $Y$ is birational to 
$\PP^2 \times \PP^1$ and hence to $\PP^3$.
\end{pf}

\begin{theorem}
\label{thm7.12} 
The division algebra $D=\Fract(\cA)$ is of degree $s$ over its center, and that center is a field of rational functions in 3 variables.
\end{theorem}
\begin{pf}
 Since the $c$-torsion-free fat points are in bijection with the simple modules over the Azumaya algebra $A[c^{-1}]_0$, 
 it follows that the center of $A[c^{-1}]_0$ is rational, and hence the center of $\cA$ is
rational. The center of $D=\Fract(\cA)$ is therefore rational.
\end{pf}

 {\bf Remark.} 
Four fibers of the map $Y \to \PP^1$, $((p,q),(p',q')) \to \pm(p+q) \in E'/\pm,$ are isomorphic to $\PP^2$ and all others 
are isomorphic to $\PP^1 \times \PP^1$.   

One can also describe $Y$ as follows.
For $p \in E$ write $\overline{p}$ for its image in $E' =E/\langle 2\t \rangle$. 
Set $Z=E' \times E'\times E$ and define an action of the dihedral group 
$D_4= \langle \a, \b \; | \;\a^2= \b^2 =(\a\b)^4 =1 \rangle$ on $Y$ by 
$ (\overline{p},\overline{q},\xi)^\a = (\overline{\xi -p}, \overline{q}, \xi)$ and $( \overline{p},\overline{q},\xi)^\b = 
(\overline{q}, \overline{p} -\xi).$
Then there is a morphism $\rho:Z \to Y$ given by
$ \rho( \overline{p},\overline{q},\xi) = ((p,\xi-p),(q,-\xi-q)) $
such that the fibers of $\rho$ are precisely the $D_4$-orbits. Hence $Y$ is
birational to $Z/D_4$.

%{\bf 8. Fat points annihilated by $c=\prod_{\omega,k} \Omega(\omega+k\tau)$.}

%\bigskip \bigskip \bigskip \bigskip

\end{document}